\documentstyle[12pt]{amsart}
\let\text=\mbox
\newcommand{\skipthistext}[1]{}

\renewcommand{\pf}{\noindent {\bf Proof} \hspace{2mm}}

              \newcommand{\W}{{\cal W}}
              \newcommand{\U}{{\cal U}}
              \newcommand{\T}{{\cal T}}
              \newcommand{\V}{{\cal V}}

              \renewcommand{\O}{{\cal O}}
              \newcommand{\N}{{\cal N}}

\newcommand{\abs}[1]{\centerline{\sc Abstract}
\vspace{3mm}
\centerline{\parbox{150mm}{\small #1}}}
\newcounter{thm}

\def\thmnumber{
\addtocounter{thm}{1}\if\arabic{section}0\relax
\else\arabic{section}.\fi%
\if\arabic{subsection}0\relax
\else\arabic{subsection}.\fi%
\arabic{thm}
\hspace{2mm}}

\def\setthmnumber{\setcounter{thm}{0}}

\newenvironment{thm}{
\vspace{3mm}\par\noindent
{\bf Theorem \thmnumber} \it}
{\vspace{3mm}\par\rm}

\newenvironment{defi}{
\vspace{3mm}\par\noindent
{\bf Definition \thmnumber}\it}
{\vspace{3mm}\par\rm}

\newenvironment{lem}{
\vspace{3mm}\par\noindent
{\bf Lemma \thmnumber} \it}
{\vspace{3mm}\par\rm}

\newenvironment{prop}{
\vspace{3mm}\par\noindent
{\bf Proposition \thmnumber}\it}
{\vspace{3mm}\par\rm}

\newenvironment{coro}{
\vspace{3mm}\par\noindent
{\bf Corollary \thmnumber}\it}
{\vspace{3mm}\par\rm}

\newenvironment{re}{
\vspace{3mm}\par\noindent
{\bf Remark \thmnumber}\it}
{\vspace{3mm}\par\rm}

\newcommand{\sectioni}[1]{\setthmnumber\setcounter{subsection}{0}
\section{#1}}

\renewcommand{\subsection}[1]{%
\vspace{7mm}\par\noindent%
\addtocounter{subsection}{1}%
{\sc\arabic{section}.\arabic{subsection} \ {#1}}
\vspace{4mm}\par\noindent\setthmnumber}

\begin{document}
\title{On a notion of maps between orbifolds\\
II. homotopy and CW-complex}
\author{Weimin Chen}
\date{\today\\
\hspace{2mm}{\it Keywords}. Orbifold. Homotopy theory. 
}
\maketitle

\abs{This is the second of a series of papers which are devoted to a
comprehensive theory of maps between orbifolds. In this 
paper, we develop a basic machinery for studying homotopy classes 
of such maps. It contains two parts: (1) the construction of a set of 
algebraic invariants -- the homotopy groups, and (2) an analog of 
CW-complex theory. As a corollary of this machinery, the classical 
Whitehead theorem which asserts that a weak homotopy equivalence is a 
homotopy equivalence is extended to the orbifold category.}

\sectioni{Introduction}

In \cite{C1} we introduced a notion of maps between orbifolds, and 
established several basic results concerning the topological structure 
of the corresponding mapping spaces. This was in an attempt to initiate 
a comprehensive study of orbifolds, which was motivated by work of Dixon, 
Harvey, Vafa and Witten \cite{DHVW} on string theories of orbifolds. 
Particularly, the aforementioned theorems about the topological structure 
of 
mapping spaces had direct applications in the theory of pseudoholomorphic
curves and Gromov-Witten invariants of symplectic orbifolds
(cf. \cite{CR1, CR2, C2, C3, C4}). See the introduction of \cite{C1} for
more detailed explanations on this aspect of the story. 

The present paper is continuation of \cite{C1}. The main objective of this 
paper is to establish the corresponding homotopy theory for studying such maps
between orbifolds. In particular, we have developed techniques of exhibiting 
a class of orbispaces (which include orbifolds) as results of gluing together 
some fundamental building blocks, ie. cells of various isotropy types. The
corresponding algebraic invariants needed for describing such constructions 
are the sets of homotopy classes of the gluing maps. Thus as homotopy groups
in this theory, we have studied in details the sets of homotopy classes of 
maps from spheres of various isotropy types into an orbispace. We remark
that these homotopy groups form a strictly larger set of invariants than the 
usual homotopy groups of an orbifold (which by definition are the homotopy 
groups of the corresponding classifying spaces, cf. \cite{Ha2}). See 
Proposition 1.4 for more details in this regard.

In this paper, we shall mainly work with a more restricted class
of maps rather than the general ones introduced in \cite{C1}. The details are
presented in the following assumption throughout. We refer the
reader to \S 2.2 of \cite{C1} for the basic definitions.

\vspace{2mm}

\noindent{\bf Convention}\hspace{1.5mm}
In this paper, we shall only consider maps of orbispaces which are
equivalence classes of groupoid homomorphisms $(\{f_\alpha\},
\{\rho_{\beta\alpha}\}):\Gamma\{U_\alpha\}\rightarrow\Gamma
\{U_{\alpha^\prime}^\prime\}$, where each $\rho_\alpha=\rho_{\alpha\alpha}:
G_{U_\alpha}\rightarrow G_{U_\alpha^\prime}$ is an injective homomorphism.
As a result of this assumption, the mappings $\{\rho_{\beta\alpha}\}$
are partially injective in the sense that if $\rho_{\beta\alpha}
(\xi_1)=\rho_{\beta\alpha}(\xi_2)$ and $\mbox{Domain }(\phi_{\xi_1})
\cap\mbox{Domain }(\phi_{\xi_2})\neq\emptyset$, then $\xi_1=\xi_2$.

\vspace{2mm}

In studying homotopy classes of maps and homotopy types of orbispaces,
we need to fix a base point in the orbispaces throughout the consideration,
in which all maps have to preserve the pre-chosen base point structures.
We shall give an introduction to this notion next, before we come to the
detailed description of the homotopy classes of maps and homotopy types
of orbispaces that are to be considered in this paper.

Let $X$ be an orbispace, with its atlas of local charts denoted by
$\U=\{U_i\}$. A base-point structure of $X$ is a triple, denoted
by $\underline{o}=(o,U_o,\hat{o})$, where $o\in X$, $U_o\in\U$ such that
$o\in U_o$, and $\hat{o}\in\pi_{U_o}^{-1}(o)\subset\widehat{U_o}$.
An orbispace $X$ with a chosen base-point structure $\underline{o}$ will
be denoted by $(X,\underline{o})$. We remark that for the special
case where $X=Y/G$, with the action of $G$ on $Y$ being trivial,
we have $X=Y$ and furthermore, a base-point structure $\underline{o}
=(o,U_o,\hat{o})$ is completely determined by the base point $o$. We
shall fix the notation $Y(G)$ for such an $X=Y/G$, and write $(Y(G),o)$ for
$(X,\underline{o})$. The orbispace $Y(G)$ will be called a space of
isotropy type $G$ throughout. The most frequently used examples of this 
type of orbispaces are $S^k(G), D^k(G)$, ie. the k-sphere and k-cell
of isotropy type $G$. Note that for the 1-cell and 0-cell of isotropy type 
$G$, we shall use different notations $I(G)$, $B_G$ respectively.

With these notations being fixed throughout the paper, we now consider the
maps from $(X,\underline{o})$ to $(X^\prime,\underline{o^\prime})$, where
$\underline{o}=(o,U_o,\hat{o})$ and $\underline{o^\prime}=(o^\prime,
U_o^\prime,\hat{o^\prime})$, which preserve the corresponding base-point
structures. Let $\rho:G_{\hat{o}}\rightarrow G_{\hat{o^\prime}}$ be any given
injective homomorphism, where $G_{\hat{o}}$, $G_{\hat{o^\prime}}$
are the subgroups of $G_{U_o}$, $G_{U_o^\prime}$ that fix $\hat{o}$
and $\hat{o^\prime}$ respectively, and let $\sigma=(\{f_\alpha\},
\{\rho_{\beta\alpha}\}):\Gamma\{U_\alpha\}\rightarrow
\Gamma\{U_{\alpha^\prime}^\prime\}$ be any groupoid homomorphism which
satisfies the following conditions: (1) $U_o\in\{U_\alpha\}$, $U_o^\prime\in
\{U_{\alpha^\prime}^\prime\}$, and $U_o\mapsto U_o^\prime$ under the
correspondence $U_\alpha\mapsto U_\alpha^\prime$, (2) $f_o:\widehat{U_o}
\rightarrow\widehat{U_o^\prime}$ satisfies $f_o(\hat{o})=\hat{o^\prime}$,
and (3) $\rho_o|_{G_{\hat{o}}}=\rho$ . Suppose $\tau=(\{f_a\},\{\rho_{ba}\})$
is induced by $\sigma$ via $\bar{\gamma}=(\theta,\{\xi_a\},\{\xi_a^\prime\})$
(cf. Lemma 2.2.4 in \cite{C1}), where Condition (1) above is satisfied for
$\tau$, and $\bar{\gamma}$ satisfies $\theta(o)=o$, $\xi_o=1\in T(U_o,U_o)
=G_{U_o}$, and $\xi_o^\prime=1\in T(U_o^\prime,U_o^\prime)=G_{U_o^\prime}$,
then one can easily verify that Conditions (2), (3) above are also satisfied
for $\tau$. In other words, the base-point structures $\underline{o}$
and $\underline{o^\prime}$ are preserved under the process of taking
equivalence classes in the sense described above. We will call such
an equivalence class a map from $(X,\underline{o})$ to
$(X^\prime,\underline{o^\prime})$, and we denote the set of all such
maps by $[(X,\underline{o});(X^\prime,\underline{o^\prime})]_\rho$.

Some variant or generalization of the preceding will also be considered.
For example, suppose $(X,\underline{o})$ and $(X^\prime,
\underline{o^\prime})$ are orbispaces with pre-chosen base-point structures,
where in the base-point structure $\underline{o}=(o,U_o,\hat{o})$, the
group $G_{U_o}$ acts trivially on $\widehat{U_o}$. In this case, we shall
define maps from $(X,\underline{o})$ to $(X^\prime,\underline{o^\prime})$
by further allowing, in the definition of equivalence relation, that $U_o$
may be replaced by a $V_o\subset U_o$ where $o\in V_o$, for which the
requirements $\xi_o=1\in G_{U_o}=G_{V_o}$ and $\xi_o^\prime=1\in
G_{U_o^\prime}$ in $(\theta,\{\xi_a\},\{\xi_a^\prime\})$ still make sense.
Correspondingly, this means that the base-point structure
$\underline{o}=(o,U_o,\hat{o})$ is reduced to the ordinary base point
$o$. We shall simply write $o$ for $\underline{o}$. An important case for
such an $(X,\underline{o})$ is $(Y(G),o)$.

On the other hand, given any $X$ and $X^\prime$, we may fix a set of
base-point structures $\{\underline{o}_i\}$ of $X$ and a set of
base-point structures $\{\underline{o^\prime}_k\}$ of $X^\prime$,
together with a correspondence $\underline{o}_i\mapsto
\underline{o^\prime}_k=\underline{o^\prime}_{\theta(i)}$.
In the same vein, we may define maps between these ``multi-based''
orbispaces, requiring that the given sets of base-point structures are
preserved under these maps with respect to the given correspondence
$\underline{o}_i\mapsto\underline{o^\prime}_k
=\underline{o^\prime}_{\theta(i)}$.

\vspace{3mm}

Now we give some detailed descriptions about the homotopy classes of maps
and homotopy types of orbispaces that will be considered in this paper.

Two maps $\Phi_1,\Phi_2:X\rightarrow X^\prime$ are said
to be homotopic (and written $\Phi_1\simeq\Phi_2$) if there exists a map
$\Psi:X\times[a,b]\rightarrow X^\prime$ such that $\Phi_1,\Phi_2$
are the restrictions of $\Psi$ to the subspaces $X\times \{a\}$
and $X\times \{b\}$ respectively. The map $\Psi$ is called a
homotopy between $\Phi_1$ and $\Phi_2$. Let $\Psi_1:X\times [a,b]
\rightarrow X^\prime$ be a homotopy between $\Phi_1$ and $\Phi_2$,
and $\Psi_2:X\times [b,c]\rightarrow X^\prime$ be a homotopy between
$\Phi_2$ and $\Phi_3$. Then there is a homotopy $\Psi_3:X\times [a,c]
\rightarrow X^\prime$ between $\Phi_1$ and $\Phi_3$, whose restrictions
to $X\times [a,b]$ and $X\times [b,c]$ are $\Psi_1$ and $\Psi_2$
respectively. A homomorphism representing $\Psi_3$ may be obtained
as follows. Take a representing homomorphism $\sigma_1$ of
$\Phi_1$, and a representing homomorphism $\sigma_2$ of $\Phi_2$.
We may assume, by passing to an induced one, that near $X\times \{b\}$,
$\sigma_1$ is given by $(\{F_{i,1}\},\{\hat{\rho}_{ji,1}\}):
\Gamma\{U_i\times (b_{-}^i,b]\}\rightarrow\Gamma\{U_{i^\prime}^\prime\}$
and $\sigma_2$ is given by $(\{F_{i,2}\},\{\hat{\rho}_{ji,2}\}):
\Gamma\{U_i\times [b,b_{+}^i)\}\rightarrow\Gamma\{U_{i^\prime}^\prime\}$
such that for each index $i$, $F_{i,1}|_{\widehat{U_i}\times \{b\}}
=F_{i,2}|_{\widehat{U_i}\times \{b\}}$, and for each pair of indexes
$(i,j)$, $\hat{\rho}_{ji,1}=\hat{\rho}_{ji,2}$. (Note that
$T(U_i\times (b_{-}^i,b],U_j\times (b_{-}^j,b])=T(U_i,U_j)
=T(U_i\times [b,b_{+}^i),U_j\times [b,b_{+}^j))$.) It is clear that
such a pair $(\sigma_1,\sigma_2)$
of homomorphisms can be patched together to form a homomorphism
whose equivalence class is a map from $X\times [a,c]$ to $X^\prime$.
We define $\Psi_3$ to be the corresponding map. Note that although
the homotopy $\Psi_3$ may not be uniquely determined, its
existence implies that the homotopy relation is an equivalence
relation on the set of maps $[X;X^\prime]$. The corresponding set
of homotopy classes will be denoted by $[[X;X^\prime]]$.

Let $\Phi:Y\rightarrow X$ be any map. There are induced mappings
$\Phi^\ast:[[X;X^\prime]]\rightarrow [[Y;X^\prime]]$ sending
$[\Psi]$ to $[\Psi\circ\Phi]$, and $\Phi_\ast:[[X^\prime;Y]]
\rightarrow [[X^\prime;X]]$ sending $[\Psi]$ to $[\Phi\circ\Psi]$.
The mappings $\Phi^\ast, \Phi_\ast$ depend only on the homotopy class of
$\Phi$. A map $\Phi:Y\rightarrow X$ is called a homotopy equivalence if
there is a map $\Psi:X\rightarrow Y$ such that $\Phi\circ\Psi\simeq Id_X$
and $\Psi\circ\Phi\simeq Id_Y$. It is routine to check that
$\Phi:Y\rightarrow X$ is a homotopy equivalence if and only if both
mappings $\Phi^\ast, \Phi_\ast$ are bijections for any orbispace $X^\prime$.

Homotopy preserving given base-point structures can be defined in the
same vein. Let $\underline{o}=(o,U_o,\hat{o})$, $\underline{o}^\prime
=(o^\prime,U_o^\prime,\hat{o}^\prime)$ be any base-point structures on
$X$ and $X^\prime$ respectively. Two maps $\Phi_1,\Phi_2\in
[(X,\underline{o});(X^\prime,\underline{o}^\prime)]_\rho$ are said to be
homotopic (and written $\Phi_1\simeq\Phi_2$) if there exists a homomorphism
$\sigma=(\{F_i\},\{\hat{\rho}_{ji}\}):\Gamma\{U_i\times I_i\}\rightarrow
\Gamma\{U_{i^\prime}^\prime\}$, where each $I_i$ is a sub-interval of
$[a,b]$, such that (1) $U_o\times [a,b]\in\{U_i\times I_i\}$, which
corresponds to $U_o^\prime$ under $U_i\times I_i\mapsto U_i^\prime$,
(2) $F_o(\{\hat{o}\}\times [a,b])=\hat{o}^\prime$, (3)
$\hat{\rho}_{o}|_{G_{\hat{o}}}=\rho$, and (4) $\Phi_1,\Phi_2$ are
the equivalence classes of the restriction of $\sigma$ to $X\times \{a\}$
and $X\times \{b\}$ respectively. The set of homotopy classes of elements
in $[(X,\underline{o});(X^\prime,\underline{o}^\prime)]_\rho$ will be
denoted by $[[(X,\underline{o});(X^\prime,\underline{o}^\prime)]]_\rho$.
Given any $\Phi\in [(Y,\underline{p});(X,\underline{o})]_\rho$,
there are induced mappings $\Phi^\ast:[[(X,\underline{o});
(X^\prime,\underline{o}^\prime)]]_\eta\rightarrow [[(Y,\underline{p});
(X^\prime,\underline{o}^\prime)]]_{\eta\circ\rho}$ and $\Phi_\ast:
[[(X^\prime,\underline{o}^\prime);(Y,\underline{p})]]_{\eta}\rightarrow
[[(X^\prime,\underline{o}^\prime);(X,\underline{o})]]_{\rho\circ\eta}$,
defined by $\Phi^\ast([\Psi])=[\Psi\circ\Phi]$ and $\Phi_\ast([\Psi])
=[\Phi\circ\Psi]$, which depend only on the homotopy class of $\Phi$.

Homotopy classes of maps between pairs of orbispaces $(X,A)$ will also
be considered, where $A$ is a subspace of $X$. If
$\underline{o}=(o,U_o,\hat{o})$ is a base-point structure
of $X$ such that $o\in A$, then there is a canonical base-point
structure of $A$, denoted by $\underline{o}|_A=(o,V_o,\hat{o})$,
where $V_o$ is the connected component of $U_o\cap A$ that contains
$o$. We will denote a pair with such a base-point structure by $(X,A,
\underline{o})$. A map from $(X,A)$ to $(X^\prime,A^\prime)$ is
an element $\Phi\in [X;X^\prime]$ such that $\Phi|_A\in
[A;A^\prime]$, and a map from $(X,A,\underline{o})$ to $(X^\prime,
A^\prime,\underline{o^\prime})$ is an element $\Phi\in
[(X,\underline{o});(X^\prime,\underline{o^\prime})]_\rho$
for some $\rho:G_{\hat{o}}\rightarrow G_{\hat{o^\prime}}$, such that
$\Phi|_A\in [(A,\underline{o}|_A);(A^\prime,\underline{o^\prime}
|_{A^\prime})]_\rho$.
The set of homotopy classes of the latter will be denoted by
$[[(X,A,\underline{o});(X^\prime,A^\prime,\underline{o^\prime})]]_\rho$.
Given any $\Phi\in [(Y,B,\underline{p});(X,A,\underline{o})]_\rho$,
there are induced mappings $\Phi^\ast:[[(X,A,\underline{o});(X^\prime,
A^\prime,\underline{o^\prime})]]_\eta\rightarrow [[(Y,B,\underline{p});
(X^\prime,A^\prime,\underline{o^\prime})]]_{\eta\circ\rho}$ and
$\Phi_\ast:[[(X^\prime,A^\prime,\underline{o^\prime});
(Y,B,\underline{p})]]_\eta\rightarrow [[(X^\prime,A^\prime,
\underline{o^\prime});(X,A,\underline{o})]]_{\rho\circ\eta}$ defined
by $\Phi^\ast([\Psi])=[\Psi\circ\Phi]$, $\Phi_\ast([\Psi])=[\Phi\circ\Psi]$.
Again the mappings $\Phi^\ast,\Phi_\ast$ depend only on the
homotopy class of $\Phi$.

\vspace{3mm}

With the preceding preparatory discussion, we now present the
main results in this paper.

\begin{defi}
For any base-point structure $\underline{o}=(o,U_o,\hat{o})$ and
injective homomorphism $\rho:G\rightarrow G_{\hat{o}}$, we define
$$
\pi_k^{(G,\rho)}(X,\underline{o})=[[(S^k(G),\ast);(X,\underline{o})]]_\rho,
\;\forall k\geq 0,
$$
and
$$
\pi_{k+1}^{(G,\rho)}(X,A,\underline{o})=[[(D^{k+1}(G),S^k(G),\ast);(X,A,
\underline{o})]]_\rho, \;\forall k\geq 0.
$$
For the special case when $G=\{1\}$, we simply write
$\pi_k(X,\underline{o})$ and $\pi_k(X,A,\underline{o})$ instead.
\end{defi}

Regarding the structure of these homotopy sets, we have

\begin{thm}
\begin{itemize}
\item [{(1)}] {\em Algebraic structures:} $\pi_k^{(G,\rho)}(X,\underline{o})$,
$\pi_{k+1}^{(G,\rho)}(X,A,\underline{o})$ are canonically identified with the
$\pi_{k-1}$ of $[(S^1(G),\ast);(X,\underline{o})]_\rho$ and the $\pi_k$ of
$[(I(G),S^0(G),\ast);(X,A,\underline{o})]_\rho$ respectively, hence have
natural group structures for $k\geq 1$ which are Abelian when $k\geq 2$.
\item [{(2)}] {\em Functoriality:} For any $\Phi\in [(X,\underline{o});
(X^\prime,\underline{o^\prime})]_\eta$ {\em (resp. $\Phi\in
[(X,A,\underline{o});(X^\prime,A^\prime,\underline{o^\prime})]_\eta$)},
there are natural homomorphisms
$\Phi_\ast:\pi_k^{(G,\rho)}(X,\underline{o})\rightarrow
\pi_k^{(G,\rho^\prime)}(X^\prime,\underline{o^\prime})$ {\em (resp.
$\Phi_\ast:\pi_k^{(G,\rho)}(X,A,\underline{o})\rightarrow
\pi_k^{(G,\rho^\prime)}(X^\prime,A^\prime,\underline{o^\prime})$)}
with $\rho^\prime=\eta\circ\rho$, which depend only on the
homotopy class of $\Phi$. For any $(H,\eta)$ where $\eta:H\rightarrow
G_{\hat{o}}$ is an injective homomorphism which factors through $\rho:
G\rightarrow G_{\hat{o}}$ by $\iota:H\rightarrow G$, there are
natural homomorphisms $\iota^\ast:\pi_k^{(G,\rho)}(X,\underline{o})
\rightarrow\pi_k^{(H,\eta)}(X,\underline{o})$ and $\iota^\ast:
\pi_{k+1}^{(G,\rho)}(X,A,\underline{o})\rightarrow\pi_{k+1}^{(H,\eta)}
(X,A,\underline{o})$.
\item [{(3)}] {\em Exact sequence:} There exist natural homomorphisms
$$
\partial:\pi_{k+1}^{(G,\rho)}(X,A,\underline{o})\rightarrow
\pi_k^{(G,\rho)}(A,\underline{o}|_A),\;\forall k\geq 0
$$
and a long exact sequence
$$
\begin{array}{l}
\cdots\rightarrow \pi_{k+1}^{(G,\rho)}(X,A,\underline{o})
\stackrel{\partial}{\rightarrow}\pi_{k}^{(G,\rho)}(A,\underline{o}|_A)
\stackrel{{i}_{\ast}}{\rightarrow}\pi_k^{(G,\rho)}(X,\underline{o})
\stackrel{{j}_{\ast}}{\rightarrow}\pi_k^{(G,\rho)}(X,A,\underline{o})
\stackrel{\partial}{\rightarrow}\cdots \\
\stackrel{{j}_{\ast}}{\rightarrow}\pi_1^{(G,\rho)}(X,A,\underline{o})
\stackrel{\partial}{\rightarrow}\pi_0^{(G,\rho)}(A,\underline{o}|_A)
\stackrel{{i}_{\ast}}{\rightarrow}\pi_0^{(G,\rho)}(X,\underline{o}).
\end{array}
$$
\end{itemize}
\end{thm}

There are natural homomorphisms
$$
C:\pi_1^{(G,\rho)}(X,\underline{o})\rightarrow \mbox{Aut}
(\pi_k^{(G,\rho)}(X,\underline{o})),\; \forall k\geq 1
$$
and
$$
C^\prime:\pi_1^{(G,\rho)}(A,\underline{o}|_A)\rightarrow \mbox{Aut}
(\pi_{k+1}^{(G,\rho)}(X,A,\underline{o})),\;\forall k\geq 1
$$
which satisfy $C(z)\circ\partial=\partial\circ C^\prime(z)$,
$\forall z\in \pi_1^{(G,\rho)}(A,\underline{o}|_A)$.

Regarding the dependence of $\pi_k^{(G,\rho)}(X,\underline{o})$ and
$\pi_{k+1}^{(G,\rho)}(X,A,\underline{o})$ on the base-point structure
$\underline{o}$ and the data $(G,\rho)$, we have

\begin{prop}
\begin{itemize}
\item [{(1)}] Given any path $u\in [(I(G),0,1);(X,\underline{o_1},
\underline{o_2})]_{(\rho,\eta)}$, there is an isomorphism, written
$u_\ast:\pi_k^{(G,\eta)}(X,\underline{o_2})\rightarrow
\pi_k^{(G,\rho)}(X,\underline{o_1})$, whose inverse $u_\ast^{-1}$ is the
isomorphism associated to the inverse path $\nu(u)\in [(I(G),0,1);
(X,\underline{o_2},\underline{o_1})]_{(\eta,\rho)}$. In fact, for any
two paths $u_1,u_2\in [(I(G),0,1);(X,\underline{o_1},\underline{o_2})]
_{(\rho,\eta)}$, $\nu(u_2)_\ast\circ (u_1)_\ast=C([\nu(u_2)\# u_1])\in
\mbox{Aut}(\pi_k^{(G,\eta)}(X,\underline{o_2}))$.
\item [{(2)}] Given any path $u\in [(I(G),0,1);(A,\underline{o_1}|_A,
\underline{o_2}|_A)]_{(\rho,\eta)}$, there is an associated mapping
$u_\ast:\pi_{k+1}^{(G,\eta)}(X,A,\underline{o_2})\rightarrow
\pi_{k+1}^{(G,\rho)}(X,A,\underline{o_1})$, which is an isomorphism for
$k\geq 1$ and a base point preserving bijection when $k=0$. Moreover,
for any $u_1,u_2\in [(I(G),0,1);(A,\underline{o_1}|_A,
\underline{o_2}|_A)]_{(\rho,\eta)}$, $\nu(u_2)_\ast\circ (u_1)_\ast
=C^\prime([\nu(u_2)\# u_1])\in\mbox{Aut}(\pi_{k+1}^{(G,\eta)}
(X,A,\underline{o_2}))$.
\item [{(3)}] The isomorphism class of $\pi_k^{(G,\rho)}(X,\underline{o})$
or $\pi_{k+1}^{(G,\rho)}(X,A,\underline{o})$ depends only on the conjugacy
class of the subgroup $\rho(G)\subset G_{\hat{o}}$.
\end{itemize}
\end{prop}

Regarding the nature of $\pi_k^{(G,\rho)}(X,\underline{o})$ and
$\pi_{k+1}^{(G,\rho)}(X,A,\underline{o})$ in certain special cases, we
have

\begin{prop}
\begin{itemize}
\item [{(1)}] $\pi_k(X,\underline{o})$, $\pi_{k+1}(X,A,\underline{o})$
are naturally isomorphic to $\pi_k(B\Gamma_X,\ast)$ and $\pi_{k+1}(B\Gamma_X,
B\Gamma_A,\ast)$, where $B\Gamma_X, B\Gamma_A$ are the classifying spaces
of the defining groupoid for $X$ and $A$ respectively, cf. Haefliger
\cite{Ha1,Ha2}.
\item [{(2)}] When $X=Y/G$ is a global quotient, for any subgroup $\rho:
H\subset G$, there are natural isomorphisms $\pi_k^{(H,\rho)}(X,\underline{o})
\cong \pi_k(Y^H,\hat{o})$ for all $k\geq 2$, and the natural exact sequence
$$
1\rightarrow \pi_1(Y^H,\hat{o})\rightarrow\pi_1^{(H,\rho)}(X,\underline{o})
\rightarrow C(H)\rightarrow \pi_0(Y^H,\hat{o})\rightarrow
\pi_0^{(H,\rho)}(X,\underline{o}),
$$
where $Y^H$ is the fixed-point set of $H$ and $C(H)$ is the centralizer of
$H$.
\end{itemize}
\end{prop}

We remark that the homotopy groups $\{\pi_k(Y^H,\hat{o})\}$ of the various
fixed-point sets appeared in the coefficient systems for the equivariant
obstruction theory in Bredon \cite{Bre}. In fact, the basic homotopy
machinery developed in this paper laid the necessary foundation for
extending the equivariant obstruction theory in Bredon \cite{Bre}
to the orbispace category.

This paper also contains a detailed, elementary treatment of the covering
theory for orbispaces. Although the theory of coverings for \'{e}tale
topological groupoids is well-understood in principle (cf. \cite{BH},
\cite{Ha3}), there is a number of results which will be used in the
later constructions in this paper that are in a specific form and can not
be found in the literature.

\vspace{2mm}

This roughly constitutes the first half of the paper, which is
concerned with the basic properties of the algebraic invariants
--- the homotopy groups --- for this homotopy theory. In the remaining
second half, we develop an analog of the CW-complex theory for
the orbispace category.

A fundamental construction in this part is given in the following

\begin{prop}
Let $Y$ be a locally path-connected and semi-locally 1-connected
orbispace {\em (cf. \S 2.4)}. For any map $\Phi:Y\rightarrow X$, there
is a canonical orbispace structure on the mapping cylinder
$M_\phi$, where $\phi:Y\rightarrow X$ is the induced map between
underlying spaces, such that {\em (1)} there are natural embeddings
$i:Y\rightarrow M_\phi$, $j:X\rightarrow M_\phi$, realizing the
orbispaces $Y$, $X$ as a subspace of $M_\phi$, and {\em (2)}
$j:X\rightarrow M_\phi$ is a strong deformation retract by a canonical
retraction $r:M_\phi\rightarrow X$ which satisfies $\Phi=r\circ i$.
\end{prop}

The above `mapping cylinder' construction is then applied to the
special case where $Y=S^{k-1}(G)$. Under further conditions on both
the map $\Phi$ and the orbispace $X$, we can collapse the subspace
$S^{k-1}(G)$ in the mapping cylinder to a point. This procedure gives
rise to a canonical orbispace structure on the mapping cone of the
induced map $\phi:S^{k-1}\rightarrow X$, which is what we will refer 
to as `attaching a k-cell of isotropy type $G$ to $X$ via $\Phi$'.

For any $n\geq 0$, denote by ${\cal{C}}^n$ the set of orbispaces
$X$, where there is a canonical filtration of subspaces
$$
X_0\subset X_1\subset\cdots\subset X_n=X
$$
such that $X_0$ is a finite subset and for each $1\leq k\leq n$,
$X_k$ is resulted from attaching finitely many k-cells of various
isotropy types to $X_{k-1}$. We let ${\cal{C}}$ be the union of
${\cal{C}}^n$ for all $n\geq 0$.

It turns out that this sub-category of orbispaces ${\cal{C}}$ is
the right one for the main objective of this paper, ie., developing a
machinery for studying
homotopy classes of maps and homotopy types of orbispaces. In this
regard, we have the following

\begin{thm}
\begin{itemize}
\item [{(1)}] Suppose $\Phi:X\rightarrow X^\prime$ is a weak
homotopy equivalence and the mapping cylinder of $\Phi$ is defined.
Then for any $Y\in {\cal{C}}$, the mapping $\Phi_\ast:
[[Y;X]]\rightarrow [[Y;X^\prime]]$ is a bijection.
\item [{(2)}] For any $X,X^\prime\in {\cal{C}}$, a map $\Phi:X\rightarrow
X^\prime$ is a homotopy equivalence if and only if it is a weak
homotopy equivalence.
\end{itemize}
\end{thm}

Set $\pi_0^G(X)\equiv [[B_G;X]]$. Then in the preceding theorem, a map 
$\Phi:X\rightarrow X^\prime$ is
called a weak homotopy equivalence if $\Phi_\ast:\pi_0^G(X)\rightarrow
\pi_0^G(X^\prime)$ is a bijection for all $G$, and $\Phi_\ast:
\pi_k^{(G,\rho)}(X,\underline{o})\rightarrow\pi_k^{(G,\rho^\prime)}
(X^\prime,\underline{o^\prime})$ is isomorphic for $k\geq 0$ for all
possible data $\underline{o}, \underline{o^\prime},(G,\rho)$ and
$(G,\rho^\prime)$.

Although the sub-category ${\cal{C}}$ is favorable for the purpose of
homotopy theory, it is not clear a priori that ${\cal{C}}$ will contain
any geometrically interesting examples, such as compact smooth orbifolds.
Next we describe a construction which will justify the sub-category
${\cal{C}}$ in this regard.

Let $K$ be a finite CW-complex such that for any attaching map, if
its image meets the interior of a cell, then it contains the whole
cell. An arrow of $K$ is an ordered pair $(\sigma_1,\sigma_2)$ of
cells in $K$ where $\sigma_2$ is a face of $\sigma_1$. We denote
arrows by lower case letters $a,b,c$, etc. Let $a=(\sigma_1,\sigma_2)$
be an arrow. Then $\sigma_1$ is called the initial point of $a$,
denoted by $i(a)$, and $\sigma_2$ is called the terminal point of
$a$, denoted by $t(a)$. Two arrows $a,b$ are composable if
$t(b)=i(a)$, and in this case, the composition of $a,b$, denoted
by $ab$, is the arrow $c$ uniquely determined by the conditions
$i(c)=i(b)$, $t(c)=t(a)$. Composition of arrows is clearly
associative.

With the preceding understood, a CW-complex of groups on $K$,
denoted by $(K,G_\sigma,\psi_a,g_{a,b})$, is given by the following
set of data:

\begin{itemize}
\item [{(1)}] Each cell $\sigma$ in $K$ is associated with a group
$G_\sigma$.
\item [{(2)}] Each arrow $a$ of $K$ is assigned with an injective
homomorphism $\psi_a:G_{i(a)}\rightarrow G_{t(a)}$.
\item [{(3)}] Each pair of composable arrows $a,b$ is assigned
with an element $g_{a,b}\in G_{t(a)}$, such that
$$
Ad(g_{a,b})\circ\psi_{ab}=\psi_a\circ\psi_b.
$$
Moreover, the following cocycle condition holds for any triple
$a,b,c$ of composable arrows
$$
\psi_a(g_{b,c})g_{a,bc}=g_{a,b}g_{ab,c}.
$$
\end{itemize}

Two CW-complexes of groups $(K,G_\sigma,\psi_a,g_{a,b})$,
$(K,G_\sigma,\psi_a^\prime,g_{a,b}^\prime)$ are said to be
equivalent if there
is a set $\{g_a\in G_{t(a)}\}$ such that
$$
\psi_a^\prime=Ad(g_a)\circ\psi_a,\hspace{2mm}
g_{a,b}^\prime=g_a\psi_a(g_b)g_{a,b}g_{ab}^{-1}.
$$

Finally, let $K_n\subset K$, $n\geq 0$, be the n-skeleton of $K$.
Then any CW-complex of groups on $K$ naturally induces one on $K_n$
by restriction.

We remark that the notion `CW-complex of groups' is a natural
extension of the notion `complex of groups' in Haefliger
\cite{Ha3}. Moreover, the following result generalizes the one to
one correspondence between equivalence classes of complexes of
groups and isomorphism classes of orbihedra therein.

\begin{prop}
Let $K$ be a finite CW-complex and $X$ be its underlying space. To
each equivalence class of CW-complexes of groups on $K$, there is
associated an orbispace structure on $X$, called the geometric
realization\footnote{there is a minor difference, although conceptual
in nature, between the notions `geometric realization' and
`orbihedron', see \S 3.3 for more detailed comments.}, which gives
a one to one correspondence, such that the orbispace $X\in
{\cal{C}}$. Moreover, the geometric realizations of the
restrictions of the CW-complex of groups to the n-skeletons of $K$
provide the natural canonical filtration of subspaces for the
geometric realization of the CW-complex of groups itself.
\end{prop}

Finally, recall the following fact about compact smooth
orbifolds, cf. e.g. \cite{Pr,Y}:

\vspace{2mm}

{\em A compact smooth orbifold may be triangulated so that it
becomes the geometric realization of a natural simplicial complex
of groups on the resulting simplicial complex.}

\vspace{2mm}

As a corollary, the classical Whitehead theorem is extended to the
orbifold category.

\vspace{2.5mm}

\noindent{\bf Corollary}\hspace{3.5mm}{\em
A map between compact smooth orbifolds is a homotopy equivalence
if and only if it is a weak homotopy equivalence.
}
\vspace{2.5mm}

The following is a glimpse at the organization of this paper.
\begin{itemize}
\item [{\S 2.1}] Space of guided loops.
\item [{\S 2.2}] Generalities of homotopy groups.
\item [{\S 2.3}] Exact sequence of a pair.
\item [{\S 2.4}] The theory of coverings.
\item [{\S 3.1}] Construction of mapping cylinders.
\item [{\S 3.2}] Orbispaces via attaching cells of isotropy.
\item [{\S 3.3}] CW-complex of groups and its geometric realization.
\end{itemize}

\vspace{1cm}

\centerline{\bf Acknowledgments}

\vspace{3mm}

This research was partially supported by NSF grants DMS-9971454, 
DMS-9729992 (through IAS), and DMS-0304956 (during the final preparation).

\vspace{1cm}

\sectioni{Homotopy groups via guided loop spaces}

\subsection{Space of guided loops}

The homotopy sets defined in Definition 1.1 will be studied through the
space of `guided' paths or loops in the orbispace. This subsection
is devoted to a preliminary study of these path or loop spaces.

Recall that $I(G)$ denotes the 1-cell of isotropy type $G$,
namely, the orbispace defined by the trivial action of $G$ on the
interval $I=[0,1]$. Let $X$ be an orbispace, given with a pair
of base-point structures $\underline{o_1}=(o_1,U_{o_1},\hat{o_1})$,
$\underline{o_2}=(o_2,U_{o_2},\hat{o_2})$, and injective homomorphisms
$\rho:G\rightarrow G_{\hat{o_1}}$, $\eta:G\rightarrow G_{\hat{o_2}}$.
We consider the set of maps $[(I(G),0,1);(X,\underline{o_1},
\underline{o_2})]_{(\rho,\eta)}$ from $(I(G),0,1)$ to $(X,\underline{o_1},
\underline{o_2})$, which by definition are equivalence classes of groupoid
homomorphisms $\sigma=(\{f_i\},\{\rho_{ji}\}):\Gamma\{I_i\}\rightarrow
\Gamma\{U_{i^\prime}\}$, where (1) $\{I_i\mid i=0,1,
\cdots,n\}$ is a cover of $I$ by sub-intervals such that $0\in
I_0$, $1\in I_n$, $I_i\cap I_j\neq\emptyset$ iff $j=i$ or $i+1$,
(2) denote by $U_i\in\{U_{i^\prime}\}$ the local chart assigned to $I_i$,
then $U_0=U_{o_1}$, $U_n=U_{o_2}$, (3) $f_0(0)=\hat{o_1}$,
$f_n(1)=\hat{o_2}$, and (4) $\rho_0=\rho$ and $\rho_n=\eta$. Note
that each mapping $\rho_{(i+1)i}:T(I_i,I_{i+1})\rightarrow T(U_i,U_{i+1})$
is completely determined by the image of $1\in G=T(I_i,I_{i+1})$
in $T(U_i,U_{i+1})$, we shall conveniently regard $\rho_{(i+1)i}$
as an element of $T(U_i,U_{i+1})$. On the other hand, the
homomorphisms $\rho_i=\rho_{ii}:G\rightarrow G_{U_i}$ can be
recovered inductively by $\rho_{i+1}=\lambda_{\rho_{(i+1)i}}\circ\rho_i$
with $\rho_0=\rho$ (cf. $(2.1.3)$ in \cite{C1} for the definition of
$\lambda_{\rho_{(i+1)i}}$). We observe that the image of $f_i:I_i\rightarrow
\widehat{U_i}$ for each index $i$ is forced to lie in the fixed-point set
of $\rho_i(G)\subset G_{U_i}$. For this reason we call each element in
$[(I(G),0,1);(X,\underline{o_1},\underline{o_2})]_{(\rho,\eta)}$ a
$(G,\rho,\eta)$-guided path in $(X,\underline{o_1},\underline{o_2})$,
or simply a guided path.

A canonical topology is given to $[(I(G),0,1);(X,\underline{o_1},
\underline{o_2})]_{(\rho,\eta)}$ by the general construction in \S 3.2
of Part I of this series \cite{C1}, as the orbispace $I(G)$ is trivially
paracompact, locally compact and Hausdorff. For this purpose, we shall
only consider the representatives $\sigma=(\{f_i\},\{\rho_{ji}\}):
\Gamma\{I_i\}\rightarrow\Gamma\{U_{i^\prime}\}$ which are admissible.
In this case, it simply means that each $f_i:I_i\rightarrow\widehat{U_i}$
is extended over to the closure $\overline{I_i}$ of $I_i$. We recall
that
$$
\O_{\{\rho_{ji}\}}=\{\sigma\mid \sigma=(\{f_i\},\{\rho_{ji}\}) \mbox{
is admissible}, [\sigma]\in [(I(G),0,1);(X,\underline{o_1},
\underline{o_2})]_{(\rho,\eta)}\}. \leqno (2.1.1)
$$
By Lemma 3.1.2 in \cite{C1}, each $\O_{\{\rho_{ji}\}}$ is embedded into
$[(I(G),0,1);(X,\underline{o_1},\underline{o_2})]_{(\rho,\eta)}$ via
$\sigma\mapsto [\sigma]$. We consider the subsets
$$
\O_{\{\rho_{ji}\}}(\{\underline{K_i}\},\{\underline{O_i}\})=\{\{f_i\}\in
\O_{\{\rho_{ji}\}}\mid f_i(K_{i,s})\subset O_{i,s},s\in\Lambda(i)\}
\leqno (2.1.2)
$$
where $\underline{K_i}=\{K_{i,s}\mid s\in\Lambda(i)\}$ is a finite set
of compact subsets of $I_i$, and $\underline{O_i}=\{O_{i,s}\mid s\in
\Lambda(i)\}$ is a finite set of open subsets of $\widehat{U_i}$.
The canonical topology on $[(I(G),0,1);(X,\underline{o_1},
\underline{o_2})]_{(\rho,\eta)}$ is the one generated by the subsets
in $(2.1.2)$ for all possible data $\{\rho_{ji}\}$, $\{\underline{K_i}\}$
and $\{\underline{O_i}\}$. By Lemma 3.2.2 in \cite{C1}, if $\{f_i\}\in
\O_{\{\rho_{ji}\}}$ is equivalent to $\{f_k\}\in\O_{\{\rho_{lk}\}}$, then
there is a local homeomorphism $\phi$ from an open neighborhood of
$\{f_i\}$ in $\O_{\{\rho_{ji}\}}$ onto an open neighborhood of $\{f_k\}$
in $\O_{\{\rho_{lk}\}}$, where $\O_{\{\rho_{ji}\}}$, $\O_{\{\rho_{lk}\}}$
are given the topology generated by the subsets in $(2.1.2)$. As a
consequence, each $\O_{\{\rho_{ji}\}}$ is an open subset of $[(I(G),0,1);
(X,\underline{o_1},\underline{o_2})]_{(\rho,\eta)}$ via the embedding
$\sigma\mapsto [\sigma]$.

Now we shall present some standard constructions on the path spaces.
First of all, given any $(G,\rho,\eta)$-guided path $u$, one obtains
its inverse, denoted by $\nu(u)$, by pre-composing $u$ by the map
$\nu:(I(G),0,1)\rightarrow (I(G),1,0)$ sending $t\mapsto 1-t,\forall
t\in I$ and $g\mapsto g, \forall g\in G$, which results in a
$(G,\eta,\rho)$-guided path. Secondly, one may compose a
$(G,\rho,\eta)$-guided path $u_1$ with a $(G,\eta,\xi)$-guided path
$u_2$ to obtain a $(G,\rho,\xi)$-guided path $u_1\# u_2$, which is done
as follows. Pick representatives $\sigma_1=(\{f_{i,1}\},\{\rho_{ji}\}):
\Gamma\{I_i\}\rightarrow\Gamma\{U_{i^\prime}\}$ of $u_1$, where $i=0,1,
\cdots,n$, and $\sigma_2=(\{f_{k,2}\},\{\eta_{lk}\}):\Gamma\{J_k\}
\rightarrow\Gamma\{V_{k^\prime}\}$ of $u_2$, where $k=0,1,\cdots,m$.
Let $a$ be the index running from $0$ to $n+m$. We define intervals $H_a$
by setting $H_a=I_a$ for $0\leq a\leq n-1$, $H_{n}=I_n\cup J_0$, and
$H_a=J_{a-n}$ for $n+1\leq a\leq n+m$. We define local chart $W_a$ on $X$
by setting $W_a=U_a$ for $0\leq a\leq n-1$, $W_{n}=U_{o_2}$, and
$W_a=V_{a-n}$ for $n+1\leq a\leq n+m$. We define $f_a$ by setting $f_a=
f_{a,1}$ for $0\leq a\leq n-1$, $f_{n}=f_{n,1}\cup f_{0,2}$, and
$f_a=f_{a-n,2}$ for $n+1\leq a \leq n+m$. We define $\xi_{(a+1)a}
=\rho_{(a+1)a}$ for $0\leq a\leq n-1$, $\xi_{(n+1)n}=\eta_{10}$,
$\xi_{(a+1)a}=\eta_{a+1-n,a-n}$ for $n+1\leq a\leq n+m-1$. Then after
the reparametrization $t\mapsto 2t$, $g\mapsto g$, we obtain a
homomorphism $\sigma_1\#\sigma_2=(\{f_a\},\{\xi_{ba}\}):
\Gamma\{H_a\}\rightarrow\Gamma\{W_{a^\prime}\}$, whose
equivalence class is a $(G,\rho,\xi)$-guided path. We define $u_1\# u_2$
to be the equivalence class of $\sigma_1\#\sigma_2$, which is
clearly well-defined. Finally, given any $\Phi\in [(X,\underline{o_1},
\underline{o_2});(X^\prime,\underline{o_1^\prime},\underline{o_2^\prime})]
_{(\eta_1,\eta_2)}$, one has the mapping $\Phi_\#:[(I(G),0,1);(X,
\underline{o_1},\underline{o_2})]_{(\rho_1,\rho_2)}\rightarrow [(I(G),0,1);
(X^\prime,\underline{o_1^\prime},\underline{o_2^\prime})]_{(\rho_1^\prime,
\rho_2^\prime)}$ defined by $u\mapsto\Phi\circ u$, where $\rho_i^\prime=
\eta_i\circ\rho_i$ for $i=1,2$, and given any injective homomorphism
$\iota:H\rightarrow G$, one has the mapping
$\iota^\#:[(I(G),0,1);(X,\underline{o_1},
\underline{o_2})]_{(\rho,\eta)}\rightarrow [(I(H),0,1);(X,\underline{o_1},
\underline{o_2})]_{(\rho\circ\iota,\eta\circ\iota)}$ defined by pre-composing
each guided path $u$ by the map $(t,h)\mapsto (t,\iota(h)), \forall t\in I,
h\in H$. The following lemma is straightforward.

\begin{lem}
The mappings $\nu,\#,\Phi_\#$ and $\iota^\#$ of path spaces are all
continuous.
\end{lem}

Now we consider a special guided path space, the guided
loop space $[(S^1(G),\ast);(X,\underline{o})]_\rho$. There is a natural
base point in $[(S^1(G),\ast);(X,\underline{o})]_\rho$,
i.e., the constant guided loop $\tilde{o}$ defined by $(t,g)\mapsto
(\hat{o},\rho(g)), \forall t\in S^1, g\in G$.

\begin{lem}
The based topological space $([(S^1(G),\ast);(X,\underline{o})]_\rho,
\tilde{o})$ is an $H$-group with the homotopy associative multiplication
$\#$ and the homotopy inverse $\nu$. Moreover, the maps $\Phi_\#$,
$\iota^\#$ are homomorphisms of $H$-groups.
\end{lem}

\pf
We refer to \cite{Sw} for the definition of $H$-group. Here we only sketch
a proof that $\nu$ is a homotopy inverse, namely, both maps
$\#\circ (\nu,Id),\;\#\circ (Id,\nu):([(S^1(G),\ast);(X,\underline{o})]_\rho,
\tilde{o})\rightarrow ([S^1(G),\ast);(X,\underline{o})]_\rho,\tilde{o})$
are homopotic to the constant map into the base point $\tilde{o}$. The
homotopy associativity of $\#$ can be proven in the same vein. The assertion
on the maps $\Phi_\#$ and $\iota^\#$ are trivial. We leave the details to
the reader.

First of all, we introduce some notations. For any $s\in [0,1]$, we set
$$
I_s=([0,\frac{1}{2}(1-s)]\cup [\frac{1}{2}(1+s),1])/
\{\frac{1}{2}(1-s)\sim\frac{1}{2}(1+s)\}.
$$
Let $\beta_s:I\rightarrow I_s$, which are homeomorphisms for $s\neq 1$,
be defined by
$$
\beta_s(t)=\left\{\begin{array}{cc}
(1-s)t & 0\leq t\leq\frac{1}{2}\\
(1-s)(t-1)+1 & \frac{1}{2}\leq t\leq 1.
\end{array} \right. \leqno (2.1.3)
$$
For any homomorphism $\sigma=(\{f_i\},\{\rho_{ji}\})$, we denote by
$\nu(\sigma)$ the homomorphism obtained from $\sigma$ by performing the
reparametrization $t\mapsto 1-t, t\in I$.

We shall define a homotopy $F:([(S^1(G),\ast);(X,\underline{o})]_\rho,
\tilde{o})\times [0,1]\rightarrow ([(S^1(G),\ast);(X,\underline{o})]_\rho,
\tilde{o})$ between $\#\circ (Id,\nu)$ and the constant map into $\tilde{o}$.
The construction of a homotopy between the map $\#\circ (\nu,Id)$ and
the constant map into $\tilde{o}$ is completely parallel.

Given any guided loop $u$ which is represented by a homomorphism $\sigma$,
and any $s\in [0,1]$, we define the guided loop $F(u,s)$ as follows.
Observe that for any $s\in [0,1]$, the restriction of the homomorphism
$\sigma\#\nu(\sigma)$ to $I_s$ is a homomorphism. We reparametrize it
by $\beta_s$ in $(2.1.3)$, and denote it by $\sigma\#_s\nu(\sigma)$. If
$\sigma$ is replaced by an equivalent homomorphism, then $\sigma\#_s
\nu(\sigma)$ will be changed to an equivalent one accordingly. Hence the
$(G,\rho)$-guided loop defined by $\sigma\#_s\nu(\sigma)$ depends only on
$u$ and $s$. We define $F(u,s)=[\sigma\#_s\nu(\sigma)]$. Clearly,
$F(u,0)=\#\circ (Id,\nu)(u)$ and $F(u,1)=\tilde{o}$ for all $u\in
[(S^1(G),\ast);(X,\underline{o})]_\rho$.

It remains to verify that for any given $u_0\in [(S^1(G),\ast);
(X,\underline{o})]_\rho$, $s_0\in [0,1]$, $F$ is continuous at $(u_0,s_0)$.
We shall only give the details for the case when $s_0\neq 1$, the remaining
case $s_0=1$ is easier and we leave it to the reader. Without loss of
generality, we may assume that $u_0$ is represented by a homomorphism
$\sigma_0=(\{f_{i,0}\},\{\rho_{ji}\}):\Gamma\{I_i\}\rightarrow\Gamma
\{U_{i^\prime}\}$, where $i=0,1,\cdots,n$, such that there exist an $i_0$
and a closed interval $[s_{-},s_{+}]$ containing $s_0$ and satisfying
$s_{+}\neq 1$, $1-s_0\in [1-s_{+},1-s_{-}]\subset I_{i_0}\setminus
\cup_{i\neq i_0}I_i$. Then for any $\sigma=\{f_i\}\in\O_{\{\rho_{ji}\}}$
and $s\in [s_{-},s_{+}]$, we have $\sigma\#_s\nu(\sigma)$ given by
$(\{f_{k,s}\},\{\xi_{lk}\}):\Gamma\{J_{k,s}\}\rightarrow
\Gamma\{V_{k^\prime}\}$, where $k=0,1,\cdots,2i_0$, and
$$
J_{k,s}=\left\{\begin{array}{l}
\{\beta_s^{-1}(\frac{1}{2}t)|t\in I_k\}, \hspace{2mm}
                  0\leq k\leq i_0-1\\
\{\beta_s^{-1}(\frac{1}{2}t)|t\in I_{i_0},t\leq 1-s,
              \mbox{ or } 2-t\leq 1-s, 2-t\in I_{i_0}\}, \; k=i_0\\
\{\beta_s^{-1}(\frac{1}{2}t)|2-t\in I_{2i_0-k}\}, \hspace{2mm}
                i_0+1\leq k\leq 2i_0,\\
\end{array} \right. \leqno (2.1.4)
$$
$$
V_k=U_k,\; 0\leq k\leq i_0, \hspace{2mm} V_k=U_{2i_0-k},\;
i_0+1\leq k\leq 2i_0, \leqno (2.1.5)
$$
$$
f_{k,s}(t)=\left\{\begin{array}{l}
f_k(2\beta_s(t)),\hspace{2mm} 0\leq k\leq i_0-1\\
\left\{\begin{array}{l}
f_k(2\beta_s(t)),\; 2\beta_s(t)\leq 1-s, 2\beta_s(t)\in I_{i_0}\\
f_k(2-2\beta_s(t)),\; 2-2\beta_s(t)\leq 1-s, 2-2\beta_s(t)\in I_{i_0}\\
\end{array} \right. \; k=i_0\\
f_k(2-2\beta_s(t)), \hspace{2mm} i_0+1\leq k\leq 2i_0\\
\end{array} \right. \leqno (2.1.6)
$$
and
$$
\xi_{k(k-1)}=\rho_{k(k-1)},\; 1\leq k\leq i_0, \hspace{2mm}
\xi_{k(k-1)}=\rho_{(2i_0-k+1)(2i_0-k)}^{-1},\; i_0+1\leq k\leq 2i_0.
\leqno (2.1.7)
$$
Let $|s_{+}-s_{-}|\ll 0$ so that
$$
J_{k,s_{-}}\cap J_{k+1,s_{+}}\neq\emptyset,\; 0\leq k\leq i_0-1,
\hspace{2mm} J_{k,s_{+}}\cap J_{k+1,s_{-}}\neq\emptyset,\;
i_0\leq k\leq 2i_0. \leqno (2.1.8)
$$
We set $H_k=J_{k,s_{-}}\cap J_{k,s_{+}}$ for $0\leq k\leq 2i_0$.
Then it is easy to see that $\{H_k\}$ is a cover of $I$ by sub-intervals
such that $H_l\cap H_k\neq \emptyset$ iff $l=k$ or $k+1$, and $H_k\subset
J_{k,s}$ for any $0\leq k\leq 2i_0$ and $s\in [s_{-},s_{+}]$. We define
$g_{k,s}=f_{k,s}|_{H_k}$ for all $0\leq k\leq 2i_0$. Then
$\sigma_s=(\{g_{k,s}\},\{\xi_{lk}\}):\Gamma\{H_k\}\rightarrow
\Gamma\{V_{k^\prime}\}$, where $k=0,1,\cdots,2i_0$, is equivalent to
$\sigma\#_s\nu(\sigma)$. Observe that in the neighborhood
$\O_{\{\rho_{ji}\}}\times [s_{-},s_{+}]$ of $(u_0,s_0)$, the map $F$ is
given by $(\sigma,s)\mapsto \sigma_s\in\O_{\{\xi_{lk}\}}$, which is clearly
continuous with respect to
the topology of $\O_{\{\rho_{ji}\}}$ and $\O_{\{\xi_{lk}\}}$.
Hence $F$ is continuous.

\hfill $\Box$

We end by introducing the guided relative loop space
$[(I(G),S^0(G),0);(X,A,\underline{o})]_\rho$ for any pair
$(X,A,\underline{o})$ and injective homomorphism $\rho:G\rightarrow
G_{\hat{o}}$. Here $S^0=\{0,1\}$. Recall that an element of
$[(I(G),S^0(G),0);(X,A,\underline{o})]_\rho$ is a map (a guided path)
$u\in [(I(G),0);(X,\underline{o})]_\rho$ such that the restriction of
$u$ to the subspace $(S^0(G),0)$ is an element of $[(S^0(G),0);
(A,\underline{o}|_A)]_\rho$.

A canonical topology can be given to $[(I(G),S^0(G),0);
(X,A,\underline{o})]_\rho$ along the general lines of \S 3.2 in Part I
of this series \cite{C1}. For this purpose, we need to digress on
the subspace structure of $A$ in some details. Let $\U=\{U_i\}$ be the
atlas of local charts on $X$ and $\Gamma$ be the defining groupoid for $X$
with space of units $\bigsqcup_i\widehat{U_i}$. Denote by $\{V_\alpha\}$
the set of connected components of all $A\cap U_i$, $U_i\in\U$, and
by $\widehat{V_\alpha}$ some fixed component of the inverse
image of $V_{\alpha}$ in $\widehat{U_i}$ if $V_\alpha$ is a component of
$A\cap U_i$. Then the orbispace structure on $A$ is given by the
restriction of $\Gamma$ to $\bigsqcup_\alpha\widehat{V_\alpha}$.
In terms of local charts, each $\widehat{V_\alpha}$ is acted upon
by a group $G_{V_\alpha}$, which is the subgroup of $G_{U_i}$ that
fixes the subset $\widehat{V_\alpha}$, with a map $\pi_{V_\alpha}=
\pi_{U_i}|_{\widehat{V_\alpha}}$ which induces $\widehat{V_\alpha}/
G_{V_\alpha}\cong V_\alpha$. Moreover, there are natural mappings
$\rho_{\beta\alpha}:T(V_\alpha,V_\beta)\rightarrow T(U_i,U_j)$, where
$V_\alpha,V_\beta$ are components of $A\cap U_i,A\cap U_j$, such that
$\rho_\alpha=\rho_{\alpha\alpha}:G_{V_\alpha}\subset G_{U_i}$.

An element $u\in [(I(G),S^0(G),0);(X,A,\underline{o})]_\rho$ is
the equivalence class of $\sigma=(\{f_i\},\{\xi_{ji}\}):\Gamma\{I_i\}
\rightarrow\Gamma\{U_{i^\prime}\}$, where (1) $\{I_i\mid i=0,1,
\cdots,n\}$ is a cover of $I$ by sub-intervals such that $0\in
I_0$, $1\in I_n$, $I_i\cap I_j\neq\emptyset$ iff $j=i$ or $i+1$,
(2) denote by $U_i\in\{U_{i^\prime}\}$ the local chart assigned to $I_i$,
then $U_0=U_{o}$, and there is a $V_\alpha$ which is a component of $A\cap
U_n$, (3) $f_0(0)=\hat{o}$, $f_n(1)\in\widehat{V_\alpha}$, and (4)
$\xi_0=\rho$ and $\xi_n(G)\subset\rho_\alpha(G_{V_\alpha})$. Note
that each mapping $\xi_{(i+1)i}:T(I_i,I_{i+1})\rightarrow T(U_i,U_{i+1})$
is completely determined by the image of $1\in G=T(I_i,I_{i+1})$
in $T(U_i,U_{i+1})$, we shall conveniently regard $\xi_{(i+1)i}$
as an element of $T(U_i,U_{i+1})$. On the other hand, the
homomorphisms $\xi_i=\xi_{ii}:G\rightarrow G_{U_i}$ can be
recovered inductively by $\xi_{i+1}=\lambda_{\xi_{(i+1)i}}\circ\xi_i$
with $\xi_0=\rho$. We observe that the image of $f_i:I_i\rightarrow
\widehat{U_i}$ for each index $i$ is forced to lie in the fixed-point set
of $\xi_i(G)\subset G_{U_i}$. Finally, the restriction of $u$ to
the subspace $(S^0(G),0)$ is the map into $(A,\underline{o}|_A)$ defined
by $(0,g)\mapsto (\hat{o},\rho(g)), (1,g)\mapsto (f_n(1),\rho_\alpha^{-1}\circ
\xi_n(g))$, $\forall g\in G$.

We introduce
$$
\O_{(\{\xi_{ji}\},V_\alpha)}=\{\sigma\mid\sigma=(\{f_i\},\{\xi_{ji}\})
\mbox{ is admissible},
f_n(1)\in\widehat{V_\alpha},\xi_n(G)\subset\rho_\alpha(G_{V_\alpha})\}.
\leqno (2.1.9)
$$
By Lemma 3.1.2 in \cite{C1}, $\O_{(\{\xi_{ji}\},V_\alpha)}$ can be regarded
as a subset of $[(I(G),S^0(G),0);(X,A,\underline{o})]_\rho$ via $\sigma
\mapsto [\sigma]$. We give each $\O_{(\{\xi_{ji}\},V_\alpha)}$ a topology
which is generated by
$$
\O_{(\{\xi_{ji}\},V_\alpha)}(\{\underline{K_i}\},\{\underline{O_i}\})=
\{\{f_i\}\in\O_{(\{\xi_{ji}\},V_\alpha)}\mid f_i(K_{i,s})\subset
O_{i,s},\forall s\in\Lambda(i)\} \leqno (2.1.10)
$$
where $\underline{K_i}=\{K_{i,s}|s\in\Lambda(i)\}$ is any finite set
of compact subsets of $I_i$, and $\underline{O_i}=\{O_{i,s}\mid s\in
\Lambda(i)\}$ is any finite set of open subsets of $\widehat{U_i}$.
Again by Lemma 3.2.2 in \cite{C1}, for any $\{f_i\}\in\O_{(\{\xi_{ji}\},
V_\alpha)}$ and $\{f_k\}\in\O_{(\{\xi_{lk}\},V_\beta)}$ which are
equivalent, there is a local homeomorphism $\phi$ from an open
neighborhood of $\{f_i\}$ in $\O_{(\{\xi_{ji}\},V_\alpha)}$ onto
an open neighborhood of $\{f_k\}$ in $\O_{(\{\xi_{lk}\},V_\beta)}$.
We give a topology to $[(I(G),S^0(G),0);(X,A,\underline{o})]_\rho$
which is generated by the subsets in $(2.1.10)$ for all possible
data $\{\xi_{ji}\},V_\alpha,\{\underline{K_i}\}$ and
$\{\underline{O_i}\}$. However, because of the existence of the said
local homeomorphisms $\{\phi\}$, each $\O_{(\{\xi_{ji}\},V_\alpha)}$
is in fact an open subset of
$[(I(G),S^0(G),0);(X,A,\underline{o})]_\rho$.

There is a special element $\tilde{o}\in [(I(G),S^0(G),0);(X,A,
\underline{o})]_\rho$ which is the equivalence class of the constant
guided relative loop, represented by any $\{f_i\}\in\O_{(\{\xi_{ji}\},
V_\alpha)}$ where each $f_i:I_i\rightarrow\widehat{U_i}=\widehat{U_o}$
has the point $\hat{o}$ as its image, each $\xi_{(i+1)i}=1\in G_{U_o}=
T(U_i,U_j)$, and $V_\alpha=V_o$. We shall fix $\tilde{o}$ as the
base point of the guided relative loop space.

We may compose a guided path with a guide relative loop. More precisely,
Let $u_1\in [(I(G),0,1);(X,\underline{o_1},\underline{o_2})]_{(\rho,\eta)}$
be a guided path, and $u_2\in [(I(G),S^0(G),0);(X,A,\underline{o_2})]_\eta$
be a guided relative loop. Then the composition $u_1\# u_2$ is
well-defined, which is an element in the guided relative loop space
$[(I(G),S^0(G),0);(X,A,\underline{o_1})]_\rho$.

The mappings $\Phi_\#:[(I(G),S^0(G),0);(X,A,\underline{o})]_\rho
\rightarrow [(I(G),S^0(G),0);(X^\prime,A^\prime,\underline{o^\prime})]
_{\eta\circ\rho}$ and $\iota^\#: [(I(G),S^0(G),0);(X,A,\underline{o})]_\rho
\rightarrow [(I(H),S^0(H),0);(X,A,\underline{o})]_{\rho\circ\iota}$
are defined for any $\Phi\in [(X,A,\underline{o});(X^\prime,A^\prime,
\underline{o^\prime})]_\eta$ and injective homomorphism $\iota:
H\rightarrow G$.  The following lemma is straightforward.

\begin{lem}
The mappings $\#,\Phi_\#$ and $\iota^\#$ are all continuous.
\end{lem}

\subsection{Generalities of homotopy groups}

We first derive two preliminary lemmas.

Let $(K,x_0)$ be a compact, locally connected topological space with
base point $x_0\in K$. The suspension of $(K,x_0)$ is the based space
$(SK,\ast)$, where
$$
SK=(I\times K)/(\{0,1\}\times K)\cup (I\times \{x_0\}), \leqno (2.2.1)
$$
with the base point $\ast\in SK$ being the image of $(\{0,1\}\times K)
\cup (I\times \{x_0\})$. The cone of $(K,x_0)$ is the based space
$(CK,\ast)$, where
$$
CK=(I\times K)/(\{0\}\times K)\cup (I\times \{x_0\}), \leqno (2.2.2)
$$
with the base point $\ast\in CK$ being the image of $(\{0\}\times K)
\cup (I\times \{x_0\})$. The based space $(K,x_0)$ is canonically a
subspace of $(CK,\ast)$ via the embedding $x\mapsto [1,x]$ where
$[1,x]$ denotes the image of $(1,x)\in I\times K$ in $CK$.

\begin{lem}
\begin{itemize}
\item [{(1)}] The set of continuous maps from $(K,x_0)$ into
$([(S^1(G),\ast);(X,\underline{o})]_\rho,\tilde{o})$ is naturally
identified with the set $[(SK(G),\ast);(X,\underline{o})]_\rho$.
\item [{(2)}] The set of continuous maps from $(K,x_0)$ into
$([(I(G),S^0(G),0);(X,A,\underline{o})]_\rho,\tilde{o})$ is
naturally identified with the set $[(CK(G),K(G),\ast);
(X,A,\underline{o})]_\rho$.
\end{itemize}
\end{lem}

\pf
(1) Given any continuous map $u:(K,x_0)\rightarrow ([(S^1(G),\ast);
(X,\underline{o})]_\rho,\tilde{o})$, we construct an element
$\Phi_u\in [(SK(G),\ast);(X,\underline{o})]_\rho$ as follows. Since
$u(K)\subset [(S^1(G),\ast);(X,\underline{o})]_\rho$ is a compact
subset, there are finitely many open subsets
$$
\{\O_\alpha|\O_\alpha=\O_{\{\rho_{ji,\alpha}\}},\alpha\in\Lambda,
\#\Lambda<+\infty\} \leqno (2.2.3)
$$
such that $u(K)\subset \cup_{\alpha\in\Lambda}\O_\alpha$. Furthermore,
we may require that (a) $u(x_0)=\tilde{o}$ is contained in $\O_{\alpha_0}$
where $U_{i,\alpha_0}=U_o$, $\rho_{ji,\alpha_0}=1\in G_{U_o}$ for
all $i,j$, and (b) $\{I_{i,\alpha}\}=\{I_i\}$ is independent of $\alpha$.
The latter can be always achieved because $\#\Lambda<+\infty$.

We take a cover $\{O_a|a\in A\}$ of $K$, where each $O_a$ is a
connected open subset, together with $\theta:a\mapsto \alpha$
such that $u(O_a)\subset\O_{\theta(a)}$. Then for each $a\in A$, there
is a set of data
$$
(\{I_i\},\{U_{i,\theta(a)}\},\{f_{i,a}\},\{\rho_{ji,\theta(a)}\}),
\leqno (2.2.4)
$$
where each $f_{i,a}:I_i\times O_a\rightarrow\widehat{U_{i,\theta(a)}}$
is continuous such that for any $x\in O_a$, the restriction of $(2.2.4)$
to $x$ is a homomorphism $\{f_{i,a}(\cdot,x)\}\in\O_{\theta(a)}$
representing $u(x)\in [(S^1(G),\ast);(X,\underline{o})]_\rho$.

Let $\{O_{ab,s}\mid s\in I_{ab}\}$ be the set of connected components of
$O_a\cap O_b$, $a,b\in A$. For any $x\in O_{ab,s}$, since both
$(\{I_i\},\{U_{i,\theta(a)}\},\{f_{i,a}(\cdot,x)\},
\{\rho_{ji,\theta(a)}\})$ and $(\{I_i\},\{U_{i,\theta(b)}\},
\{f_{i,b}(\cdot,x)\},\{\rho_{ji,\theta(b)}\})$ represent
$u(x)\in [(S^1(G),\ast);(X,\underline{o})]_\rho$, we may apply
Lemma 2.2.4 and then Lemma 3.1.2 of \cite{C1} to conclude that there
exists a set of elements $\xi_i^{ba,s}(x)\in T(U_{i,\theta(a)},
U_{i,\theta(b)})$ such that
$$
\begin{array}{l}
f_{i,b}(\cdot,x)=\phi_{\xi_i^{ba,s}(x)}\circ f_{i,a}(\cdot,x), \;
\forall i,\\
\rho_{ji,\theta(b)}=\xi_j^{ba,s}(x)\circ\rho_{ji,\theta(a)}\circ
\xi_i^{ba,s}(x)^{-1}({\bf a}),\; \forall i,j,
\end{array} \leqno (2.2.5)
$$
where ${\bf a}\in\Lambda(\xi_i^{ba,s}(x)^{-1},\rho_{ji,\theta(a)},
\xi_j^{ba,s}(x))$ is the element containing $f_{i,b}
(I_i\cap I_j\times \{x\})$. Moreover, the elements $\xi_0^{ba,s}(x)\in
T(U_{0,\theta(a)},U_{0,\theta(b)})$ have to satisfy $\xi_0^{ba,s}(x)=1
\in G_{U_o}$ (note that $U_{0,\alpha}=U_o$ for all $\alpha\in\Lambda$).
Now the second equation in $(2.2.5)$ implies inductively that each
$\xi_i^{ba,s}(x)$ is constant in $x$, depending only on the connected
component $O_{ab,s}$ that contains $x$. We denote the constant value by
$\xi_i^{ba,s}$. It also follows from the second equation in $(2.2.5)$
that if a connected component ${\bf b}$ of $O_{ab,s}\cap O_{bc,t}$ is
contained in $O_{ac,r}$, then
$$
\xi^{cb,t}_i\circ\xi^{ba,s}_i(f_{i,a}({\bf b}))=\xi_i^{ca,r}.
\leqno (2.2.6)
$$

It is easily seen that $\tau(\theta,\{O_a\},\{\O_\alpha\})=(\{f_{i,a}\},
\{\rho_{ji,\theta(b)}\circ\xi_i^{ba,s}\})$, where $\rho_{ji,\theta(b)}
\circ\xi_i^{ba,s}$ is the element of $T(U_{i,\theta(a)},U_{j,\theta(b)})$
determined by the component containing $f_{i,a}(I_j\cap I_i\times O_{ab,s})$,
is a groupoid homomorphism from $\Gamma\{I_i\times O_a\}$
to $\Gamma\{U_{i^\prime,\alpha}\}$, whose equivalence class is an element
of $[(SK(G),\ast);(X,\underline{o})]_\rho$. We define $\Phi_u$ to be the
equivalence class of $\tau(\theta,\{O_a\},\{\O_\alpha\})$.

It remains to verify that $\Phi_u$ is well-defined. For this
purpose, let $\tau(\theta,\{O_a\},\{\O_\alpha\})$ and
$\tau(\theta^\prime,\{O_{a^\prime}\},\{\O_{\alpha^\prime}\})$ be
any two such homomorphisms. We shall construct a common induced homomorphism
of them as follows. For any $x\in K$, we pick an $\alpha$ and an
$\alpha^\prime$ such that $u(x)$ is contained in $\O_\alpha\cap
\O_{\alpha^\prime}$. Then by Lemma 3.2.2 in \cite{C1}, there is an open
subset $\O_x$ of some $\O_{\{\rho_{lk,x}\}}$, together with
open embeddings $\phi_\alpha:\O_x\rightarrow\O_\alpha$,
$\phi_{\alpha^\prime}:\O_x\rightarrow\O_{\alpha^\prime}$, such that
$u(x)$ is contained in $\O_x$ and each $\kappa\in\O_x$ is a common
induced homomorphism of $\phi_\alpha(\kappa)\in\O_\alpha$ and
$\phi_{\alpha^\prime}(\kappa)\in \O_{\alpha^\prime}$. Since $u(K)$ is
compact in $[(S^1(G),\ast);(X,\underline{o})]_\rho$, we can cover it by
finitely many such open subsets $\O_r=\O_{x_r}, r=1,2,\cdots,N$. Note that
there are mappings $\imath:r\mapsto\alpha$, $\imath^\prime:r\mapsto
\alpha^\prime$ such that $\O_r$ is mapped into $\O_{\imath(r)}$,
$\O_{\imath^\prime(r)}$ under the open embeddings $\phi_{\imath(r)}$,
$\phi_{\imath^\prime(r)}$ respectively. On the other hand, we take a
cover $\{O_s\}$ of $K$ by connected open subsets, together with mappings
$\jmath:s\mapsto a$, $\jmath^\prime:s\mapsto a^\prime$ such that $O_s\subset
O_{\jmath(s)}\cap O_{\jmath^\prime(s)}$, and with a mapping $\bar{\theta}:
s\mapsto r$ satisfying $u(O_s)\subset\O_{\bar{\theta}(s)}$ and
$\imath\circ\bar{\theta}=\theta\circ\jmath$, $\imath^\prime\circ
\bar{\theta}=\theta^\prime\circ\jmath^\prime$. Then the set of data
$(\bar{\theta},\{O_s\},\{\O_r\})$ gives rise to a homomorphism
$\tau(\bar{\theta},\{O_s\},\{\O_r\})$, which is induced by both
$\tau(\theta,\{O_a\},\{\O_\alpha\})$ and
$\tau(\theta^\prime,\{O_{a^\prime}\},\{\O_{\alpha^\prime}\})$. Hence
$\Phi_u$ is well-defined.

We denote the correspondence $u\mapsto \Phi_u$ by $\phi$. Then the
inverse of $\phi$ is obtained as follows. Given any $\Phi\in
[(SK(G),\ast);(X,\underline{o})]_\rho$, we define a map $u_{\Phi}:
(K,x_0)\rightarrow ([(S^1(G),\ast);(X,\underline{o})]_\rho,\tilde{o})$
by the rule that for any $x\in K$, $u_{\Phi}(x)$ is the restriction to
the subspace $\{[t,x]|t\in I\}$ of $(SK(G),\ast)$, where $[t,x]$
denotes the image of $(t,x)\in I\times K$ in $SK$. It is easily seen
that the map $u_{\Phi}$ is continuous, and the correspondence
$\psi:\Phi\mapsto u_{\Phi}$ is the inverse of $\phi$, namely,
$u_{\Phi_u}=u$ and $\Phi_{u_{\Phi}}=\Phi$. Hence (1) of the lemma.

\vspace{2mm}

(2) We shall only sketch the proof here since it is completely parallel
to the one above.

Given any continuous map $u:(K,x_0)\rightarrow ([(I(G),S^0(G),0);(X,A,
\underline{o})]_\rho,\tilde{o})$, we cover $u(K)$ by finitely many
$\O_k,k=0,1,\cdots,m$, where $\O_k=\O_{(\{\xi_{ji,k}\},V_{\alpha(k)})}$
with $i=0,1,\cdots,n$, such that $u(x_0)=\tilde{o}$ is contained in $\O_0$
with $U_{i,0}=U_o$, $\xi_{ji,0}=1\in G_{U_o}$ for all $i,j$. We then take
a cover $\{O_a\}$ of $K$ by connected open subsets, together with a mapping
$\theta:a\mapsto k$ satisfying $u(O_a)\subset\O_{\theta(a)}$. Then for each
index $a$, there is a set of data
$$
(\{I_i\},\{U_{i,\theta(a)}\},V_{\alpha(\theta(a))},\{f_{i,a}\},
\{\xi_{ji,\theta(a)}\}), \; i=0,1,\cdots,n, \leqno (2.2.7)
$$
where each $f_{i,a}:I_i\times O_a\rightarrow\widehat{U_{i,\theta(a)}}$
is continuous such that the restriction of $(2.2.7)$ to each $x\in O_a$
is an element $\{f_{i,a}(\cdot,x)\}\in\O_{\theta(a)}$, which
represents $u(x)\in [(I(G),S^0(G),0);(X,A,\underline{o})]_\rho$.
In particular, $f_{n,a}(1,x)\in\widehat{V_{\alpha(\theta(a))}}$ for all
$x\in O_a$. We define $g_a:O_a\rightarrow\widehat{V_{\alpha(\theta(a))}}$
by $x\mapsto f_{n,a}(1,x)$.

For each connected component $O_{ab,s}$ of $O_a\cap O_b$, there exists
a set of elements $\eta_i^{ba,s}\in T(U_{i,\theta(a)},U_{i,\theta(b)})$,
$i=0,1,\cdots,n$, where $\eta_0^{ba,s}=1\in G_{U_o}$ and $\eta_n^{ba,s}
=\rho_{\alpha(\theta(b))\alpha(\theta(a))}(\zeta^{ba,s})$ for some
$\zeta^{ba,s}\in T(V_{\alpha(\theta(a))},V_{\alpha(\theta(b))})$,
such that $g_a(O_{ab,s})$ is contained in $\mbox{Domain }(
\phi_{\zeta^{ba,s}})$, and the following equations are satisfied:
$$
\begin{array}{l}
f_{i,b}(\cdot,x)=\phi_{\eta_i^{ba,s}}\circ f_{i,a}(\cdot,x), \;
\forall x\in O_{ab,s}, \forall i,\\
\xi_{ji,\theta(b)}=\eta_j^{ba,s}\circ\xi_{ji,\theta(a)}\circ
(\eta_i^{ba,s})^{-1}\; \forall i,j,
\end{array} \leqno (2.2.8)
$$
where $\eta_j^{ba,s}\circ\xi_{ji,\theta(a)}\circ (\eta_i^{ba,s})^{-1}$ is
the element in $T(U_{i,\theta(b)},U_{j,\theta(b)})$ that is determined by
the component containing $f_{i,b}(I_j\cap I_i\times O_{ab,s})$. Analogous
to $(2.2.6)$, we have
$$
\eta^{cb,t}_i\circ\eta^{ba,s}_i(f_{i,a}({\bf b}))=\eta_i^{ca,r}
\leqno (2.2.9)
$$
for any component ${\bf b}$ of $O_{ab,s}\cap O_{bc,t}$ which is contained
in $O_{ac,r}$. Now $\tau(\theta,\{O_a\},\{\O_k\})=(\{f_{i,a}\},
\{\xi_{ji,\theta(b)}\circ\eta_i^{ba,s}\})$, where $\xi_{ji,\theta(b)}
\circ\eta_i^{ba,s}$ is the element in $T(U_{i,\theta(a)},U_{j,\theta(b)})$
determined by the component containing $f_{i,a}(I_j\cap I_i\times O_{ab,s})$,
is a homomorphism from $\Gamma\{I_i\times O_a\}$ to
$\Gamma\{U_{i^\prime,k}\}$. We define $\Phi_u$ to be the equivalence class
of $\tau(\theta,\{O_a\},\{\O_k\})$, which is well-defined by a similar
argument as in (1) above, and is clearly an element of
$[(CK(G),\ast);(X,\underline{o})]_\rho$.

We need to verify that the restriction of $\Phi_u$ to the subspace
$(K(G),\ast)$ is an element of $[(K(G),\ast);(A,\underline{o}|_A)]_\rho$
so that $\Phi_u\in [(CK(G),K(G),\ast);(X,A,\underline{o})]_\rho$.
Now the restriction of $\tau(\theta,\{O_a\},\{O_k\})$ to $(K(G),\ast)$
is clearly the homomorphism $(\{g_a\},\{\zeta^{ba,s}\}):\Gamma\{O_a\}
\rightarrow\Gamma\{V_{\alpha^\prime(\theta(a))}\}$, which defines
an element of $[(K(G),\ast);(A,\underline{o}|_A)]_\rho$. In order
to justify that $(\{g_a\},\{\zeta^{ba,s}\})$ is indeed a homomorphism,
a number of equations needs to be verified: The first equation in
$(2.2.8)$ for the case $i=n$ implies that $g_b(x)=\phi_{\zeta^{ba,s}}\circ
g_a(x),\;\forall x\in O_{ab,s}$, and $(2.2.9)$ for the case $i=n$ implies
that $\zeta^{cb,t}\circ\zeta^{ba,s}(g_a({\bf b}))=\zeta^{ca,r}$. Moreover,
for each index $a$, recall that $\{\xi_{i,a}:G\rightarrow
G_{U_{i,\theta(a)}}\}$ may be inductively determined by $\xi_{i+1,a}=
\lambda_{\xi_{(i+1)i,\theta(a)}}\circ\xi_{i,a}$ with $\xi_{0,a}=\rho:
G\rightarrow G_{U_o}$ and $\xi_{n,a}(G)\subset\rho_{\alpha(\theta(a))}
(G_{V_{\alpha(\theta(a))}})$. We define $\zeta^a=\rho_{\alpha(\theta(a))}^{-1}
\circ\xi_{n,a}:G\rightarrow G_{V_{\alpha(\theta(a))}}$. Then the second
equation in $(2.2.8)$ implies that $\zeta^b=\lambda_{\zeta^{ba,s}}\circ
\zeta^a$.

The correspondence $\phi:u\mapsto \Phi_u$ is a bijection whose inverse
$\psi$ is defined as follows. For any $\Phi\in [(CK(G),K(G),\ast);
(X,A,\underline{o})]_\rho$, let $u_{\Phi}(x)$ be the restriction of
$\Phi$ to the subspace $\{[t,x]\mid t\in I\}$, which is an element
of $[(I(G),S^0(G),0);(X,A,\underline{o})]_\rho$. The mapping $x\mapsto
u_{\Phi}(x)$ is easily seen to be continuous. The identities $u_{\Phi_u}=u$
and $\Phi_{u_{\Phi}}=\Phi$ are clear from the construction.
Hence (2) of the lemma.

\hfill $\Box$

The second lemma is concerned with the special case when $X=Y/G$ is a
global quotient.

We begin by fixing the notations. Let $\underline{o}=(o,U_o,\hat{o})$ be
any base-point structure of $X$. We denote by $U^o$ the connected component
of $X$ that contains $o$. The space $\widehat{U^o}$ is a connected component
of $Y$ that contains $\widehat{U_o}$, hence $\hat{o}\in\widehat{U_o}$ may
be regarded as a point in $\widehat{U^o}\subset Y$. Moreover, $T(U_o,U^o)$
is canonically identified with $G_{U^o}$. We assume $H\subset G$ is a subgroup
and denote by $\rho:H\rightarrow G$ the inclusion. We set $Y^H=\{y\in Y\mid
h\cdot y=y,\;\forall h\in H\}$ and $C(H)=\{g\in G\mid gh=hg,\;\forall h\in H\}$.

Let $P(Y^H,\hat{o})=\{\gamma\mid\gamma:I\rightarrow Y^H \mbox{ continuous },
\gamma(0)=\hat{o}\}$ given with the usual ``compact-open'' topology,
and let $P(Y^H,C(H),\hat{o})=\{(\gamma,g)\mid\gamma\in P(Y^H,\hat{o}),
g\in C(H), \mbox{ s.t. } \gamma(1)=g\cdot\hat{o}\}$ given with the
relative topology as a subspace of $P(Y^H,\hat{o})\times C(H)$. The
space $P(Y^H,C(H),\hat{o})$ has a natural base point $\tilde{o}=
(\gamma_o,1)$ where $\gamma_o$ is the constant map to $\hat{o}$.

Let $\Phi:(X,\underline{o})\rightarrow (X^\prime,\underline{o^\prime})$
be any map\footnote{here $\Phi$ is allowed to be a general map, i.e., $\Phi$
is not in the more restricted class specified by {\bf Convention} in
Introduction.} defined by a pair $(f,\lambda):(Y,G)\rightarrow (Y^\prime,G^\prime)$
where $f$ is $\lambda$-equivariant, and $\hat{o^\prime}=f(\hat{o})$ in the
base-point structure $\underline{o^\prime}=(o^\prime,U_o^\prime,
\hat{o^\prime})$. For any subgroup $H\subset G$, we set $H^\prime=\lambda(H)$
and denote the inclusion $H^\prime\subset G^\prime$ by $\rho^\prime$. There is
an induced map $P(f,\lambda):(P(Y^H,C(H),\hat{o}),\tilde{o})\rightarrow
(P((Y^\prime)^{H^\prime},C(H^\prime),\hat{o}^\prime),\tilde{o}^\prime)$
defined by $(\gamma,g)\mapsto (f\circ\gamma,\lambda(g))$.

\begin{lem}
There exists a $\phi_X:([(S^1(H),\ast);(X,\underline{o})]_\rho,
\tilde{o})\cong (P(Y^H,C(H),\hat{o}),\tilde{o})$ satisfying
$P(f,\lambda)\circ\phi_X=\phi_{X^\prime}\circ\Phi_\#$ for any map
$\Phi$ defined by $(f,\lambda):(Y,G)\rightarrow (Y^\prime,G^\prime)$.
\end{lem}

\pf
First of all, we define $\phi_X:([(S^1(H),\ast);(X,\underline{o})]_\rho,
\tilde{o})\rightarrow (P(Y^H,C(H),\hat{o}),\tilde{o})$ as follows.
Let $u\in [(S^1(H),\ast);(X,\underline{o})]_\rho$ be any
guided loop which is represented by $\sigma=\{f_i\}\in\O_{\{\xi_{ji}\}}$,
$i=0,1,\cdots,n$. Since the image of $u$ in the underlying space $X$
lies in $U^o$, we may assume that $U_i=U^o$ for each $i\neq 0,n$ by
replacing $\sigma$ with an equivalent homomorphism. Consequently,
each $\xi_{ji}\in T(U_i,U_j)$ is an element of $G_{U^o}$.
We define a path $\gamma:I\rightarrow Y$ as follows: we set
$\gamma=f_0$ on $I_{0}$, $\gamma=\phi_{\xi_{10}}^{-1}\circ f_1$ on
$I_1$, $\cdots, \gamma=\phi_{\xi_{10}}^{-1}\circ\cdots\circ
\phi_{\xi_{n(n-1)}}^{-1}\circ f_n$ on $I_n$, and we define $g=
\xi_{10}^{-1}\circ\cdots\circ\xi_{n(n-1)}^{-1}\in G_{U^o}\subset G$.
Clearly $\gamma(0)=\hat{o}$, $\gamma(t)\in Y^H,\forall t\in I$ and
$\gamma(1)=g\cdot\hat{o}$. Moreover, $Ad(g)=\lambda_{\xi_{10}}^{-1}\circ
\cdots\circ\lambda_{\xi_{n(n-1)}}^{-1}$ which satisfies $\rho=Ad(g)^{-1}
\circ\rho$. Hence $Ad(g)(h)=h,\forall h\in H$ so that $g\in C(H)$. We define
$\phi_X$ by setting $\phi_X(u)=(\gamma,g)$, which is clearly an
element of $P(Y^H,C(H),\hat{o})$. It is easily seen that $\phi_X$ is
well-defined, i.e., $\phi_X(u)$ is independent of the choice on $\sigma$.
The map $\phi_X$ is continuous as well. Furthermore, $\phi_X$ is a base-point
preserving map, and satisfies $P(f,\lambda)\circ\phi_X=\phi_{X^\prime}\circ
\Phi_\#$ for any map $\Phi$ defined by $(f,\lambda):(Y,G)\rightarrow
(Y^\prime,G^\prime)$.

It remains to show that $\phi_X$ is a homeomorphism. We construct a
$\psi_X:P(Y^H,C(H),\hat{o})\rightarrow [(S^1(H),\ast);(X,\underline{o})]_\rho$
as follows. Given any $(\gamma,g)\in P(Y^H,C(H),\hat{o})$, it is obvious
that $\gamma(t)\in \widehat{U^o}$ for all $t\in I$, and $g\in G_{U^o}$
because $g\cdot \widehat{U^o}=\widehat{U^o}$. We define $\psi_X$ by
sending $(\gamma,g)$ to the element in $[(S^1(H),\ast);(X,\underline{o})]_\rho$
which is the equivalence class of $\{f_i\}\in\O_{\{\xi_{ji}\}}$, $i=0,1,2$,
where $I_0\subset\gamma^{-1}(\widehat{U_o})$, $I_1=(0,1)$,
$I_2\subset\gamma^{-1}(g\cdot\widehat{U_o})$, $U_0=U_o$, $U_1=U^o$, $U_2=U_o$,
$f_0=\gamma|_{I_0}$, $f_1=\gamma|_{I_1}$, $f_2=g^{-1}\circ\gamma|_{I_2}$, and
$\xi_{10}=1$, $\xi_{21}=g^{-1}$. Clearly $\psi_X$ is continuous, and
$\phi_X\circ\psi_X=Id$, $\psi_X\circ\phi_X=Id$. Hence $\phi_X$ is a
homeomorphism.

\hfill $\Box$

Now we present the first list of properties of the homotopy sets
defined in Definition 1.1.

\begin{prop}
\begin{itemize}
\item [{(1)}] There are natural group structures on $\pi_k^{(G,\rho)}
(X,\underline{o})$ and $\pi_{k+1}^{(G,\rho)}(X,A,\underline{o})$
for $k\geq 1$ which are Abelian when $k\geq 2$.
\item [{(2)}] For any $\Phi\in [(X,\underline{o});
(X^\prime,\underline{o^\prime})]_\eta$ (resp. $\Phi\in [(X,A,\underline{o});
(X^\prime,A^\prime,\underline{o^\prime})]_\eta$), there are natural
homomorphisms $\Phi_\ast:\pi_k^{(G,\rho)}(X,\underline{o})\rightarrow
\pi_k^{(G,\rho^\prime)}(X^\prime,\underline{o^\prime})$ (resp.
$\Phi_\ast:\pi_k^{(G,\rho)}(X,A,\underline{o})\rightarrow
\pi_k^{(G,\rho^\prime)}(X^\prime,A^\prime,\underline{o^\prime})$)
with $\rho^\prime=\eta\circ\rho$, which depend only on the
homotopy class of $\Phi$. For any $(H,\eta)$ where $\eta:H\rightarrow
G_{\hat{o}}$ is an injective homomorphism which factors through $\rho:
G\rightarrow G_{\hat{o}}$ by $\iota:H\rightarrow G$, there are
natural homomorphisms $\iota^\ast:\pi_k^{(G,\rho)}(X,\underline{o})
\rightarrow\pi_k^{(H,\eta)}(X,\underline{o})$ and $\iota^\ast:
\pi_{k+1}^{(G,\rho)}(X,A,\underline{o})\rightarrow\pi_{k+1}^{(H,\eta)}
(X,A,\underline{o})$.
\item [{(3)}] For the case $G=\{1\}$, $\pi_k(X,\underline{o})$,
$\pi_{k+1}(X,A,\underline{o})$ are naturally isomorphic to
$\pi_k(B\Gamma_X,\ast)$ and $\pi_{k+1}(B\Gamma_X,B\Gamma_A,\ast)$, where
$B\Gamma_X, B\Gamma_A$ are the classifying spaces of the defining groupoid
for $X$ and $A$ respectively.
\item [{(4)}] When $X=Y/G$ is a global quotient, for any subgroup $\rho:
H\subset G$, there are natural isomorphisms $\theta_X:\pi_k^{(H,\rho)}
(X,\underline{o})\cong \pi_k(Y^H,\hat{o})$ for all $k\geq 2$, and the
natural exact sequence
$$
1\rightarrow \pi_1(Y^H,\hat{o})\rightarrow\pi_1^{(H,\rho)}(X,\underline{o})
\rightarrow C(H)\rightarrow \pi_0(Y^H,\hat{o})\rightarrow
\pi_0^{(H,\rho)}(X,\underline{o}).\leqno (2.2.10)
$$
\end{itemize}
\end{prop}

\pf
(1) For any $k\geq 1$, we fix an identification $S^k=SS^{k-1}$ and
$D^{k+1}=CS^k$. Then by Lemma 2.2.1, $\pi_k^{(G,\rho)}(X,\underline{o})$
and $\pi_{k+1}^{(G,\rho)}(X,A,\underline{o})$ are the $\pi_{k-1}$ and
$\pi_k$ of $([(S^1(G),\ast);(X,\underline{o})]_\rho,\tilde{o})$ and
$([(I(G),S^0(G),0);(X,A,\underline{o})]_\rho,\tilde{o})$ respectively.
Hence they have natural group structures for $k\geq 1$ which are Abelian
when $k\geq 2$, cf. Lemma 2.1.2.

\vspace{1.5mm}

(2) Straightforward.

\vspace{1.5mm}

(3) Note that a map from a topological space $S$ (which is regarded
as an orbispace trivially) to an orbispace $X$ is an equivalence
class of $\Gamma_X$-structures on $S$, where $\Gamma_X$ is ranging in
a certain set of \'{e}tale topological groupoids which define
the orbispace structure on $X$. It is easily seen that when a
choice of $\Gamma_X$ is fixed, a map from $S$ to $X$ can be naturally
identified with a $\Gamma_X$-structure on $S$ (cf. Lemma 3.1.2 in \cite{C1}).
On the other hand, by a theorem of Haefliger in \cite{Ha1}, the set of
homotopy classes of $\Gamma_X$-structures on $S$ is naturally in one to
one correspondence with the set of homotopy classes of continuous maps
from $S$ into the classifying space $B\Gamma_X$ of $\Gamma_X$. The
assertion follows easily.

\vspace{1.5mm}

(4) We consider the continuous map $\pi:(P(Y^H,C(H),\hat{o}),\tilde{o})
\rightarrow (C(H),1)$ by sending $(\gamma,g)$ to $g$. The fiber at
$1\in C(H)$ is identified with the loop space $\Omega(Y^H,\hat{o})$ via
the embedding $\Omega(Y^H,\hat{o})\hookrightarrow P(Y^H,C(H),\hat{o})$
sending $\gamma$ to $(\gamma,1)$. With Lemma 2.2.2, the assertion follows
essentially from the fact that $\pi:P(Y^H,C(H),\hat{o})\rightarrow C(H)$
is a fibration (note that $C(H)$ has a discrete topology), and that
$\pi\circ\phi_X(u_1\# u_2)=\pi\circ\phi_X(u_1)\cdot\pi\circ\phi_X(u_2)$
for any guided loops $u_1,u_2\in [(S^1(H),\ast);(X,\underline{o})]_\rho$.
What remains is to examine the part $C(H)\rightarrow \pi_0(Y^H,\hat{o})
\rightarrow\pi_0^{(H,\rho)}(X,\underline{o})$ in $(2.2.10)$.

The map $C(H)\rightarrow \pi_0(Y^H,\hat{o})$ is defined by sending
$g\in C(H)$ to the class of $g\cdot\hat{o}$ in $\pi_0(Y^H,\hat{o})$,
and the map $\pi_0(Y^H,\hat{o})\rightarrow\pi_0^{(H,\rho)}
(X,\underline{o})$ is defined by sending the class of $y\in Y^H$
to the class of $u_y\in [(S^0(H),0);(X,\underline{o})]_\rho$, where
$u_y$ is the map defined by $(f,\rho):(S^0,H)\rightarrow (Y,G)$
where $f(0)=\hat{o}$, $f(1)=y$. The exactness of $(2.2.10)$ at
$C(H)$ is trivial, we shall focus on the remaining case at
$\pi_0(Y^H,\hat{o})$. First of all, suppose $y\in Y^H$ is path-connected
to $g\cdot \hat{o}$ in $Y^H$ for some $g\in C(H)$ by a path $\gamma$.
We set $z=g^{-1}\cdot y$. Then it is easily seen that $u_y=u_z$,
which is homotopic to the base point in $\pi_0^{(H,\rho)}(X,\underline{o})$
via the guided path defined by $(g^{-1}\circ\gamma,\rho)$. On the other
hand, suppose $u_y$ is homotopic to the base point in $\pi_0^{(H,\rho)}
(X,\underline{o})$. Then the homotopy defines an element $\Phi\in
[(I(H),0);(X,\underline{o})]_\rho$ such that the restriction of $\Phi$
to $S^0=\partial I$ is $u_y$. We apply the construction in Lemma 2.2.2
to represent $\Phi$ by a pair $(\gamma,\rho):(I,H)\rightarrow (Y,G)$ such
that $\gamma(0)=\hat{o}$. Then there is a $g\in G$ such that $y=g\cdot
\gamma(1)$ and $\rho=Ad(g)^{-1}\circ\rho$. The last equation implies that
$g\in C(H)$. Moreover, $y$ is path-connected to $g\cdot\hat{o}$ in $Y^H$
through $g\circ\gamma$. Hence $(2.2.10)$ is exact at $\pi_0(Y^H,\hat{o})$.

\hfill $\Box$

Now we turn our attention to the natural homomorphisms
$$
C:\pi_1^{(G,\rho)}(X,\underline{o})\rightarrow
\mbox{Aut}(\pi_k^{(G,\rho)}(X,\underline{o})),\; k\geq 1,
\leqno (2.2.11)
$$
and
$$
C^\prime:\pi_1^{(G,\rho)}(A,\underline{o}|_A)\rightarrow\mbox{Aut}
(\pi_{k+1}^{(G,\rho)}(X,A,\underline{o})), \; k\geq 1.
\leqno (2.2.12)
$$
The homomorphism in $(2.2.11)$ is defined as follows. Let $u_0\in
[(S^1(G),\ast);(X,\underline{o})]_\rho$ be any guided loop. There
is a continuous map $F_{u_0}:[(S^1(G),\ast);(X,\underline{o})]_\rho
\rightarrow [(S^1(G),\ast);(X,\underline{o})]_\rho$ defined by
$u\mapsto (u_0\# u)\#\nu(u_0)$. Moreover, if $u_s$, $s\in [0,1]$,
is a path in $[(S^1(G),\ast);(X,\underline{o})]_\rho$, then the map
$F:[(S^1(G),\ast);(X,\underline{o})]_\rho\times [0,1]\rightarrow
[(S^1(G),\ast);(X,\underline{o})]_\rho$, where $F(u,s)=F_{u_s}(u)$,
is continuous. We simply define $(2.2.11)$ by setting
$$
C([u])=(F_u)_\ast:\pi_{k-1}([(S^1(G),\ast);
(X,\underline{o})]_\rho,\tilde{o})
\rightarrow\pi_{k-1}([(S^1(G),\ast);(X,\underline{o})]_\rho,
F_u(\tilde{o}))
$$
and identifying $\pi_{k-1}([(S^1(G),\ast);(X,\underline{o})]_\rho,
F_u(\tilde{o}))$ with $\pi_{k-1}([(S^1(G),\ast);
(X,\underline{o})]_\rho,\tilde{o})$ by a canonically
chosen path between $F_u(\tilde{o})$ and $\tilde{o}$ (cf. Lemma 2.1.2).
The map $C([u])$ is a homomorphism because $([(S^1(G),\ast);
(X,\underline{o})]_\rho,\tilde{o})$ is an $H$-group so that the
multiplication in $\pi_{k-1}([(S^1(G),\ast);(X,\underline{o})]_\rho,
\tilde{o})$ is also given by the homotopy associative multiplication
$\#$. On the other hand, it is easy to see that $C([u])\circ C([u^\prime])
=C([u\# u^\prime])$ and $C([\tilde{o}])=Id$. Hence $C$ is a homomorphism
from $\pi_1^{(G,\rho)}(X,\underline{o})$ into $\mbox{Aut}
(\pi_{k}^{(G,\rho)}(X,\underline{o}))$. By nature of construction,
the homomorphism in $(2.2.11)$ is natural with respect to the
homomorphisms $\Phi_\ast$ and $\iota^\ast$ in Proposition 2.2.3 (2).

The homomorphism in $(2.2.12)$ is defined as follows. Let $u\in [(S^1(G),\ast);
(A,\underline{o}|_A)]_\rho$ be any guided loop. We consider a family
of elements $u(s)\in [(I(G),S^0(G),0);(X,A,\underline{o})]_\rho$, $s\in
[0,1]$, which is defined by restricting $u\#\tilde{o}$ to $[0,s]$
and then reparametrizing it by $t\mapsto ts,\forall t\in I$. Note
that $u(0)=\tilde{o}$ and $u(1)=u\#\tilde{o}$. With $u(s)$, a
continuous map $F_u$ from the loop space $\Omega([(I(G),S^0(G),0);
(X,A,\underline{o})]_\rho,\tilde{o})$ to itself is defined as follows.
Given any loop $v:[0,1]\rightarrow [(I(G),S^0(G),0);(X,A,\underline{o})]_\rho$
where $v(0)=v(1)=\tilde{o}$, we define the loop $F_u(v)$ by
$$
F_u(v)(s)=\left\{\begin{array}{ll}
u(3s) & 0\leq s\leq\frac{1}{3}\\
u\# v(3s-1) &
\frac{1}{3}\leq s\leq \frac{2}{3}\\
u(3-3s) & \frac{2}{3}\leq s\leq 1.\\
\end{array} \right. \leqno (2.2.13)
$$
The map $F_u$ is clearly continuous. Let $v_0$ be the constant
loop into $\tilde{o}$, and set $v_0^\prime=F_u(v_0)$. Then $F_u$
induces a homomorphism
$$
C^\prime_u:\pi_k(\Omega([(I(G),S^0(G),0);(X,A,\underline{o})]_\rho,\tilde{o}),
v_0)\rightarrow\pi_k(\Omega([(I(G),S^0(G),0);(X,A,\underline{o})]_\rho,
\tilde{o}),v_0^\prime)
$$
for $k\geq 0$. On the other hand, $v_0^\prime=F_u(v_0)$ is canonically
homotopic to $v_0$, so that $C^\prime_u$ gives rise to a homomorphism
$C^\prime(u):\pi_{k+1}^{(G,\rho)}(X,A,\underline{o})\rightarrow
\pi_{k+1}^{(G,\rho)}(X,A,\underline{o})$ for $k\geq 1$. Finally, we observe
that $C^\prime(u)$ depends only on the homotopy class of $u$, and that
$C^\prime(u)\circ C^\prime(u^\prime)=C^\prime(u\# u^\prime)$,
$C^\prime(\tilde{o})=Id$. Hence the homomorphism in $(2.2.12)$.

We end this subsection with

\vspace{1.5mm}

\noindent{\bf Proof of Proposition 1.3}

\vspace{1.5mm}

(1) Given any guided path $u\in [(I(G),0,1);(X,\underline{o_1},
\underline{o_2})]_{(\rho,\eta)}$, the mapping $v\mapsto (u\# v)\#\nu(u)$
which is from $([(S^1(G),\ast);(X,\underline{o_2})]_\eta,\tilde{o_2})$
to $([(S^1(G),\ast);(X,\underline{o_1})]_\rho,
(u\#\tilde{o_2})\#\nu(u))$ is continuous. The induced homomorphism
between the homotopy groups is defined to be $u_\ast$, where we
take a canonical path between $\tilde{o_1}$ and $(u\#\tilde{o_2})\#\nu(u)$
to identify the homotopy groups of $([(S^1(G),\ast);(X,\underline{o_1})]_\rho,
\tilde{o_1})$ with the corresponding ones of $([(S^1(G),\ast);
(X,\underline{o_1})]_\rho,(u\#\tilde{o_2})\#\nu(u))$. By the
nature of definition, for any guided paths $u_1, u_2$, we have
$\nu(u_2)_\ast\circ (u_1)_\ast=C([v(u_2)\# u_1])$, from which it
follows that $u_\ast$ is isomorphic.

\vspace{1.5mm}

(2) Given any guided path $u\in [(I(G),0,1);(A,\underline{o_1}|_A,
\underline{o_2}|_A)]_{(\rho,\eta)}$, we define a continuous map $F_u:
v\mapsto u\# v$ from $[(I(G),S^0(G),0);(X,A,\underline{o_2})]_\eta$ to
$[(I(G),S^0(G),0);(X,A,\underline{o_1})]_\rho$, sending the base
point $\tilde{o_2}$ to $u\#\tilde{o_2}$, which is canonically
path-connected to the base point $\tilde{o_1}$ as follows.
Consider the restriction of $u\#\tilde{o_2}$ to $[0,s]$, and then
reparametrize it by $t\mapsto ts,\forall t\in I$. We obtain a path
$u(s)$, $s\in [0,1]$, in $[(I(G),S^0(G),0);(X,A,\underline{o_1})]_\rho$
which satisfies $u(0)=\tilde{o_1}$ and $u(1)=u\#\tilde{o_2}$.
With the canonical isomorphisms provided by the path $u(s)$, the
map $F_u$ induces a map $u_\ast$ from the homotopy groups of $([(I(G),
S^0(G),0);(X,A,\underline{o_2})]_\eta,\tilde{o_2})$ to the
corresponding ones of $([(I(G),S^0(G),0);(X,A,\underline{o_1})]_\rho,
\tilde{o_1})$, which is a homomorphism for $k\geq 1$ and base
point preserving for $k=0$. By the nature of definition, for any
guided paths $u_1,u_2$, we have $\nu(u_2)_\ast\circ (u_1)_\ast=
C^\prime([\nu(u_2)\# u_1])$ when $k\geq 1$, from which it follows
that $u_\ast$ is isomorphic. When $k=0$, one can easily check that
$\nu(u)_\ast\circ (u)_\ast$ is the identity map. Hence $u_\ast$ is
a bijection when $k=0$.

\vspace{1.5mm}

(3) The existence of natural mappings $\iota^\ast$ in Proposition
2.2.3 (2) implies that as far as the isomorphism class of
$\pi_k^{(G,\rho)}(X,\underline{o})$ or $\pi_{k+1}^{(G,\rho)}(X,A,
\underline{o})$ is concerned, one may always assume that $G$ is a
subgroup of $G_{\hat{o}}$ and $\rho$ is the inclusion $G\subset
G_{\hat{o}}$. In order to see that for a different but conjugate
subgroup, there is a canonical isomorphism intervening, we
consider the following more general situation: Suppose $\underline{o_1}$,
$\underline{o_2}$ are base-point structures such that $o_1=o_2=o$, and
$G_1\subset G_{\hat{o_1}}$ and $G_2\subset G_{\hat{o_2}}$ such that
there is a $\xi\in T(U_{o_1},U_{o_2})$ satisfying $\lambda_\xi(G_1)=G_2$.
Then there is a guided path $u_\xi\in [(I(G_1),0,1);(X,\underline{o_1},
\underline{o_2})]_{(\rho_1,\rho_2\circ\lambda_\xi)}$, or in the
second case, in $[(I(G_1),0,1);(A,\underline{o_1}|_A,
\underline{o_2}|_A)]_{(\rho_1,\rho_2\circ\lambda_\xi)}$, which is
defined by $(\{f_0,f_1\},\{\xi\})$ where $f_0([0,\frac{2}{3}))=\hat{o_1}$
and $f_1((\frac{1}{3},1])=\hat{o_2}$. The isomorphism $(u_\xi)_\ast$
in (1) or (2) will then do.

\hfill $\Box$

\subsection{Exact sequence of a pair}

Let $(X,A,\underline{o})$ be any pair, and $\rho:G\rightarrow G_{\hat{o}}$
be any injective homomorphism. Denote by $i:[(S^k(G),\ast);
(A,\underline{o}|_A)]_\rho\rightarrow ([(S^k(G),\ast);(X,\underline{o})]_\rho$
the mapping induced by the inclusion $(A,\underline{o}|_A)\subset
(X,\underline{o})$, and by $j:([(S^1(G),\ast);(X,\underline{o})]_\rho,\tilde{o})
\rightarrow ([(I(G),S^0(G),0);(X,A,\underline{o})]_\rho,\tilde{o})$ the
continuous map induced by the mapping $[(S^1(G),\ast);(X,\underline{o})]_\rho
\rightarrow [(I(G),0);(X,\underline{o})]_\rho$ of forgetting the second base-point
structure. Denote by $i_\ast$ (for the case when $k=1$), $j_\ast$ the
corresponding homomorphisms between the homotopy groups. Then we have

\begin{thm}
Let $(X,A,\underline{o})$ be any pair, and $\rho:G\rightarrow G_{\hat{o}}$
be any injective homomorphism.
\begin{itemize}
\item [{(1)}] For any $k\geq 0$, there exists a mapping
$$
\partial:\pi_{k+1}^{(G,\rho)}(X,A,\underline{o})\rightarrow
\pi_k^{(G,\rho)}(A,\underline{o}|_A) \leqno (2.3.1)
$$
which is a homomorphism for $k\geq 1$, and base point preserving for
$k=0$. Moreover, for any $z\in\pi_1^{(G,\rho)}(A,\underline{o}|_A)$,
$\partial\circ C^\prime(z)=C(z)\circ\partial$ holds where $C$ and
$C^\prime$ are given in $(2.2.11)$ and $(2.2.12)$, and for any map
$\Phi$ between pairs, we have $\partial\circ\Phi_\ast=(\Phi|_A)_\ast
\circ\partial$.
\item [{(2)}] There exists a natural long exact sequence
$$
\begin{array}{l}
\cdots\rightarrow \pi_{k+1}^{(G,\rho)}(X,A,\underline{o})
\stackrel{\partial}{\rightarrow}\pi_{k}^{(G,\rho)}(A,\underline{o}|_A)
\stackrel{{i}_{\ast}}{\rightarrow}\pi_k^{(G,\rho)}(X,\underline{o})
\stackrel{{j}_{\ast}}{\rightarrow}\pi_k^{(G,\rho)}(X,A,
\underline{o})\stackrel{\partial}{\rightarrow}\cdots \\
\stackrel{{j}_{\ast}}{\rightarrow}\pi_1^{(G,\rho)}(X,A,
\underline{o})\stackrel{\partial}{\rightarrow}\pi_0^{(G,\rho)}
(A,\underline{o}|_A)\stackrel{{i}_{\ast}}{\rightarrow}
\pi_0^{(G,\rho)}(X,\underline{o}).
\end{array}
\leqno (2.3.2)
$$
\end{itemize}
\end{thm}

\pf
(1) The case $k=0$. Given any guided relative loop $u\in
[(I(G),S^0(G),0);(X,A,\underline{o})]_\rho$, the restriction $u|_A\in
[(S^0(G),0);(A,\underline{o}|_A)]_\rho$, such that if $u_1,u_2$ are
path-connected, so are $(u_1)|_A,(u_2)|_A$. We define $\partial [u]=
[u|_A]$, which is clearly base point preserving and satisfies $\partial
\circ\Phi_\ast=(\Phi|_A)_\ast\circ\partial$ for any map $\Phi$ between
pairs.

To define $(2.3.1)$ for $k\geq 1$, we regard $(S^k,\ast)$ as the suspension
$(SS^{k-1},\ast)$. Then for any continuous map $u:(S^k,\ast)\rightarrow
([(I(G),S^0(G),0);(X,A,\underline{o})]_\rho,\tilde{o})$, we consider
the corresponding map $\Phi_u\in [(CS^k(G),S^k(G),\ast);(X,A,
\underline{o})]_\rho$ constructed in Lemma 2.2.1 (2). We define $\partial
([u])=[(\Phi_u)|_{S^k(G)}]\in\pi_k^{(G,\rho)}(A,\underline{o}|_A)$, which
is also the homotopy class of the continuous map $\partial u=
u_{(\Phi_u)|_{S^k(G)}}:(S^{k-1},\ast)\rightarrow ([(S^1(G),\ast);
(A,\underline{o}|_A)]_\rho,\widetilde{o|_A})$ constructed in Lemma 2.2.1 (1).
To see that $(2.3.1)$ is a homomorphism, we simply observe that $(SS^{k-1},\ast)$
is an $H$-cogroup (cf. \cite{Sw}) and the $H$-cogroup structure of
$(SS^{k-1},\ast)$, which defines the group structure of
$\pi_{k+1}^{(G,\rho)}(X,A,\underline{o})$, corresponds to the $H$-group
structure of $([(S^1(G),\ast);(A,\underline{O}|_A)]_\rho,\widetilde{o|_A})$
under the correspondence $u\mapsto \partial{u}$. Finally, it is
straightforward from the definitions that $\partial\circ C^\prime(z)
=C(z)\circ\partial$ holds for any $z\in\pi_1^{(G,\rho)}(A,\underline{o}|_A)$,
and $\partial\circ\Phi_\ast=(\Phi|_A)_\ast\circ\partial$ holds for any map
$\Phi$ between pairs.

\vspace{1.5mm}

(2) We begin by showing that the composition of any two consecutive
homomorphisms of $(2.3.2)$ is zero. First of all, it is straightforward
from the definition that $\partial\circ {j}_\ast=0$. To see
${j}_\ast\circ {i}_\ast=0$, suppose $u:(S^k,\ast)\rightarrow
([(S^1(G),\ast);(A,\underline{o}|_A)]_\rho,\widetilde{o|_A})$, $k\geq 0$,
is any continuous map. We shall prove ${j}_\ast\circ i_\ast([u])=0$.
We take a homomorphism $\tau(\theta,\{O_a\},\{\O_\alpha\})$ constructed
in the proof of Lemma 2.2.1 (1), which represents the map $\Phi_u\in
[(SS^k(G),\ast);(A,\underline{o}|_A)]$ that corresponds to $u$ in
Lemma 2.2.1 (1). Now for any $s\in [0,1]$, denote by $\kappa_s$ the
homomorphism which is obtained by rescaling the restriction of
$\tau(\theta,\{O_a\},\{\O_\alpha\})$ to $\{[t,x]|t\leq s,x\in S^k\}$
under $t\mapsto ts$. Each $\kappa_s$ defines a continuous map
$u_s:(S^k,\ast)\rightarrow ([(I(G),S^0(G),0);(X,A,\underline{o})]_\rho,
\tilde{o})$ such that $s\mapsto u_s$ is continuous. Moreover, we have
$u_1=j\circ i\circ u$ and $u_0(S^k)=\tilde{o}$. Hence $j_\ast\circ
i_\ast([u])=0$. To see $i_\ast\circ\partial=0$, we first look at the case
$k\geq 1$. For any $u:(S^k,\ast)\rightarrow ([(I(G),S^0(G),0);
(X,A,\underline{o})]_\rho,\tilde{o})$, we consider the map
$\Phi_u\in [(CS^k(G),S^k(G),\ast);(X,A,\underline{o})]_\rho$ constructed
in Lemma 2.2.1 (2). We identify $(CS^k,\ast)$ with $(SCS^{k-1},\ast)$,
and then apply Lemma 2.2.1 (1) to obtain a continuous map
$H=u_{\Phi_u}:(CS^{k-1},\ast)\rightarrow ([(S^1(G),\ast);
(X,\underline{o})]_\rho,\tilde{o})$. Now observe $H|_{(S^{k-1},\ast)}=
i\circ (\partial u)$. Hence $i_\ast\circ\partial ([u])=0$ for any
continuous map $u:(S^k,\ast)\rightarrow ([(I(G),S^0(G),0);
(X,A,\underline{o})]_\rho,\tilde{o})$. The case $k=0$ is similar,
where we start with a guided relative loop $u\in [(I(G),S^0(G),0);
(X,A,\underline{o})]_\rho$ instead of the map $u$, and replace $\Phi_u$
by the corresponding path $u\in [(I(G),0);(X,\underline{o})]_\rho$ in the
argument.

The exactness at $\pi_0^{(G,\rho)}(A,\underline{o}|_A)$: Let
$u\in [(S^0(G),0),(A,\underline{o}|_A)]_\rho$ be any element such that
$i_\ast([u])=[\tilde{o}]$. Then there is a homotopy between $i(u)$ and
$\tilde{o}$, which may be interpreted as a guided path $\Phi\in
[(I(G),0);(X,\underline{o})]_\rho$. Now $\Phi|_{S^0(G)}=u$ so that
$\Phi\in [(I(G),S^0(G),0);(X,A,\underline{o})]_\rho$. Clearly
$\partial([\Phi])=[u]$. Hence $(2.3.2)$ is exact at $\pi_0^{(G,\rho)}
(A,\underline{o}|_A)$.

The exactness at $\pi_1^{(G,\rho)}(X,A,\underline{o})$: Let
$u\in [(I(G),S^0(G),0);(X,A,\underline{o})]_\rho$ be any guided relative
loop such that $\partial([u])=[\widetilde{o|_A}]$. Then $u|_{(S^0(G),0)}$
is homotopic to $\widetilde{o|_A}$ in $[(S^0(G),0);(A,\underline{o}|_A)]
_\rho$, which gives rise to an element $v\in [(I(G),0);(A,
\underline{o}|_A)]_\rho$ satisfying $v|_{(S^0(G),0)}=u|_{(S^0(G),0)}$.
We may `compose' $u$ with $i(\nu(v))$ to obtain a guided loop
$u\# i(\nu(v))\in [(S^1(G),\ast);(X,\underline{o})]_\rho$ as follows.
We pick a homomorphism $\{f_i\}\in\O_{(\{\xi_{ji}\},V_\alpha)}$
representing $u$ and pick a homomorphism $\{h_k\}\in\O_{\{\eta_{lk}\}}$
representing $\nu(v)$ where each $\eta_{lk}\in T(V_{\alpha(k)},V_{\alpha(l)})$
for some correspondence $k\mapsto \alpha=\alpha(k)$. We compose
$(\{f_i\},\{\xi_{ji}\})$ with $(\{h_k\},\{\rho_{\alpha(l)\alpha(k)}
(\eta_{lk})\})$ to obtain a homomorphism $\tau$, and define
$u\# i(\nu(v))=[\tau]$. Note that $[\tau]$ may not be uniquely determined
by $u$ and $v$. But because of its specific form, $\tau$ provides
a natural homotopy between $j(u\# i(\nu(v)))$ and $u$. Hence
$j_\ast([u\# i(\nu(v))])=[u]$, and $(2.3.2)$ is exact at
$\pi_1^{(G,\rho)}(X,A,\underline{o})$.

The exactness at $\pi_k^{(G,\rho)}(X,\underline{o})$ with $k\geq 1$:
Let $u:(S^{k-1},\ast)\rightarrow ([(S^1(G),\ast);(X,\underline{o})]_\rho,
\tilde{o})$ be any continuous map such that $[u]\in\ker\;j_\ast$.
Then there is a homotopy $H:(CS^{k-1},\ast)\mapsto ([(I(G),S^0(G),0);
(X,A,\underline{o})]_\rho,\tilde{o})$ with $H|_{(S^{k-1},\ast)}=j\circ u$.
As shown in Lemma 2.2.1 (2), $H$ determines an element $\Phi_H\in
[(CCS^{k-1}(G),CS^{k-1}(G),\ast);(X,A,\underline{o})]_\rho$, which
is the equivalence class of a canonically constructed homomorphism
$\tau=\tau(\theta,\{O_a\},\{\O_\alpha\})$. Let $[t_1,t_2,x],t_1\in I, t_2\in I,
x\in S^{k-1}$ be the coordinates of $CCS^{k-1}$. Then we observe that
the restriction of $\tau$ to $\{[t_1,1,x]|t_1\in I, x\in S^{k-1}\}$
represents the element $\Phi_u\in [(SS^{k-1}(G),\ast);(X,\underline{o})]_\rho$
determined by $u$ in Lemma 2.2.1 (1), and the restriction of $\tau$ to
$\{[1,t_2,x]|t_2\in I,x\in S^{k-1}\}$ represents $i((\Phi_H)|
_{(CS^{k-1}(G),\ast)})\in [(CS^{k-1}(G),\ast);(X,\underline{o})]_\rho$.
Moreover, $(\Phi_H)|_{(CS^{k-1}(G),\ast)}$ is in fact defined over
$(SS^{k-1}(G),\ast)$ because its restriction to $(S^{k-1}(G),\ast)$
is $(\Phi_{j\circ u})|_{(S^{k-1}(G),\ast)}$, where $\Phi_{j\circ u}$
is the element in $[(CS^{k-1}(G),S^{k-1}(G),\ast);(X,A,\underline{o})]_\rho$
that is constructed in Lemma 2.2.1 (2) for $j\circ u$.
$(\Phi_{j\circ u})|_{(S^{k-1}(G),\ast)}$ is obviously a constant map.
Let $u_{(\Phi_H)|_{(CS^{k-1}(G),\ast)}}$, $u_{i((\Phi_H)|_{(CS^{k-1}
(G),\ast)})}$ be the continuous maps defined in Lemma 2.2.1 (1)
which correspond to $(\Phi_H)|_{(CS^{k-1}(G),\ast)}$ and $i((\Phi_H)|
_{(CS^{k-1}(G),\ast)})$ respectively. Then the map $x\mapsto u(x)\#
\nu(u_{i((\Phi_H)|_{(CS^{k-1}(G),\ast)})}(x))$, $x\in S^{k-1}$, is
homotopic to the constant map $x\mapsto\tilde{o}$ via a homotopy $F(s):
(S^{k-1},\ast)\rightarrow ([(S^1(G),\ast);(X,\underline{o})]_\rho,
\tilde{o})$, where $F(s)$ is defined by the restriction of $\tau$ to
$$
\{[t_1,s,x]|t_1\leq s,x\in S^{k-1}\}\cup \{[s,t_2,x]|
t_2\leq s,x\in S^{k-1}\}.
$$
This implies that $[u]=[u_{i((\Phi_H)|_{(CS^{k-1}(G),\ast)})}]=
[i(u_{(\Phi_H)|_{(CS^{k-1}(G),\ast)}})]=i_\ast(
[u_{(\Phi_H)|_{(CS^{k-1}(G),\ast)}}])$. Hence $(2.3.2)$ is exact at
$\pi_k^{(G,\rho)}(X,\underline{o})$.

The exactness at $\pi_k^{(G,\rho)}(A,\underline{o}|_A)$ with $k\geq 1$:
Let $u:(S^{k-1},\ast)\rightarrow ([(S^1(G),\ast);(A,\underline{o}|_A)]_\rho,
\widetilde{o|_A})$ be any continuous map such that $[u]\in \ker\;i_\ast$.
Then there is a homotopy $H:(CS^{k-1},\ast)\rightarrow ([(S^1(G),\ast);
(X,\underline{o})]_\rho,\tilde{o})$ with $H|_{(S^{k-1},\ast)}=i\circ u$.
We shall repeat the construction in Lemma 2.2.1 (1) to the map
$H$, but for the specific purpose here, we choose the cover $\{\O_\alpha\}$
in the construction of the homomorphism $\tau(\theta,\{O_a\},\{\O_\alpha\})$
as follows. First of all, we fix a homomorphism $\sigma=(\{i_\alpha\},
\{\rho_{\beta\alpha}\})$ which represents the subspace inclusion
$(A,\underline{o}|_A)\subset (X,\underline{o})$. Then we cover
$u(S^{k-1})$ by finitely many $\{\O_s\}$, where each $\O_s$ has the form
$\O_{\{\delta_{ji,s}\}}$. Consequently, $H(S^{k-1})=i\circ u(S^{k-1})$
is covered by $\{\sigma(\O_s)\}$, where $\sigma(\O_s)=\O_{\{\rho_{\alpha(s,j)
\alpha(s,i)}(\delta_{ji,s})\}}$. Now we cover $H(CS^{k-1})$ by finitely many
$\{\O_k\}$ such that $\{\sigma(\O_s)\}\subset \{\O_k\}$ and covers
$H(S^{k-1})$. We proceed to construct the homomorphism
$\tau(\theta,\{O_a\},\{\O_k\})$, which defines the element $\Phi_H
\in [(SCS^{k-1}(G),\ast);(X,\underline{o})]_\rho$ associated to $H$
in Lemma 2.2.1 (1). Now we identify $(SCS^{k-1},\ast)$ with
$(CSS^{k-1},\ast)=(CS^k,\ast)$. Then the special choice of the
cover $\{\O_k\}$ in $\tau(\theta,\{O_a\},\{\O_k\})$ implies that
$\tau(\theta,\{O_a\},\{\O_k\})$ also defines an element $\Psi\in
[(CS^k(G),S^k(G),\ast);(X,A,\underline{o})]_\rho$, which is associated
with a continuous map $u_\Psi:(S^k,\ast)\rightarrow
([(I(G),S^0(G),0);(X,A,\underline{o})]_\rho,\tilde{o})$ in Lemma 2.2.1 (2),
such that $\partial u_\Psi:(S^{k-1},\ast)\rightarrow
([(S^1(G),\ast);(A,\underline{o}|_A)]_\rho,\widetilde{o|_A})$ equals $u$.
Hence $[u]=\partial([u_\Psi])$ and $(2.3.2)$ is exact at
$\pi_k^{(G,\rho)}(A,\underline{o}|_A)$.

The exactness at $\pi_{k+1}^{(G,\rho)}(X,A,\underline{o})$, $k\geq 1$:
Let $u:(S^k,\ast)\rightarrow ([(I(G),S^0(G),0);
(X,A,\underline{o})]_\rho,\tilde{o})$ be any continuous map such that
$[u]\in \ker\;\partial$. Then there is a homotopy $H:(CS^{k-1},\ast)
\rightarrow ([(S^1(G),\ast);(A,\underline{o}|_A)]_\rho,\widetilde{o|_A})$
with $H|_{(S^{k-1},\ast)}=\partial u$. We fix a homomorphism $\sigma=
(\{i_\alpha\},\{\rho_{\beta\alpha}\})$ which represents the subspace
inclusion $(A,\underline{o}|_A)\subset (X,\underline{o})$, and cover
$H(CS^{k-1})$ by finitely many $\{\O_s\}$ where each $\O_s
=\O_{\{\xi_{ji,s}\}}$. Then $i\circ H(CS^{k-1})$ is covered
by $\{\sigma(\O_s)\}$ where $\sigma(\O_s)=\O_{\{\rho_{\alpha(j,s)
\alpha(i,s)}(\xi_{ji,s})\}}$. Now we construct a homomorphism
$\tau_1=\tau(\theta,\{O_a\},\{\sigma(\O_s)\})$ which represents the
element $\Phi_{i\circ H}\in [(SCS^{k-1}(G),\ast);(X,\underline{o})]_\rho$
associated to $i\circ H$ in Lemma 2.2.1 (1), and construct a homomorphism
$\tau_2=\tau(\theta^\prime,\{O_{a^\prime}\},\{\O_{\alpha^\prime}\})$
which represents the element $\Phi_u\in [(CS^k(G),S^k(G),\ast);
(X,A,\underline{o})]_\rho$ associated to $u$ in Lemma 2.2.1 (2). By
identifying $(SCS^{k-1},\ast)$ with $(CSS^{k-1},\ast)=(CS^k,\ast)$,
we regard $\tau_1$ as a homomorphism which defines a map from $(CS^k(G),\ast)$
to $(X,\underline{o})$. Since the restriction of $\tau_1$ to the subspace
$(S^k(G),\ast)$ represents $\Phi_{i\circ\partial u}=i\circ (\Phi_u)|_{(S^k,\ast)}$,
we may arrange to join $\tau_2$ with $\tau_1$ along $(S^k(G),\ast)$ to define
a homomorphism $\tau$ where $[\tau]\in [(SS^k(G),\ast);(X,\underline{o})]_\rho$,
where we identify $(SS^k(G),\ast)$ with $(CS^k(G),\ast)\bigcup_{(S^k(G),\ast)}
(CS^k(G),\ast)$. We then apply Lemma 2.2.1 (1) to $[\tau]$ to obtain
a continuous map $u_{[\tau]}:(S^k,\ast)\rightarrow ([(S^1(G),\ast);
(X,\underline{o})]_\rho,\tilde{o})$. The specific form of $\tau$
gives rise to a natural homotopy between $j\circ u_{[\tau]}$ and $u$.
Hence $[u]=j_\ast([u_{[\tau]}])$ and $(2.3.2)$ is exact at
$\pi_{k+1}^{(G,\rho)}(X,A,\underline{o})$.

\vspace{1.5mm}

The proof of Theorem 2.3.1 is thus completed.

\hfill $\Box$

\subsection{The theory of coverings}

In this subsection we give a detailed, elementary presentation of
the theory of coverings of orbispaces. We also derive the following
useful corollary which was involved in the proof of the Arzela-Ascoli
precompactness theorem for $C^r$ maps between smooth orbifolds in \S 3.4
of \cite{C1}. Its proof is given at the end of this subsection.

\begin{lem}
Let $\Phi:X\rightarrow X^\prime$ be a map\footnote{here $\Phi$ is
allowed to be a general map, not restricted by {\bf Convention} in
Introduction.} between global quotients,
where $X=Y/G$, $X^\prime=Y^\prime/G^\prime$, such that $Y$ is connected
and locally path-connected. Then $\Phi$ is defined by a pair $(f,\lambda):
(Y,G)\rightarrow (Y^\prime,G^\prime)$ where $f$ is $\lambda$-equivariant
if and only if there are base-point structures $\underline{o}$,
$\underline{o^\prime}$ such that under $\Phi_\ast:\pi_1(X,\underline{o})
\rightarrow\pi_1(X^\prime,\underline{o^\prime})$, $\pi_1(Y,\hat{o})$ is
mapped into $\pi_1(Y^\prime,\hat{o}^\prime)$, where $\pi_1(Y,\hat{o})$,
$\pi_1(Y^\prime,\hat{o}^\prime)$ are regarded as subgroups of
$\pi_1(X,\underline{o})$, $\pi_1(X^\prime,\underline{o^\prime})$
via the exact sequence $(2.2.10)$ with $H=\{1\}$.
\end{lem}

The following definition generalizes the notion
`orbifold covering' in Thurston \cite{Th}.

\begin{defi}
Let $\Pi:Y\rightarrow X$ be a map of orbispaces and
$\pi: Y\rightarrow X$ be the induced map between underlying spaces.
We call $\Pi$ a covering map, if there exists a homomorphism
$(\{\pi_\alpha\},\{\rho_{\beta\alpha}\}):\Gamma\{V_\alpha\}
\rightarrow\Gamma\{U_{\alpha^\prime}\}$, $\alpha\in\Lambda$,
whose equivalence class is the map $\Pi$, such that
\begin{itemize}
\item [{(a)}] each $\pi_\alpha:\widehat{V_\alpha}\rightarrow
\widehat{U_\alpha}$ is a $\rho_\alpha$-equivariant homeomorphism,
\item [{(b)}] for any $U\in\{U_\alpha\}$, $\{V_\alpha|\alpha\in
\Lambda(U)\}$ is the set of connected components of $\pi^{-1}(U)$,
where $\Lambda(U)$ is the subset of $\Lambda$ defined by
$\Lambda(U)=\{\alpha\in\Lambda\mid U=U_\alpha\}$, and
\item [{(c)}] the map $\pi:Y\rightarrow X$ between underlying spaces
is surjective.
\end{itemize}

The orbispace $Y$ together with the covering map $\Pi:Y\rightarrow X$
is called a covering space of $X$. Each element $U\in\{U_\alpha\}$, which
all together form a cover of $X$ by ${\em (c)}$, is called an elementary
neighborhood of $X$ with respect to the covering map $\Pi$.
\end{defi}

\begin{re}{\em

\vspace{1.5mm}

(1) By passing to an induced homomorphism of
$(\{\pi_\alpha\},\{\rho_{\beta\alpha}\})$ if necessary, we may
assume that a connected open subset of an elementary
neighborhood is an elementary neighborhood.

\vspace{1.5mm}

(2) By the general assumption made in {\bf Convention} in
Introduction, each $\rho_\alpha:G_{V_\alpha}\rightarrow G_{U_\alpha}$
is injective, and each $\rho_{\beta\alpha}:T(V_\alpha,V_\beta)
\rightarrow T(U_\alpha,U_\beta)$ is partially injective in the
sense that if $\rho_{\beta\alpha}(\xi_1)=\rho_{\beta\alpha}(\xi_2)$
and $\mbox{Domain }(\phi_{\xi_1})\cap\mbox{Domain }(\phi_{\xi_2})\neq
\emptyset$, then $\xi_1=\xi_2$. But in the case of a covering map,
one can easily verify that $\rho_{\beta\alpha}$ is in fact
injective.

\vspace{1.5mm}

(3) The composition of two covering maps is a covering map.

\vspace{1.5mm}

(4) Let $\underline{q}=(q,V_o,\hat{q})$, $\underline{p}=(p,U_o,\hat{p})$
be base-point structures of $Y$ and $X$ respectively. Given any
covering map $\Pi:Y\rightarrow X$, we may always assume that
$V_o\in\{V_\alpha\}$, $U_o\in\{U_{\alpha^\prime}\}$, and $\pi_o:\widehat{V_o}
\rightarrow\widehat{U_o}$ sends $\hat{q}$ to $\hat{p}$.
}
\end{re}

First of all, we investigate the path-lifting property of a
covering map. To simplify the notation, we write $P(X,\underline{o_1},
\underline{o_2})$ for the path space $[(I,0,1);(X,\underline{o_1},
\underline{o_2})]$. Let $\Pi:(Y,\underline{q})\rightarrow
(X,\underline{p})$ be a covering map, represented by $(\{\pi_\alpha\},
\{\rho_{\beta\alpha}\})$ as described in Definition 2.4.2. Given any
$\underline{p^\prime}=(p^\prime,U,\hat{p}^\prime)$ where $U\in
\{U_{\alpha^\prime}\}$, and any connected component $V_\alpha$ of
$\pi^{-1}(U)$, $(\{\pi_\alpha\},\{\rho_{\beta\alpha}\})$ canonically
determines a $\underline{q^\prime}=(q^\prime, V_{\alpha},\hat{q}^\prime)$
by setting $\hat{q}^\prime=\pi^{-1}_\alpha(\hat{p}^\prime)$ and
$q^\prime=\pi_{V_{\alpha}}(\hat{q}^\prime)\in Y$. For any $g\in
G_{U}$ and $\underline{p^\prime}=(p^\prime,U,\hat{p}^\prime)$, we
denote the base-point structure $(p^\prime,U,g\cdot\hat{p}^\prime)$ by
$g\cdot\underline{p}^\prime$. Note that there is a natural mapping
from $P(X,\underline{p},\underline{p^\prime})$ to $P(X,\underline{p},
g\cdot\underline{p^\prime})$, denoted by $u\mapsto g\cdot u$, which is
defined by sending a representative $(\{f_i\},\{\xi_{ji}\})$,
$i=0,1,\cdots,n$, of $u$ to $(\{f_i^\prime\},\{\xi_{ji}^\prime\})$,
$i=0,1,\cdots,n$, where $f_i^\prime=f_i$ for $i\leq n-1$,
$\xi_{(i+1)i}^\prime=\xi_{(i+1)i}$ for $i\leq n-2$, and $f_n^\prime=
g\circ f_n$, $\xi_{n(n-1)}^\prime=g\circ\xi_{n(n-1)}$.

\begin{lem}
Let $\Pi:(Y,\underline{q})\rightarrow (X,\underline{p})$ be
a covering map defined by $(\{\pi_\alpha\},\{\rho_{\beta\alpha}\})$,
and $\underline{p^\prime}=(p^\prime,U,\hat{p}^\prime)$ where $U\in
\{U_{\alpha^\prime}\}$. For any path $u\in P(X,\underline{p},
\underline{p^\prime})$, there is a unique component $V_{\alpha(u)}$ of
$\pi^{-1}(U)$ and a unique left coset $\delta(u)\in G_{U}/
\rho_{\alpha(u)}(G_{V_{\alpha(u)}})$, such that after fixing a
representative $g\in G_{U}$ of $\delta(u)$, there is a unique path
$\ell(u)_g\in P(Y,\underline{q},\underline{q^\prime}(g))$ where
$\underline{q^\prime}(g)$ is the base-point structure canonically
associated to $g\cdot\underline{p}^\prime$ and $V_{\alpha(u)}$ by
$(\{\pi_\alpha\},\{\rho_{\beta\alpha}\})$, such that
$(\{\pi_\alpha\},\{\rho_{\beta\alpha}\})\circ\ell(u)_g=g\cdot u$.
For any continuous family $u_s\in P(X,\underline{p},
\underline{p^\prime})$, the component $V_{\alpha(u_s)}$
and the coset $\delta(u_s)$ are locally constant with respect
to $s$. When the parameter space of $s$ is connected, one can choose
representatives of $\delta(u_s)$ independent of $s$, and in this case,
the family of paths $\ell(u_s)_g$ is continuous in $s$, where $g$
is any such a representative of $\delta(u_s)$.
\end{lem}

\pf
Given any $u\in P(X,\underline{p},\underline{p^\prime})$,
we pick a representative $\sigma=(\{f_i\},\{\xi_{ji}\}):\Gamma\{I_i\}
\rightarrow\Gamma\{U_{i^\prime}\}$, $i=0,1,\cdots,n$, where each
$U_{i}\in\{U_{i^\prime}\}$ is an elementary neighborhood of $X$ with
respect to the covering map $\Pi$. We shall construct a
$\ell(\sigma)=(\{\ell(f_i)\},
\{\ell(\xi_{ji})\}):\Gamma\{I_i\}\rightarrow\Gamma\{V_{i^\prime}\}$,
such that (1) each $V_i\in\{V_{i^\prime}\}$ is a connected component
of $\pi^{-1}(U_i)$, and (2) there exists $\{g_i\in G_{U_i}\mid g_0=1,i=0,1,
\cdots,n\}$ satisfying
$$
\pi_i\circ\ell(f_i)=g_i\circ f_i, \;
\rho_{ji}(\ell(\xi_{ji}))=g_j\circ\xi_{ji}\circ g_i^{-1},
\leqno (2.4.1)
$$
where $\pi_i:\widehat{V_i}\rightarrow \widehat{U_i}$ and $\rho_{ji}:
T(V_i,V_j)\rightarrow T(U_i,U_j)$ are given in $(\{\pi_\alpha\},
\{\rho_{\beta\alpha}\})$. $\ell(\sigma)$ is defined inductively as
follows. For $i=0$, we simply define $\ell(f_0)=\pi_o^{-1}\circ f_0:
I_0\rightarrow \widehat{V_o}$, which clearly satisfies $\ell(f_0)(0)
=\hat{q}$. Now suppose that $\ell(\xi_{k(k-1)})$ and $\ell(f_k)$ are
defined for all $0\leq k\leq i$. We shall define $\ell(\xi_{(i+1)i})$
and $\ell(f_{i+1})$ as follows. First of all, it follows from
$\pi_i\circ\ell(f_i)=g_i\circ f_i$ that $\pi(\pi_{V_i}(\ell(f_i)
(I_i\cap I_{i+1})))\subset U_{i+1}$, hence by (b) of Definition
2.4.2, there is a unique connected component $V_{i+1}$ of $\pi^{-1}(U_{i+1})$
which is determined by the condition $\pi_{V_i}(\ell(f_i)(I_i\cap I_{i+1}))
\cap V_{i+1}\neq\emptyset$. Secondly, let $V$ be the connected component
of $V_i\cap V_{i+1}$ that contains $\pi_{V_i}(\ell(f_i)(I_i\cap I_{i+1}))$,
and let $W$ be the connected component of $U_i\cap U_{i+1}$ that contains
$\pi(\pi_{V_i}(\ell(f_i)(I_i\cap I_{i+1})))=\pi_{U_i}(f_i
(I_i\cap I_{i+1}))$. Then $\rho_{(i+1)i}$ maps $T_V(V_i,V_{i+1})$
injectively into $T_W(U_i,U_{i+1})$, and $\xi_{(i+1)i}$ lies in
$T_W(U_i,U_{i+1})$. We pick a $\eta_{(i+1)i}\in T_V(V_i,V_{i+1})$
such that $\ell(f_i)(I_i\cap I_{i+1})\subset\mbox{Domain }
(\phi_{\eta_{(i+1)i}})$. Then since $\pi_i\circ\ell(f_i)=g_i\circ f_i$,
$\rho_{(i+1)i}(\eta_{(i+1)i})$ has the same domain with
$\xi_{(i+1)i}\circ g_i^{-1}$, and therefore there exists a
$g_{i+1}\in G_{U_{i+1}}$ such that $\rho_{(i+1)i}(\eta_{(i+1)i})=g_{i+1}
\circ\xi_{(i+1)i}\circ g_i^{-1}$. Finally, we define $\ell(f_{i+1})=
\pi_{i+1}^{-1}\circ (g_{i+1}\circ f_{i+1})$, $\ell(\xi_{(i+1)i})
=\eta_{(i+1)i}$. Then it is clear that $\ell(\xi_{(i+1)i})\circ
\ell(f_i)=\ell(f_{i+1})$ on $I_i\cap I_{i+1}$, and $\pi_{i+1}\circ
\ell(f_{i+1})=g_{i+1}\circ f_{i+1}$, $\rho_{(i+1)i}(\ell(\xi_{(i+1)i}))
=g_{i+1}\circ\xi_{(i+1)i}\circ g_i^{-1}$. By induction, $\ell(\sigma)$
is defined with the claimed properties.

We observe that in each step of the above construction, $V_i$ is uniquely
determined, and if $\ell^\prime(\sigma)=(\{\ell^\prime(f_i)\},
\{\ell^\prime(\xi_{ji})\})$ and $\{g_i^\prime\}$ is another choice,
then there exists a set $\{h_i\mid h_i\in G_{V_i}\}$ such that
$g_i^\prime=\rho_i(h_i)g_i$, and $\ell^\prime(f_i)=h_i\circ\ell(f_i)$,
$\ell^\prime(\xi_{ji})=h_j\circ\ell(\xi_{ji})\circ h_i^{-1}$. In
other words, $\ell^\prime(\sigma)$ is conjugate to $\ell(\sigma)$.

Let $\tau=(\{f_k\},\{\xi_{lk}\}):\Gamma\{J_k\}\rightarrow
\Gamma\{U_{k^\prime}\}$, $k=0,1,\cdots,m$, be any induced homomorphism
of $\sigma$, defined via $(\theta,\{\xi_k\},\{\jmath_k\})$ where
$\theta:\{0,1,\cdots,m\}\rightarrow \{0,1,\cdots,n\}$ satisfies
$\theta(0)=0$, $\theta(m)=n$, $\xi_k:J_k\rightarrow I_{\theta(k)}$
is the inclusion for all $k$, and $\jmath_k\in T(U_k,U_{\theta(k)})$
satisfies $\jmath_0=1$, $\jmath_m=1$, such that
$$
f_k=\jmath^{-1}_k\circ f_{\theta(k)}|_{J_k},\;
\xi_{lk}=\jmath_l^{-1}\circ\xi_{\theta(l)\theta(k)}\circ\jmath_k.
\leqno (2.4.2)
$$
For any choice of $(\ell(\sigma),\{g_i\})$, $(\ell(\tau),\{g_k\})$
with $(2.4.1)$ satisfied, if we set $\bar{\jmath_k}=g_{\theta(k)}\circ
\jmath_k\circ g_k^{-1}\in T(U_k,U_{\theta(k)})$, then it is easy
to check that the following hold:
$$
\pi_k\circ\ell(f_k))=\bar{\jmath_k}^{-1}\circ (\pi_{\theta(k)}
\circ\ell(f_{\theta(k)})), \;
\rho_{lk}(\ell(\xi_{lk}))=\bar{\jmath_l}^{-1}\circ
\rho_{\theta(l)\theta(k)}(\ell(\xi_{\theta(l)\theta(k)}))\circ
\bar{\jmath_k} \leqno (2.4.3)
$$
On the other hand, observe that $\bar{\jmath_0}=1$ and
$V_0=V_{\theta(0)}=V_o$, so that $(2.4.3)$ implies inductively
that $V_k\cap V_{\theta(k)}\neq\emptyset$, and $\bar{\jmath_k}
=\rho_{\theta(k)k}(\imath_k)$ for a unique $\imath_k\in
T(V_k,V_{\theta(k)})$ such that
$$
\ell(f_k)=\imath_k^{-1}\circ\ell(f_{\theta(k)}),\;
\ell(\xi_{lk})=\imath_l^{-1}\circ\ell(\xi_{\theta(l)\theta(k)})
\circ \imath_k. \leqno (2.4.4)
$$

With the preceding understood, we now determine $V_{\alpha(u)}$,
$\delta(u)$, and define $\ell(u)_g$ after having chosen a representative
$g$ of $\delta(u)$. First of all, observe that $V_m,V_{\theta(m)}
=V_n$ are connected components of $\pi^{-1}(U)$, so that
$V_m\cap V_{\theta(m)}\neq\emptyset$ implies $V_m=V_{\theta(m)}$.
We define $V_{\alpha(u)}=V_m$. Secondly, $\imath_m\in T(V_m,V_{\theta(m)})
=G_{V_m}$ and $\rho_m(\imath_m)=\bar{\jmath_m}=g_{\theta(m)}\circ\jmath_m
\circ g_m^{-1}=g_{\theta(m)}g_m^{-1}$ implies that the left coset of
$g_m$ in $G_U/\rho_{\alpha(u)}(G_{V_{\alpha(u)}})$ depends only on
$u$. We define $\delta(u)$ to be the coset of $g_m$ in
$G_U/\rho_{\alpha(u)}(G_{V_{\alpha(u)}})$. Finally, if we fix a
representative $g$ of $\delta(u)$ for the $g_m,g_n$ in the above
consideration, then $\bar{\jmath_m}=1$ and $\imath_m=1$, so that
$\ell(\sigma)$, $\ell(\tau)$ define the same element in
$P(Y,\underline{q},\underline{q^\prime}(g))$, which is defined to
be $\ell(u)_g$. By the nature of construction, $\Pi\circ\ell(u)_g
=g\cdot u$, where $\Pi$ is regarded as a map from $(Y,\underline{q},
\underline{q^\prime}(g))$ to $(X,\underline{p},g\cdot\underline{p^\prime})$.

As for the uniqueness of $\ell(u)_g$, suppose $u_1,u_2\in
P(Y,\underline{q},\underline{q^\prime}(g))$ are two paths such that
$\Pi\circ u_1=\Pi\circ u_2$. Then in particular, the paths in the
underlying space of $X$ coincide, and hence there exist
$\tau_1=(\{f_{i,1}\},\{\eta_{ji,1}\}):\Gamma\{I_i\}\rightarrow
\Gamma\{V_{i(1)}\}$, $\tau_2=(\{f_{i,2}\},\{\eta_{ji,2}\}):\Gamma\{I_i\}
\rightarrow\Gamma\{V_{i(2)}\}$ representing $u_1, u_2$ respectively,
where for each $i$, $V_{i(1)}, V_{i(2)}$ are connected components of
$\pi^{-1}(U_i)$ for some elementary neighborhood $U_i$. The compositions
of $\tau_1$, $\tau_2$ with $(\{\pi_\alpha\},\{\rho_{\beta\alpha}\})$
are equivalent, because they represent the same $\Pi\circ u_1=
\Pi\circ u_2$. Hence by Lemma 3.1.2 in \cite{C1}, they must be conjugate,
and furthermore, since they preserve the base-point structures, they
are actually equal, which means
$$
\pi_{i(1)}\circ f_{i,1}=\pi_{i(2)}\circ f_{i,2},\;
\rho_{j(1)i(1)}(\eta_{ji,1})=\rho_{j(2)i(2)}(\eta_{ji,2}).
\leqno (2.4.5)
$$
We shall derive from $(2.4.5)$ that $V_{i(1)}=V_{i(2)}$ for all $i$
by induction. For $i=0$, $V_{i(1)}=V_{i(2)}=V_o$, so that by $(2.4.5)$,
we obtain $f_{0,1}=f_{0,2}$ and $\eta_{10,1}$ and $\eta_{10,2}$
have the same domain. The latter then implies that $V_{i(1)}\cap V_{i(2)}
\neq\emptyset$ for $i=1$, or equivalently $V_{i(1)}=V_{i(2)}$ for $i=1$.
By induction $V_{i(1)}=V_{i(2)}$ for all $i$, with which $(2.4.5)$
implies $\tau_1=\tau_2$. Hence $u_1=u_2$ and $\ell(u)_g$ is unique.

Finally, let $u_s$ be any continuous family of paths.
For any $s_0$, when $s$ is sufficiently close to $s_0$,
$u_s$ is represented by a $\sigma_s=\{f_i(\cdot,s)\}\in
\O_{\{\xi_{ji}\}}$, where each $f_i(t,s)$ is continuous in both
variables. Thus when $s$ is sufficiently close to $s_0$,
$V_{\alpha(u_{s})}=V_{\alpha(u_{s_0})}$, and $\delta(u_{s})=
\delta(u_{s_0})$. In other words, $V_{\alpha(u_s)}$ and $\delta(u_s)$
are locally constant in $s$. When the parameter space of $s$ is connected,
$V_{\alpha(u_s)}$ and $\delta(u_s)$ are constant so that we can choose
representatives of $\delta(u_s)$ independent of $s$. The family of
liftings $\ell(u_s)_g$ is continuous in $s$ for any choice of
such a representative $g$, because it can be locally represented by
$\ell(\sigma_s)$ which is continuous in $s$.

\hfill $\Box$

As in the classical covering theory, we can deduce the following
immediate corollary.

\begin{coro}
Let $\Pi:(Y,\underline{q})\rightarrow (X,\underline{p})$ be
a covering map. Then the induced homomorphism
$\Pi_\ast:\pi_k(Y,\underline{q})\rightarrow
\pi_k(X,\underline{p})$ is isomorphic when $k\geq 2$ and injective
when $k=1$.
\end{coro}

We introduce the following definitions.

\begin{itemize}
\item [{(1)}] An orbispace is locally path-connected if the
associated \'{e}tale topological groupoid is locally
path-connected, or equivalently, each local chart $\widehat{U_i}$
is locally path-connected.
\item [{(2)}] An orbispace $X$ is semi-locally 1-connected if for any
point $p$, there is a local chart $U_i$ containing $p$, such that
the composition of homomorphisms $\pi_1(\widehat{U_i},\hat{p})
\rightarrow\pi_1(U_i,\underline{p})\rightarrow\pi_1(X,\underline{p})$
has trivial image for any base-point structure $\underline{p}=
(p,U_i,\hat{p})$.
\item [{(3)}] A connected, locally path-connected covering space
$\Pi:Y\rightarrow X$ is universal if $\pi_1(Y)$ is trivial.
\end{itemize}

The last definition is justified by the following lemma.

\begin{lem}
Let $\Pi:(Y,\underline{q})\rightarrow (X,\underline{p})$ be a
covering map. For any map\footnote{here both $\Phi,\ell(\Phi)$
are general maps, i.e., not restricted by {\bf Convention} in
Introduction.} $\Phi:(Z,\underline{z})\rightarrow
(X,\underline{p})$ where $Z$ is connected and locally
path-connected, there exists a unique map $\ell(\Phi):
(Z,\underline{z})\rightarrow (Y,\underline{q})$ such that
$\Pi\circ\ell(\Phi)=\Phi$ if and only if
$\Phi_\ast(\pi_1(Z,\underline{z}))\subset\Pi_\ast(\pi_1(Y,\underline{q}))$
in $\pi_1(X,\underline{p})$.
\end{lem}

\pf
The `only if' part is trivial as usual. We shall prove the
`if' part next.

We introduce some notations first. Let $(\{\pi_\alpha\},
\{\rho_{\beta\alpha}\})$ be a representative of $\Pi$ as in
Definition 2.4.2, which will be fixed throughout the proof. Let
$\underline{q}=(q,V_o,\hat{q})$, $\underline{p}=(p,U_o,\hat{p})$ and
$\underline{z}=(z,W_o,\hat{z})$ be the base-point structures.

We pick a representative of $\Phi$, denoted by $\sigma=(\{\phi_i\},
\{\eta_{ji}\}):\Gamma\{W_i\}\rightarrow\Gamma\{U_{i^\prime}\}$, where
we may assume that each $U_i\in\{U_{i^\prime}\}$ is an elementary
neighborhood. The strategy of the proof is to find a set $\{\delta_i\mid
\delta_i\in G_{U_i},\delta_o=1\in G_{U_o}\}$, such that $\delta(\sigma)
=(\{\delta(\phi_i)\},\{\delta(\eta_{ji})\})$, where $\delta(\phi_i)
=\delta_i\circ\phi_i$, $\delta(\eta_{ji})=\delta_j\circ\eta_{ji}\circ
\delta_i^{-1}$, can be lifted to $(Y,\underline{q})$
through $(\{\pi_\alpha\},\{\rho_{\beta\alpha}\})$ directly.

Such a set $\{\delta_i\}$ may be obtained as follows. For each $i$, we
pick a point $z_i\in W_i$ and a $\hat{z_i}\in\widehat{W_i}$ such that
$\pi_{W_i}(\hat{z_i})=z_i$, with $z_o=z$ and $\hat{z_o}=\hat{z}$ where
$z,\hat{z}$ are given in $\underline{z}=(z,W_o,\hat{z})$. This gives rise
to a set of base-point structures $\underline{z_i}=(z_i,W_i,\hat{z_i})$
of $Z$, and a corresponding set of base-point structures of $X$:
$\underline{p_i}=(p_i,U_i,\hat{p_i})$ where each $p_i=\phi(z_i)$ and
$\hat{p_i}=\phi_i(\hat{z_i})$, with $\underline{z_o}=\underline{z}$,
$\underline{p_o}=\underline{p}$. (Here $\phi:Z\rightarrow X$ is
the induced map of $\Phi$ between underlying spaces.) Since $Z$ is
connected and locally path-connected, it follows that for each $W_i$, there
exists a path $u_i\in P(Z,\underline{z},\underline{z_i})$, and therefore
a push-forward path $u_i^\prime=\sigma\circ u_i\in P(X,\underline{p},
\underline{p_i})$. By Lemma 2.4.4, there exist a connected component
$V_{\alpha(u_i^\prime)}$ of $\pi^{-1}(U_i)$ and a coset $\delta(u_i^\prime)$
such that after choosing a representative $\delta_i$ of $\delta(u_i^\prime)$,
the path $\delta_i\cdot u_i^\prime$ in $P(X,\underline{p},\delta_i\cdot
\underline{p_i})$ can be lifted to a path in $Y$ through $(\{\pi_\alpha\},
\{\rho_{\beta\alpha}\})$. The upshot is that the component
$V_{\alpha(u_i^\prime)}$ and the coset $\delta(u_i^\prime)$ are
independent of the choice on the path $u_i$, because of the assumption
that $\Phi_\ast(\pi_1(Z,\underline{z}))\subset\Pi_\ast(\pi_1(Y,
\underline{q}))$ in $\pi_1(X,\underline{p})$. Now we set
$V_i=V_{\alpha(u_i^\prime)}$ and choose a representative
$\delta_i\in G_{U_i}$ of the coset $\delta(u_i^\prime)$
for each $i$. We define $\delta(\sigma)$ using the set $\{\delta_i\}$
thus obtained. Note that $\delta(\sigma)$ has the following
property: for any path $u_i\in P(Z,\underline{z},
\underline{z_i})$, the push-forward path $\delta(\sigma)\circ u_i\in
P(X,\underline{p},\delta_i\cdot\underline{p_i})$ can be directly lifted
to a path in $(Y,\underline{q},\underline{q_i}(\delta_i))$ by
$(\{\pi_\alpha\},\{\rho_{\beta\alpha}\})$, where $\underline{q_i}(\delta_i)
=(q_i,V_i,\hat{q_i})$ is the base-point structure defined by
$\hat{q_i}=\pi_i^{-1}(\delta_i\cdot\hat{p_i})$ and
$q_i=\pi_{V_i}(\hat{q_i})$ (here $\pi_i:\widehat{V_i}\rightarrow
\widehat{U_i}$ is the homeomorphism given in $(\{\pi_\alpha\},
\{\rho_{\beta\alpha}\})$).

We prove next that $\delta(\eta_{ji}):T(W_i,W_j)\rightarrow
T(U_i,U_j)$ has its image contained in $\rho_{ji}(T(V_i,V_j))$,
in particular, each $\delta(\eta_i)=Ad(\delta_i)\circ\eta_i:G_{W_i}
\rightarrow G_{U_i}$ has its image contained in $\rho_i(G_{V_i})$,
so that we can define $\ell(\sigma)=(\{\ell(\phi_i)\},
\{\ell(\eta_{ji})\})$ by setting $\ell(\phi_i)=\pi_i^{-1}\circ
\delta(\phi_i)$ and $\ell(\eta_{ji})=\rho_{ji}^{-1}
\circ\delta(\eta_{ji})$. The desired map $\ell(\Phi):(Z,\underline{z})
\rightarrow (Y,\underline{q})$ is the equivalence class of $\ell(\sigma)$.

Given any $\xi\in T(W_i,W_j)$, we pick a point $\hat{z}^\prime
\in\mbox{Domain }(\phi_\xi)$. Since $Z$ is connected and locally
path-connected, there is a path $u_i\in P(Z,\underline{z},\underline{z_i})$,
a path $\gamma$ in $\widehat{W_i}$ joining $\hat{z_i}$ to $\hat{z}^\prime$,
and a path $\gamma^\prime$ in $\widehat{W_j}$ joining
$\phi_\xi(\hat{z}^\prime)$ to $\hat{z_j}$. Let $u_j=u_i\#\gamma\#_\xi
\gamma^\prime$ be the path in $P(Z,\underline{z},\underline{z_j})$ which is
defined by composing $u_i$ with $\gamma$ and then with
$\gamma^\prime$ through $\xi$. The push-forward path $\delta(\sigma)
\circ u_j\in P(X,\underline{p},\delta_j\cdot\underline{p_j})$ can be
lifted directly to $(Y,\underline{q},\underline{q}_j(\delta_j))$ through
$(\{\pi_\alpha\},\{\rho_{\beta\alpha}\})$, which implies that
$\delta(\eta_{ji})(\xi)$ lies in the image of $\rho_{ji}$.

Finally, we address the uniqueness of $\ell(\Phi)$. Let
$\Psi_1,\Psi_2:(Z,\underline{z})\rightarrow (Y,\underline{q})$ be any
maps such that $\Pi\circ\Psi_1=\Pi\circ\Psi_2=\Phi$. Since in
particular the induced maps between the underlying spaces
coincide, we may represent $\Psi_1$, $\Psi_2$ by $\tau_1
=(\{\psi_i^{(1)}\},\{\theta_{ji}^{(1)}\}):\Gamma\{W_i\}\rightarrow
\Gamma\{V_{i(1)}\}$ and $\tau_2=(\{\psi_i^{(2)}\},\{\theta_{ji}^{(2)}\}):
\Gamma\{W_i\}\rightarrow\Gamma\{V_{i(2)}\}$ respectively, where
$V_{i(1)},V_{i(2)}$ are connected components of $\pi^{-1}(U_i)$ for
some elementary neighborhood $U_i$. The compositions of $\tau_1$,
$\tau_2$ with $(\{\pi_\alpha\},\{\rho_{\beta\alpha}\})$ are equivalent,
because they represent the same map $\Phi$. Hence by Lemma 3.1.2 in
\cite{C1}, they are conjugate, and furthermore, because they also
preserve the base-point structures and $Z$ is connected, they are
actually equal. If $V_{i(1)}=V_{i(2)}$ for all $i$, then it is easily
seen that $\tau_1=\tau_2$ so that $\Psi_1=\Psi_2$. Suppose
$V_{i(1)}\neq V_{i(2)}$ for some index $i$. We pick a $\underline{z_i}
=(z_i,W_i,\hat{z}_i)$ and a path $u\in P(Z,\underline{z},\underline{z_i})$.
Then it is easily seen that $\Psi_1\circ u$ and $\Psi_2\circ u$ are
two different liftings of $\Phi\circ u$, contradicting the uniqueness of
path-lifting in Lemma 2.4.3. Hence $V_{i(1)}=V_{i(2)}$ for all $i$.

\hfill $\Box$

Next we discuss the existence of universal covering and deck
transformations.

\begin{prop}
Let $X$ be a connected, locally path-connected and
semi-locally 1-connected orbispace. Fix a base-point structure
$\underline{p}=(p,U_o,\hat{p})$ of $X$. Then for any subgroup $H$ of
$\pi_1(X,\underline{p})$, there is a connected covering space
$\Pi:(Y,\underline{q})\rightarrow (X,\underline{p})$, which is
unique up to isomorphisms, such that $\Phi_\ast(\pi_1(Y,\underline{q}))
=H$.
\end{prop}

\pf
We first construct the underlying space of $Y$. To this end, we
consider the path space $P(X,\underline{p})=[(I,0);(X,\underline{p})]$,
which is given a natural topology by the general method described in
\S 3.2 of \cite{C1}. We introduce an equivalence relation
$\sim_H$ in $P(X,\underline{p})$ as follows. Given any $u_1,u_2\in
P(X,\underline{p})$, we define $u_1\sim_H u_2$, if there are
$\sigma_1$, $\sigma_2$ representing $u_1, u_2$ respectively, and
there is a $\xi\in\bigsqcup_{i,j}T(U_i,U_j)$, where $\{U_i\}$ is
the atlas of local charts on $X$, such that $\sigma_1$ may compose
with $\nu(\sigma_2)$ via $\xi$ to form a homomorphism
$\sigma_1\#_\xi\nu(\sigma_2)$, which defines a loop in $(X,\underline{p})$
whose homotopy class lies in $H$. We define the underlying space of $Y$ to
be $P(X,\underline{p})/\sim_H$, which is obviously path-connected, hence
connected. As a notational convention, for any $u\in P(X,\underline{p})$,
we denote its equivalence class under $\sim_H$ by $u_H\in Y$. There is
a natural surjective continuous map $\pi:Y\rightarrow X$ defined by
$u_H\mapsto u(1)$, where $u(1)$ is the terminal point of the path $u$
in $X$. (Here the surjectivity of $\pi$ relies on the fact that $X$
is path-connected.) The space $Y$ has a natural base point
$q=\tilde{p}_H$, where $\tilde{p}\in P(X,\underline{p})$ is
the constant path, defined by the constant map into $\hat{p}\in
\widehat{U_o}$. Clearly $\pi(q)=p$.

For any point $y=u_H\in Y$, where $u$ is represented by $\sigma
=(\{\gamma_k\},\{\xi_{lk}\}):\Gamma\{I_k\}\rightarrow
\Gamma\{U_{k^\prime}\}$, $k=0,1,\cdots,n$, we set $U_\sigma=U_n$
and $\hat{\sigma}=\gamma_n(1)\in \widehat{U_\sigma}$. We may assume that
the semi-locally 1-connectedness holds for $U_\sigma$ without loss of
generality. We define a map $\pi_\sigma:\widehat{U_\sigma}\rightarrow Y$
as follows. For each $z\in\widehat{U_\sigma}$, we connect $\hat{\sigma}$
to $z$ by a path $\gamma_{z}$ in $\widehat{U_\sigma}$ (the existence of
$\gamma_z$ is ensured by the locally path-connectedness of $X$), and define
$\pi_\sigma(z)=[\sigma\#\gamma_z]_H$. The assumption that $X$ is
semi-locally 1-connected ensures that the map $\pi_\sigma$ is well-defined,
and the assumption that $X$ is locally path-connected implies
that $\pi_\sigma$ is continuous. We set $V_\sigma=
\pi_\sigma(\widehat{U_\sigma})
\subset Y$, which obviously satisfies $\pi(V_\sigma)=U_\sigma$.

We will show:
(1) $V_\sigma$ is a connected open neighborhood of $y$ in $Y$.
(2) There is a subgroup $G_{V_\sigma}$ of $G_{U_\sigma}$ such that
$\widehat{U_\sigma}/G_{V_\sigma}$ is homeomorphic to $V_\sigma$ under
$\pi_\sigma$. (3) There is a set $\T=\{T(V_\sigma,V_\tau)\}$, where
each $T(V_\sigma,V_\tau)$ is a subset of $T(U_\sigma,U_\tau)$.
Together with the atlas of local charts $\{(\widehat{V_\sigma},G_{V_\sigma},
\pi_\sigma)\}$, where $\widehat{V_\sigma}=\widehat{U_\sigma}$,
it defines an orbispace structure on $Y$ (cf. Proposition 2.1.1 in
\cite{C1}). (4) The maps $\pi^\sigma:\widehat{V_\sigma}=\widehat{U_\sigma}$
and $\rho_{\tau\sigma}:T(V_\sigma,V_\tau)\subset T(U_\sigma,U_\tau)$
define a covering map $\Pi:Y\rightarrow X$. (5) We have
$\Pi_\ast(\pi_1(Y,\underline{q}))=H$ where in $\underline{q}
=(q,V_o,\hat{q})$, $V_o=\pi_{\sigma_q}(\widehat{U_o})$, and $\hat{q}
=\hat{p}\in\widehat{U_o}=\widehat{V_o}$. Here $\sigma_q$ is the
canonical representative of $q$, the constant map into $\hat{p}$.

For (1), we first prove that $V_\sigma$ is open. It suffices to show
that the inverse image of $V_\sigma$ under the projection
$P(X,\underline{p})\rightarrow Y$ is open in $P(X,\underline{p})$.
To this end, we need to show that given any $u_0\in P(X,\underline{p})$
with $(u_0)_H=\pi_\sigma(z_0)$ for some $z_0\in\widehat{U_\sigma}$, and
for any $u\in P(X,\underline{p})$ sufficiently close to $u_0$, there is
a $z\in\widehat{U_\sigma}$ such that $\pi_\sigma(z)=u_H$. Let $\gamma_{z_0}$
be a path in $\widehat{U_\sigma}$ connecting $\hat{\sigma}$ to $z_0$,
such that $\pi_\sigma(z_0)=[\sigma\#\gamma_{z_0}]_H$. We take a
$\tau_0=\{f_{i,0}\}\in\O_{\{\xi_{ji}\}}$, $0\leq i\leq m$, such
that $[\tau_0]=u_0$, where without loss of generality, we assume that
$U_m=U_\sigma$, $f_{m,0}(1)=z_0$, and the homotopy class of the loop
$[\tau_0\#\nu(\sigma\#\gamma_{z_0})]$ lies in $H$. Now we pick a
$t_i\in I_i\cap I_{i+1}$ for each $i=0,\cdots,m-1$, and for each
$\tau=\{f_{i}\}$ in the open neighborhood $\O_{\{\xi_{ji}\}}$ of $\tau_0$,
we take a path $\gamma_i$ in $\mbox{Domain }(\phi_{\xi_{(i+1)i}})$ running
from $f_i(t_i)$ to $f_{i,0}(t_i)$ for each $i=0,\cdots,m-1$, and take a
path $\gamma$ in $\widehat{U_\sigma}$ running from $f_m(1)$ to $f_{m,0}(1)=z_0$.
We define $\tau^\prime=\{f_{i}^\prime\}\in\O_{\{\xi_{ji}\}}$ where
for each $0\leq i\leq m-1$, $f_i^\prime$ is obtained from pre-composing
$f_i$ by $\phi_{\xi_{i(i-1)}}\circ \nu(\gamma_{i-1})$ and post-composing
$f_i$ by $\gamma_i$, and $f_m^\prime$ is obtained from
pre-composing $f_m$ by $\phi_{\xi_{m(m-1)}}\circ \nu(\gamma_{m-1})$ and
post-composing $f_m$ by $\gamma\# \nu(\gamma)$. Then $\tau^\prime$ is
homotopic to $\tau$ so that $[\tau^\prime]_H=[\tau]_H$. On the other hand,
by the semi-locally 1-connectedness of $X$, the loop
$[\tau^\prime\#\nu(\tau_0\#\nu(\gamma))]$ is null-homotopic, so that
$[\tau^\prime]_H=[\tau_0\#\nu(\gamma)]_H$. Now observe that
$\tau_0\#\nu(\gamma)\#\nu(\sigma\# (\gamma_{z_0}\#\nu(\gamma)))$ is
homotopic to $\tau_0\#\nu(\sigma\#\gamma_{z_0})$, which defines a loop
whose homotopy class lies in $H$. If we let $z=f_m(1)$, then we have
$\pi_\sigma(z)=[\sigma\# (\gamma_{z_0}\#\nu(\gamma))]_H=
[\tau_0\#\nu(\gamma)]_H=[\tau^\prime]_H=[\tau]_H$. Hence $V_\sigma$
is open in $Y$. Finally, $V_\sigma$ is connected because $\widehat{U_\sigma}$
is and $V_\sigma=\pi_\sigma(\widehat{U_\sigma})$.
This concludes the proof of (1).

For (2), we obtain the subgroup $G_{V_\sigma}$ as follows. We denote
by $\underline{\sigma}$ the base-point structure $(\pi(y),U_\sigma,
\hat{\sigma})$, and consider the isomorphism $[\sigma]_\ast:
\pi_1(X,\underline{\sigma})\rightarrow\pi_1(X,\underline{p})$, where
$[\sigma]$ is the path in $P(X,\underline{p},\underline{\sigma})$ defined
by $\sigma$, cf. Proposition 1.3 (1). Set $H_\sigma=[\sigma]_\ast^{-1}(H)$.
Now for each $g\in G_{U_\sigma}$, we take a path $\gamma_g$ in
$\widehat{U}_\sigma$
connecting $\hat{\sigma}$ and $g\cdot\hat{\sigma}$.  Then the pair
$(\gamma_g,g)$ determines an element $[(\gamma_g,g)]$ in $\pi_1(U_\sigma,
\underline{\sigma})$, cf. Lemma 2.2.2, whose image $[g]$ under
$\pi_1(U_\sigma,\underline{\sigma})\rightarrow\pi_1(X,\underline{\sigma})$
is independent of the choice on the path $\gamma_g$ by the semi-locally
1-connectedness of $X$. We simply put $G_{V_\sigma}=\{g\in G_{U_\sigma}\mid [g]
\in H_\sigma\}$.  The fact that $(\gamma_g,g)\# (\gamma_h,h)=(\gamma_{gh},gh)$
for some path $\gamma_{gh}$ connecting $\hat{\sigma}$ to $gh\cdot \hat{\sigma}$
in $\widehat{U}_\sigma$ shows that $G_{V_\sigma}$ is a subgroup of
$G_{U_\sigma}$. It remains to show that $\pi_\sigma$ induces a
homeomorphism between $\widehat{U_\sigma}/G_{V_\sigma}$ and
$V_\sigma$. First of all, $\pi_\sigma: \widehat{U_\sigma}\rightarrow
V_\sigma$ is $G_{V_\sigma}$-invariant. This can be seen as follows. Let
$g$ be any element in $G_{V_\sigma}$. Given any $z\in\widehat{U_\sigma}$,
let $\gamma_z$ be a path in $\widehat{U_\sigma}$ running from
$\hat{\sigma}$ to $z$ such that $\pi_\sigma(z)=[\sigma\#\gamma_z]_H$.
Then $\pi_\sigma(g\cdot z)=[\sigma\# (\gamma_g\# g\circ\gamma_z)]_H$ because
$\gamma_g\# g\circ\gamma_z$ is a path in $\widehat{U_\sigma}$ connecting
$\hat{\sigma}$ and $g\cdot z$. Now it is easily seen that
$\sigma\# (\gamma_g\# g\circ\gamma_z)$ may compose with $\nu(\sigma\#\gamma_z)$
through $g^{-1}\in G_{U_\sigma}$, and the resulting homomorphism
represents a loop whose homotopy class is $[\sigma]_\ast([g])\in
H$. Hence $\pi_\sigma(g\cdot z)=\pi_\sigma(z)$.  Secondly, suppose
$\pi_\sigma(z)=\pi_\sigma(z^\prime)$. Take paths $\gamma_z$,
$\gamma_{z^\prime}$ connecting $\hat{\sigma}$ to $z$ and $z^\prime$
respectively. Then the assumption $\pi_\sigma(z)=\pi_\sigma(z^\prime)$
implies that there is some $g\in G_{U_\sigma}$ with $z^\prime=g\cdot z$
such that $\sigma\#\gamma_z$ may compose with $\nu(\sigma\#\gamma_{z^\prime})$
through $g$ to form a homomorphism which defines a loop in $X$, whose homotopy
class is an element $h\in H$. If we join $\hat{\sigma}$ and
$g\cdot\hat{\sigma}$ by the path $\gamma_g=\gamma_{z^\prime}\#
\nu(g\circ\gamma_z)$, we see that the class $[(\gamma_g,g)]$ in
$\pi_1(U_\sigma,\underline{\sigma})$ has its image $[g]=
[\sigma]_\ast^{-1}(h^{-1})\in H_\sigma$ under $\pi_1(U_\sigma,
\underline{\sigma})
\rightarrow \pi_1(X,\underline{\sigma})$.  Hence $g$ lies in $G_{V_\sigma}$.
From here it is easy to see that $\pi_\sigma$ induces a homeomorphism
between $\widehat{U_\sigma}/G_{V_\sigma}$ and $V_\sigma$.

For (3), suppose
$(\widehat{V_{\sigma_1}},G_{V_{\sigma_1}},\pi_{V_{\sigma_1}})$,
$(\widehat{V_{\sigma_2}},G_{V_{\sigma_2}},\pi_{V_{\sigma_2}})$
are two local charts constructed from $\sigma_1$, $\sigma_2$ respectively,
such that $V_{\sigma_1}\cap V_{\sigma_2}\neq\emptyset$. The assumption
$V_{\sigma_1}\cap V_{\sigma_2}\neq \emptyset$ implies that $U_{\sigma_1}
\cap U_{\sigma_2}\neq \emptyset$ since $\pi(V_\sigma)=U_\sigma$. We shall
define $T(V_{\sigma_1},V_{\sigma_2})$ as a subset of $T(U_{\sigma_1},
U_{\sigma_2})$ as follows. Given any $\xi\in T(U_{\sigma_1},U_{\sigma_2})$,
if there exists a $z\in\mbox{Domain }(\phi_\xi)$ such that the composition
of $\sigma_1\#\gamma_z$ with $\nu(\sigma_2\#\gamma_{\phi_\xi(z)})$ through
$\xi$ defines a loop whose homotopy class lies in $H$, then we put $\xi$ in
$T(V_{\sigma_1},V_{\sigma_2})$. Here $\gamma_z$, $\gamma_{\phi_\xi(z)}$
are paths in $\widehat{U_{\sigma_1}}$, $\widehat{U_{\sigma_2}}$ connecting
$\hat{\sigma}_1$, $\hat{\sigma}_2$ to $z$, $\phi_\xi(z)$ respectively.
Note that $T(V_\sigma,V_\sigma)=G_{V_\sigma}$ according to this
definition. On the other hand, since $\mbox{Domain }(\phi_\xi)$ is
path-connected, $\xi\in T(V_{\sigma_1},V_{\sigma_2})$ implies
that $\pi_{\sigma_1}=\pi_{\sigma_2}\circ\phi_\xi$ on $\mbox{Domain }(\phi_\xi)$
and that $\pi_{\sigma_1}(\mbox{Domain }(\phi_\xi))$ is a connected component of
$V_{\sigma_1}\cap V_{\sigma_2}$. This allows us to assign $\phi_\xi$
to $\xi$ even if $\xi$ is regarded as an element of
$T(V_{\sigma_1},V_{\sigma_2})$. Finally, we
observe that $\bigsqcup_{\sigma_1,\sigma_2}
T(V_{\sigma_1},V_{\sigma_2})$ is closed under
taking inverse and composition as a subset of $\bigsqcup_{i,j}
T(U_i,U_j)$. Hence by Proposition 2.1.1 in \cite{C1}, it together with the
atlas of local charts $\{(\widehat{V_\sigma},G_{V_\sigma},\pi_\sigma)\}$
defines an orbispace structure on $Y$. This concludes the proof of (3).

For (4), in order to show that $(\{\pi^\sigma\},\{\rho_{\tau\sigma}\})$
defines a covering map $\Pi:Y\rightarrow X$, it suffices to verify (b) of
Definition 2.4.2 for it. More concretely, we need to show that for any
$U\in\{U_\sigma\}$, a connected component of $\pi^{-1}(U)$ is of the form
$V_\sigma$. Let $V$ be any connected component of $\pi^{-1}(U)$. Then
for any $y\in V$, since $\pi(y)\in U$, there is a $\sigma_y$
representing $y$, with a canonically constructed neighborhood
$V_{\sigma_y}$ satisfying $\widehat{V_{\sigma_y}}=\widehat{U}$. Since
$V_{\sigma_y}$ is connected and $\pi(V_{\sigma_y})=U$, we have
$V_{\sigma_y}\subset V$. Hence $V=\cup_{y\in V}V_{\sigma_y}$. On the
other hand, it is easy to see that if $V_\sigma$ and $V_\tau$ have
non-empty intersection and $\pi(V_\sigma)=\pi(V_\tau)$, then
$V_\sigma=V_\tau$. Hence $V=V_{\sigma_y}$ for any $y\in V$, and
(b) is verified.

For (5), we first show that $H\subset\Pi_\ast(\pi_1(Y,\underline{q}))$.
Let $u$ be any loop in $X$ whose homotopy class lies in $H$. We
represent it by a homomorphism $\sigma\in\O_{\{\xi_{ji}\}}$,
$\sigma:\Gamma\{I_i\}\rightarrow\Gamma\{U_{i^\prime}\}$, $0\leq i\leq n$.
For each $i=1,\cdots,n$, we pick a $t_i\in I_{i-1}\cap I_i$ and let
$\tau_i$ be the restriction of $\sigma$ on $[t,t_i]$. We let $V_{\sigma_i}$
be the component of $\pi^{-1}(U_i)$ that contains $[\tau_i]_H\in Y$,
denote by $z_i\in\widehat{V_{\sigma_i}}$ a point satisfying $\pi_{\sigma_i}
(z_i)=[\tau_i]_H$, and denote by $\gamma_{z_i}$ a path in $\widehat
{V_{\sigma_i}}$ that connects $\hat{\sigma}_i$ to $z_i$. Then there
exist $g_i\in G_{U_i}$, $i=1,\cdots,n$, such that the composition
of $\tau_i$ with $\nu(\sigma_i\#\gamma_{z_i})$ through $g_i\circ\xi_{i(i-1)}$
defines a loop whose homotopy class lies in $H$. This implies that
$g_1\circ\xi_{10}\in T(V_o,V_{\sigma_1})$, $g_i\circ\xi_{i(i-1)}\circ
g_{i-1}^{-1}\in T(V_{\sigma_{i-1}},V_{\sigma_i})$ for $i\geq 2$. Note that
because the homotopy class of the loop $u$ lies in $H$, $V_{\sigma_n}=V_o$
and $g_n$ may be taken to be $1\in G_{U_o}$. Thus according to Lemma 2.4.4,
the loop $u$ can be lifted to a loop in $(Y,\underline{q})$.
Hence $H\subset\Pi_\ast(\pi_1(Y,\underline{q}))$. It remains to show
that $\Pi_\ast(\pi_1(Y,\underline{q}))\subset H$. Given any $\tau\in
\O_{\{\eta_{ji}\}}$, $i=0,1,\cdots,n$, where $\eta_{ji}\in T(V_{\sigma_i},
V_{\sigma_j})$ such that the equivalence class $[\tau]\in
[(S^1,\ast);(Y,\underline{q})]$. Then there are, for $i=0,\cdots,n-1$,
points $z_i\in\mbox{Domain }(\phi_{\eta_{(i+1)i}})$, $z_i^\prime=
\phi_{\eta_{(i+1)i}}(z_i)$, and paths $\gamma_{z_i}$ in
$\widehat{V_{\sigma_i}}$, $\gamma_{z_i^\prime}$ in
$\widehat{V_{\sigma_{i+1}}}$,
which connect $\hat{\sigma}_i$, $\hat{\sigma}_{i+1}$ to $z_i$, $z_i^\prime$
respectively, such that the composition of $\sigma_i\#\gamma_{z_i}$ with
$\nu(\sigma_{i+1}\#\gamma_{z_i^\prime})$ through $\eta_{(i+1)i}$
defines a loop whose homotopy class is an element $h_i\in H$. By the
semi-locally 1-connectedness of $X$ and the path-connectedness of each
$\mbox{Domain }(\phi_{\eta_{(i+1)i}})$, we see that $\Pi\circ [\tau]$
is a loop whose homotopy class equals $h_0\cdots h_{n-1}\in H$. Hence
$\Pi_\ast(\pi_1(Y,\underline{q}))\subset H$.

\hfill $\Box$

We shall derive a short exact sequence, which relates the $\pi_1$
of an orbispace with the group of deck transformations of its
universal covering, whose existence was established in the
preceding proposition. By definition, a deck transformation of a
covering space $\Pi:Y\rightarrow X$ is an isomorphism of orbispaces
$\Phi:Y\rightarrow Y$ such that $\Pi\circ\Phi=\Pi$. The group of
deck transformations of $\Pi:Y\rightarrow X$ is denoted by
$\mbox{Deck }(\Pi)$.

Let $\Pi:(Y,\underline{q})\rightarrow (X,\underline{p})$ be a covering
map with a representative $(\{\pi_\alpha\},\{\rho_{\beta\alpha}\})$ fixed
throughout. For each $\alpha$, we set $K_{U_\alpha}=\{g\in G_{U_\alpha}\mid
g\cdot x=x,\forall x\in\widehat{U_\alpha}\}$. We consider
$$
K_\Pi=\{\{g_\alpha\}\mid g_\alpha\in K_{U_\alpha} \mbox{ and }
\rho_{\beta\alpha}(\xi)=g_\beta\circ\rho_{\beta\alpha}(\xi)
\circ g_\alpha^{-1} \;\forall \xi\in T(V_\alpha,V_\beta)\}, \leqno (2.4.6)
$$
which is a group under the multiplication $\{g_\alpha\}\{h_\alpha\}=
\{g_\alpha h_\alpha\}$. The subgroup
$$
C_\Pi=\{\{g_\alpha\}\in K_\Pi\mid g_\alpha\in\rho_\alpha(G_{V_\alpha})\}
\leqno (2.4.7)
$$
is contained in the center of $K_\Pi$, hence is a normal subgroup
of $K_\Pi$.

Let $\tilde{p}$ and $\tilde{q}$ be the constant map from $I$ into
$\hat{p}$ and $\hat{q}$ respectively. For any $g\in K_{U_o}$,
$h\in K_{V_o}$, we set $[g]=[(\tilde{p},g)]\in\pi_1(X,\underline{p})$,
$[h]=[(\tilde{q},h)]\in\pi_1(Y,\underline{q})$. Clearly $\Pi_\ast([h])
=[\rho_o(h)]$. For any $\{g_\alpha\}\in K_\Pi$, we observe that $[g_o]$
actually lies in the center of $\Pi_\ast(\pi_1(Y,\underline{q}))$ in
$\pi_1(X,\underline{p})$, where $g_o$ is the component of $\{g_\alpha\}$
corresponding to $U_o$. To see this, let $u\in \pi_1(Y,\underline{q})$
be an element which is represented by $\{f_i\}\in\O_{\{\xi_{ji}\}}$,
$0\leq i\leq n$. Then $\Pi_\ast(u)$ is represented by $\{\pi_i\circ\gamma_i\}
\in\O_{\{\rho_{ji}(\xi_{ji})\}}$, and $[g_o]^{-1}\#\Pi_\ast(u)\#[g_o]$ is
represented by $\{\pi_i\circ\gamma_i\}\in\O_{\{\eta_{ji}\}}$, where
$\eta_{10}=\rho_{10}(\xi_{10})\circ g_o$, $\eta_{ji}=\rho_{ji}(\xi_{ji})$
for $2\leq j\leq n-1$, and $\eta_{n(n-1)}=g_o^{-1}\circ\rho_{n(n-1)}
(\xi_{n(n-1)})$. Observe that if we let $g_i$ be the component of
$\{g_\alpha\}$ corresponding to $U_i$, then $\rho_{ji}(\xi_{ji})=
g_j\circ\rho_{ji}(\xi_{ji})\circ g_i^{-1}$ are satisfied, which
imply that $g_1^{-1}\circ\eta_{10}=\rho_{10}(\xi_{10})$, $g_j^{-1}\circ
\eta_{ji}\circ g_i=\rho_{ji}(\xi_{ji})$ for $2\leq j\leq n-1$, and
$\eta_{n(n-1)}\circ g_n=\rho_{n(n-1)}(\xi_{n(n-1)})$. In other words,
$(\{\pi_i\circ\gamma_i\},\{\eta_{ji}\})$ is conjugate to
$(\{\pi_i\circ\gamma_i\},\{\rho_{ji}(\xi_{ji})\})$. Hence
$[g_o]^{-1}\#\Pi_\ast(u)\#[g_o]=\Pi_\ast(u)$.

Thus there is a homomorphism
$$
\Xi_{p,q}:K_\Pi/C_\Pi\rightarrow N(\Pi_\ast(\pi_1(Y,\underline{q})))
/\Pi_\ast(\pi_1(Y,\underline{q})), \leqno (2.4.8)
$$
induced by $\{g_\alpha\}\mapsto [g_o]$, where $N(\Pi_\ast(\pi_1(Y,
\underline{q})))$ is the normalizer of $\Pi_\ast(\pi_1(Y,\underline{q}))$
in $\pi_1(X,\underline{p})$. Note that $\Xi_{p,q}$ is injective when $Y$ is
connected. This is because for any $\{g_\alpha\}\in K_\Pi$, $g_o\in
\rho_o(G_{V_o})$ implies $g_\alpha\in\rho_\alpha(G_{V_\alpha})$ by
$\rho_{\beta\alpha}(\xi)=g_\beta\circ\rho_{\beta\alpha}(\xi)
\circ g_\alpha^{-1} \;\forall \xi\in T(V_\alpha,V_\beta)$ in $(2.4.6)$,
and by the connectedness of $Y$.

\begin{prop}
Let $\Pi:(Y,\underline{q})\rightarrow (X,\underline{p})$ be a
connected, locally path-connected covering space of $X$. Set
$H=\Pi_\ast(\pi_1(Y,\underline{q}))$ and denote by $N(H)$ the
normalizer of $H$ in $\pi_1(X,\underline{p})$. Then there is a
homomorphism $\Theta_{p,q}:N(H)/H\rightarrow\mbox{Deck }(\Pi)$,
such that $\Xi_{p,q}$ and $\Theta_{p,q}$ fit into a short exact
sequence
$$
1\rightarrow K_{\Pi}/C_{\Pi}\stackrel{\Xi_{p,q}}{\longrightarrow}
N(H)/H \stackrel{\Theta_{p,q}}{\longrightarrow}\mbox{Deck }(\Pi)
\rightarrow 1. \leqno (2.4.9)
$$
Moreover, if $\Pi:(Y,\underline{q})\rightarrow (X,\underline{p})$,
$\Pi^\prime:(Y^\prime,\underline{q^\prime})\rightarrow
(X^\prime,\underline{p^\prime})$ are connected, locally
path-connected covering spaces, where $\Pi_\ast(\pi_1(Y,\underline{q}))$,
$\Pi_\ast^\prime(\pi_1(Y^\prime,\underline{q^\prime}))$ are normal
subgroups, and there are maps\footnote{here $\Phi,\Psi$ are
general maps, i.e., not in the restricted class specified by
{\bf Convention} in Introduction.}
$\Phi:(X,\underline{p})\rightarrow (X^\prime,
\underline{p^\prime})$, $\Psi:(Y,\underline{q})\rightarrow
(Y^\prime,\underline{q^\prime})$ such that $\Pi^\prime\circ
\Psi=\Phi\circ\Pi$, then we have
$$
\Psi\circ \Theta_{p,q}(z)=\Theta_{p^\prime,q^\prime}(z^\prime)\circ\Psi,\;
\forall z\in \pi_1(X,\underline{p})/\Pi_\ast(\pi_1(Y,\underline{q})),
\leqno (2.4.10)
$$
where in $(2.4.10)$, $z^\prime$ is the image of $z$ in
$\pi_1(X^\prime,\underline{p^\prime})/\Pi_\ast^\prime(\pi_1(Y^\prime,
\underline{q^\prime}))$ under $\Phi_\ast$, and $\Psi$ is regarded
as an element of $[Y;Y^\prime]$.
\end{prop}

\pf
Let $u$ be a loop in $(X,\underline{p})$ such that its homotopy
class $[u]\in N(H)$. By Lemma 2.4.4, we associate $u$ with a pair
$(V_{\alpha(u)},\delta(u))$, where $V_{\alpha(u)}$ is a connected
component of $\pi^{-1}(U_o)$ and $\delta(u)$ is a left coset in
$G_{U_o}/\rho_{\alpha(u)}(G_{V_{\alpha(u)}})$, such that for any
chosen representative $g$ of $\delta(u)$, there is a unique path
$\ell(u)_g\in P(Y,\underline{q},\underline{q}(g))$ satisfying
$\Pi\circ\ell(u)_g=g\cdot u$. We consider the covering maps
$\Pi^{(1)}=\Pi:(Y,\underline{q})\rightarrow (X,\underline{p})$ and
$\Pi^{(2)}:(Y,\underline{q}(g))\rightarrow (X,\underline{p})$.
The assumption that $[u]\in N(H)$ implies that
$\Pi^{(1)}_\ast(\pi_1(Y,\underline{q}))=\Pi^{(2)}_\ast
(\pi_1(Y,\underline{q}(g)))$ in $\pi_1(X,\underline{p})$. Hence by
Lemma 2.4.6, there is a unique map $\Phi_u:(Y,\underline{q})\rightarrow
(Y,\underline{q}(g))$ such that $\Pi^{(2)}\circ\Phi_u=\Pi^{(1)}$. If
we pick a different representative $g^\prime$ of $\delta(u)$,
where $g^\prime=\rho_{\alpha(u)}(h)g$ for some $h\in G_{V_{\alpha(u)}}$,
then $\underline{q}(g)$ is changed to $\underline{q}(g^\prime)=
h\cdot\underline{q}(g)$. This implies that the corresponding map
$\Phi_u:Y\rightarrow Y$ is independent of the choice of the
representative $g$. As for the dependence on $(\{\pi_\alpha\},
\{\rho_{\beta\alpha}\})$, we go back to the proof of Lemma 2.4.4 and
observe that in the construction of $\ell(\sigma)$, if we change
$(\{\pi_\alpha\},\{\rho_{\beta\alpha}\})$ by conjugation, then the
set $\{g_i\}$ in $(2.4.1)$ will change accordingly such that overall,
the orbit of $\underline{q^\prime}(g)$ under the action of
$G_{V_{\alpha(u)}}$ remains the same. On the other hand, by
Remark 2.4.3 (1), two different
choices of $(\{\pi_\alpha\},\{\rho_{\beta\alpha}\})$ differ essentially
by a conjugation after we replace $\sigma$ by an appropriate induced one.
Hence $\Phi_u:Y\rightarrow Y$ is also independent of the choice on
$(\{\pi_\alpha\},\{\rho_{\beta\alpha}\})$. It is obvious that $\Phi_u$
depends only on the class of $[u]$ in $N(H)/H$, and is a deck transformation
of $\Pi:Y\rightarrow X$. The map $\Theta_{p,q}:N(H)/H\rightarrow
\mbox{Deck }(\Pi)$ is defined by $[u]\mapsto \Phi_u$.

In order to verify $(2.4.10)$, we pick representatives $\sigma$, $\tau$
for $\Phi$ and $\Psi$ respectively, and fix $(\{\pi_\alpha\},
\{\rho_{\beta\alpha}\})$, $(\{\pi^\prime_{\alpha^\prime}\},
\{\rho^\prime_{\beta^\prime\alpha^\prime}\})$ for $\Pi$, $\Pi^\prime$,
such that
$$
\sigma\circ (\{\pi_\alpha\},\{\rho_{\beta\alpha}\})
=(\{\pi^\prime_{\alpha^\prime}\},\{\rho^\prime_{\beta^\prime
\alpha^\prime}\})\circ\tau. \leqno (2.4.11)
$$
Given any $z\in\pi_1(X,\underline{p})/\Pi_\ast(\pi_1(Y,\underline{q}))$,
we represent $z$ by a loop $u$ in $(X,\underline{p})$. We pick a
representative $g$ of $\delta(u)$ (w.r.t $(\{\pi_{\alpha}\},
\{\rho_{\beta\alpha}\})$), and then use $\sigma$ to
determine a representative $g^\prime$ of $\delta(\Phi\circ u)$ (w.r.t
$(\{\pi^\prime_{\alpha^\prime}\},\{\rho^\prime_{\beta^\prime
\alpha^\prime}\})$), such that $\tau$ preserves $\underline{q}(g)$ and
$\underline{q^\prime}(g^\prime)$ (because of $(2.4.11)$). Now we consider
the problem of factoring the map $\Phi\circ\Pi:(Y,\underline{q})
\rightarrow (X^\prime,\underline{p^\prime})$ through the covering map
$(Y^\prime,\underline{q^\prime}(g^\prime))\rightarrow (X^\prime,
\underline{p^\prime})$. There are apparently two solutions:
$\Psi\circ \Theta_{p,q}(z)$ and $\Theta_{p^\prime,q^\prime}(z^\prime)
\circ\Psi$ (in their appropriate based versions). By the uniqueness in
Lemma 2.4.6, these two solutions must be the same. Hence $(2.4.10)$.

Next we verify that $\Theta_{p,q}$ is a homomorphism. Let $u_1,u_2$
be two loops in $(X,\underline{p})$ such that both $[u_1],[u_2]\in N(H)$.
We pick representatives $g_1$ of $\delta(u_1)$, $g_2$ of $\delta(u_2)$,
and $g$ of $\delta(u_1\# u_2)$. Then by the uniqueness of path-lifting
in Lemma 2.4.4, the following holds for the based version
$\Theta_{p,q}([u_1]):(Y,\underline{q},\underline{q}(g_2))\rightarrow
(Y,\underline{q}(g_1),\underline{q}(g))$,
$$
\ell(u_1)_{g_1}\# (\Theta_{p,q}([u_1])\circ\ell(u_2)_{g_2})=
\ell(u_1\# u_2)_g. \leqno (2.4.12)
$$
In particular, there is a based version
$$
\Theta_{p,q}([u_1]):(Y,\underline{q}(g_2))\rightarrow
(Y,\underline{q}(g)). \leqno (2.4.13)
$$
Now by the uniqueness in Lemma 2.4.6, the composition of
$\Theta_{p,q}([u_2]):(Y,\underline{q})\rightarrow (Y,\underline{q}(g_2))$
with $(2.4.13)$ must be equal to $\Theta_{p,q}([u_1\# u_2]):(Y,\underline{q})
\rightarrow (Y,\underline{q}(g))$. The corresponding maps satisfy
$$
\Theta_{p,q}([u_1\# u_2])=\Theta_{p,q}([u_1])\circ\Theta_{p,q}([u_2]),
\leqno (2.4.14)
$$
which shows that $\Theta_{p,q}$ is a homomorphism.\footnote{note that
$\Theta_{p,q}$ induces a left-action of $N(H)/H$ on $Y$.}

It remains to verify $(2.4.9)$. As for the surjectivity of $\Theta_{p,q}$,
given any $\Phi\in\mbox{Deck }(\Pi)$, we take a path $u\in P(Y,\underline{q},
\underline{q^\prime})$ where $\underline{q^\prime}$ is the image of
$\underline{q}$ under a choice of based versions of $\Phi$. Then
$\Pi\circ u$ determines a loop $v$ in $(X,\underline{p})$, and
$$
\Pi_\ast(\pi_1(Y,\underline{q}))=\Pi_\ast(\pi_1(Y,\underline{q^\prime}))
=\Pi_\ast(\nu(u)_\ast(\pi_1(Y,\underline{q})))=[v]^{-1}\cdot\Pi_\ast
(\pi_1(Y,\underline{q}))\cdot [v], \leqno (2.4.15)
$$
which implies that $[v]\in N(H)$. The uniqueness in Lemma 2.4.6 then
asserts that $\Theta_{p,q}([v])=\Phi$. Hence $\Theta_{p,q}$ is surjective.

Finally, we determine the kernel of $\Theta_{p,q}$. Suppose
$\Theta_{p,q}([u])=1$ for some loop $u$ in $(X,\underline{p})$. Then
from the construction in Lemma 2.4.6 and by the definition of
$\Theta_{p,q}$ above, we see that $\Theta_{p,q}([u])$ is defined by
$\ell((\{\pi_\alpha\},\{\rho_{\beta\alpha}\}))$ as constructed in
Lemma 2.4.6, where there exists a set $\{g_\alpha\mid g_\alpha\in
G_{U_\alpha}\}$, such that $\ell((\{\pi_\alpha\},\{\rho_{\beta\alpha}\}))
=(\{f_\alpha\},\{\eta_{\beta\alpha}\}):\Gamma\{V_\alpha\}\rightarrow
\Gamma\{V_\alpha\}$ with $f_\alpha=\pi_\alpha^{-1}\circ g_\alpha\circ
\pi_\alpha$, $\eta_{\beta\alpha}(\xi)=\rho_{\beta\alpha}^{-1}(g_\beta\circ
\rho_{\beta\alpha}(\xi)\circ g_\alpha^{-1})$. Moreover, a different choice
of $\{g_\alpha\}$ has the form $\{\rho_\alpha(h_\alpha)g_\alpha\}$ for
some $h_\alpha\in G_{V_\alpha}$. Now the assumption $\Theta_{p,q}([u])=1$
implies that $\ell((\{\pi_\alpha\},\{\rho_{\beta\alpha}\}))$ is conjugate
to the identity, and because of this, we can actually choose a set
$\{g_\alpha\}$ such that $\ell((\{\pi_\alpha\},\{\rho_{\beta\alpha}\}))$
is the identity. It then follows easily that $\{g_\alpha\}\in K_\Pi$.
Furthermore, one is ready to check that $u\# (\tilde{p},g_o^{-1})$,
where $g_o$ is the component of $\{g_\alpha\}$ that corresponds to $U_o$,
can be lifted to a loop in $(Y,\underline{q})$, hence $[u]=[g_o] \bmod
\Pi_\ast(\pi_1(Y,\underline{q}))$. This exactly means that $\ker\;
\Theta_{p,q}\subset\mbox{Im}\;\Xi_{p,q}$. It remains to show that
$\mbox{Im}\;\Xi_{p,q}\subset\ker\;\Theta_{p,q}$, so that $\ker\;
\Theta_{p,q}=\mbox{Im}\;\Xi_{p,q}$ and $(2.4.9)$ is verified. To this end,
let $\{g_\alpha\}\in K_\Pi$ be any element, we consider $\Theta_{p,q}([u])$
where $[u]=[g_o]$. It is clear that the component $V_{\alpha(u)}$ associated
to $u$ is $V_o$, and $g_o$ can be chosen to serve as a representative $g$
of $\delta(u)$. With $g=g_o$, we see that $\underline{q}(g)$ in the
definition of $\Theta_{p,q}([u])$ is actually $\underline{q}$ since
$g_o\in K_{U_o}$. By the uniqueness in Lemma 2.4.6, we conclude
that $\Theta_{p,q}([u])$ is the identity map. Hence $\mbox{Im}\;\Xi_{p,q}
\subset\ker\;\Theta_{p,q}$, and the proof of Proposition 2.4.8 is completed.

\hfill $\Box$

\noindent{\bf Proof of Lemma 2.4.1}

\vspace{1.5mm}

The `only if' part follows from the fact that $(2.2.10)$ is natural
with respect to $(f,\lambda)$. We shall prove the `if' part next.

First of all, we observe that for any global quotient $X=Y/G$,
there is a natural covering map $\Pi:Y\rightarrow X$ defined by
$(Id,1):(Y,\{1\})\rightarrow (Y,G)$, with the orbit map $\pi:Y
\rightarrow Y/G=X$ being the induced map between underlying spaces.
Moreover, for any base-point structure $\underline{o}=(o,U_o,\hat{o})$,
the injective homomorphism $\Pi_\ast:\pi_1(Y,\hat{o})\rightarrow
\pi_1(X,\underline{o})$ coincides with the one in $(2.2.10)$.

To define the map $f:Y\rightarrow Y^\prime$, we denote
$\Pi^\prime:Y^\prime\rightarrow X^\prime$ the canonical covering
map for $X^\prime$. Then the assumption for the ``if'' part of the
lemma can be rephrased as that there are base-point structures
$\underline{o}$, $\underline{o^\prime}$, such that $(\Phi\circ\Pi)_\ast
(\pi_1(Y,\hat{o}))\subset \Pi_\ast^\prime(\pi_1(Y^\prime,\hat{o}^\prime))$.
By Lemma 2.4.6, there is a unique map $f:(Y,\hat{o})\rightarrow
(Y^\prime,\hat{o}^\prime)$ such that $\Phi\circ\Pi=\Pi^\prime\circ f$.

As for the homomorphism $\lambda:G\rightarrow G^\prime$, we observe
that $\pi_1(X,\underline{o})\rightarrow G$ in $(2.2.10)$ is surjective
since $Y$ is path-connected, hence by the assumption that
$\Phi_\ast(\pi_1(Y,\hat{o}))\subset\pi_1(Y^\prime,\hat{o}^\prime)$,
$\Phi_\ast$ induces a homomorphism $\lambda:G\rightarrow G^\prime$.

It is easily seen that the action of $\pi_1(X,\underline{o})/
\pi_1(Y,\hat{o})$ on $Y$ through deck transformations coincides with
the action of $G$ on $Y$ under the isomorphism $\pi_1(X,\underline{o})
/\pi_1(Y,\hat{o})\cong G$ given by $(2.2.10)$. Hence $f$ is
$\lambda$-equivariant by virtue of $(2.4.10)$. Note that here we need
not to assume that $Y^\prime$ is connected, locally path-connected
because the action of $\pi_1(X^\prime,\underline{o^\prime})
/\pi_1(Y^\prime,\hat{o}^\prime)$ on $Y^\prime$ through deck transformations
is given a priori by the homomorphism $\pi_1(X^\prime,\underline{o^\prime})
/\pi_1(Y^\prime,\hat{o}^\prime)\rightarrow G^\prime$ in $(2.2.10)$ and
the action $G^\prime$ on $Y^\prime$.

It remains to show that $\Phi$ is defined by $(f,\lambda)$. To this end,
we first observe that an induced homomorphism $(\{f_a\},\{\lambda_{ba}\}):
\Gamma\{U_a\}\rightarrow\Gamma\{U_{a^\prime}^\prime\}$ of $(f,\lambda)$
has the form $f_a=\delta_a^{-1}\circ f|_{\widehat{U_a}}$ and $\lambda_{ba}
(g)=\delta_b^{-1}\lambda(g)\delta_a$ for a set $\{\delta_a\in G^\prime\}$,
where each $\widehat{U_a}$ is an open subset of $Y$, and $g\in T(U_a,
U_b)\subset G$. Secondly, we observe that a more concrete description
of a representative of the canonical covering map $\Pi:Y\rightarrow X$,
which is induced by $(Id,1):(Y,\{1\})\rightarrow (Y,G)$, may be
given as
$$
(\{V_{a,i}|i\in I_a\},\{U_a\},\{(\pi_{a,i},
\delta_{a,i})\},\{\rho_{(b,j)(a,i)}^s\}), \leqno (2.4.16)
$$
where (1) $\{U_a\}$ is a cover of $X$, (2) for each $a$,
$\{V_{a,i}|i\in I_a\}$ is the set of connected components
of $\pi^{-1}(U_a)$ in $Y$, with $V_{a,i_0}$ taken to be
$\widehat{U_a}$ for some index $i_0$, (3) $\pi_{a,i}:
V_{a,i}\rightarrow\widehat{U_a}$ is the homeomorphism induced
by an element $\delta_{a,i}\in G$, and (4) for any $V_{a,i}$,
$V_{b,j}$ with $V_{a,i}\cap V_{b,j}\neq\emptyset$, there
is a set $\{\rho_{(b,j)(a,i)}^s\in T(U_a,U_b)\}$
labeled by the set of components of $V_{a,i}\cap V_{b,j}$,
denoted by $\{s\}$, where each $\rho_{(b,j)(a,i)}^s=
\delta_{b,j}\delta_{a,i}^{-1}\in G$. For simplicity, we
shall require $\delta_{a,i_0}=1$ for each index $a$ without
loss of generality.

Now according to the construction in Lemma 2.4.6, the map $f:(Y,\hat{o})
\rightarrow (Y^\prime,\hat{o}^\prime)$ is defined as follows. First of
all, we fix a representative $\kappa^\prime$ of $\Pi^\prime:
(Y^\prime,\hat{o}^\prime)\rightarrow (X^\prime,\underline{o}^\prime)$
as described in $(2.4.16)$. Then we can choose a representative $\kappa$
for $\Pi:(Y,\hat{o})\rightarrow (X,\underline{o})$ as in $(2.4.16)$,
and a representative $\sigma$ of $\Phi:(X,\underline{o})\rightarrow
(X^\prime,\underline{o}^\prime)$, such that $\sigma\circ\kappa$ can be
lifted by $\kappa^\prime$ to $Y^\prime$, which is defined to be $f$.
In particular, $\sigma\circ\kappa=\kappa^\prime\circ f$. To fix the
notations, we write $\sigma=(\{f_a\},\{\lambda_{ba}\})$. We denote
the component of $(\pi^\prime)^{-1}(U_a^\prime)$ that contains
$f(V_{a,i})$ by $V_{a,i}^\prime$ for any $V_{a,i}$, and denote by
$\delta_{a,i}^\prime$ the element of $G^\prime$ that induces the
homeomorphism $\pi_{a,i}^\prime:V_{a,i}^\prime\rightarrow \widehat
{U_a^\prime}$. Note that $V_{a,i_0}^\prime$, which is the component
containing $f(V_{a,i_0})=f(\widehat{U_a})$, is not necessarily
$\widehat{U_a^\prime}$.

With these notational conventions understood, we can easily see that
$f_a=\delta_{a,i_0}^\prime\circ f|_{\widehat{U_a}}$ as a consequence
of $\sigma\circ\kappa=\kappa^\prime\circ f$ and the assumption
$\delta_{a,i_0}=1$. We set $\delta_a=(\delta_{a,i_0}^\prime)^{-1}$.
Then to show that $\sigma$ is induced by $(f,\lambda)$, we only need
to check
$$
\lambda_{ba}(g)=\delta_b^{-1}\lambda(g)\delta_a. \leqno (2.4.17)
$$
It suffices to check $(2.4.17)$ for two special cases: (a) $g=\rho_{(b,j)
(a,i)}^s$ for any indexes $i\in I_a,j\in I_b$, and (b) the indexes $a=b$
and $g\in G_{U_a}$.

For case (a) where $g=\rho_{(b,j)(a,i)}^s$ for some indexes $i\in I_a,
j\in I_b$, we shall prove the relation $(\delta_{a,i_0}^\prime)^{-1}
\delta_{a,i}^\prime=\lambda(\delta_{a,i})$ for any index $a$ and $i
\in I_a$, which implies case (a) because $\rho_{(b,j)(a,i)}^s=\delta_{b,j}
\delta_{a,i}^{-1}$, and $\lambda_{ba}(\rho_{(b,j)(a,i)}^s)=\delta_{b,j}^\prime
(\delta_{a,i}^\prime)^{-1}$ as a consequence of $\sigma\circ\kappa=
\kappa^\prime\circ f$. To see $(\delta_{a,i_0}^\prime)^{-1}
\delta_{a,i}^\prime=\lambda(\delta_{a,i})$, we pick paths $\gamma_0$,
$\gamma$ in $Y$ satisfying $\gamma_0(0)=\gamma(0)=\hat{o}$, $\gamma_0(1)
\in\widehat{U_a}$, and $\gamma(1)\in V_{a,i}$ with $\gamma_0(1)=
\delta_{a,i}\cdot\gamma(1)$. Then $\kappa\circ\gamma_0$ may compose with
$\nu(\kappa\circ\gamma)$ to define a loop $u$ in $(X,\underline{o})$.
The lifting of $u$ to $(Y,\hat{o})$ by $\kappa$ is easily seen to be the
path $\gamma_0\#\nu(\delta_{a,i}\circ\gamma)$ whose terminal point is
$\delta_{a,i}\cdot\hat{o}$. Hence the image of $[u]$ under the homomorphism
$\pi_1(X,\underline{o})\rightarrow G$ in $(2.2.10)$ is $\delta_{a,i}$. On
the other hand, the push-forward $\sigma\circ u$ is lifted to
$(f\circ\gamma_0)\#\nu((\delta_{a,i_0}^\prime)^{-1}\delta_{a,i}^\prime
\circ (f\circ\gamma))$ by $\kappa^\prime$, whose terminal point is
$(\delta_{a,i_0}^\prime)^{-1}\delta_{a,i}^\prime\cdot\hat{o}^\prime$. Hence
the image of $\Phi_\ast([u])$ under $\pi_1(X^\prime,\underline{o^\prime})
\rightarrow G^\prime$ in $(2.2.10)$ is $(\delta_{a,i_0}^\prime)^{-1}
\delta_{a,i}^\prime$. This gives $(\delta_{a,i_0}^\prime)^{-1}
\delta_{a,i}^\prime=\lambda(\delta_{a,i})$ by the definition of $\lambda$.

For case (b) where the indexes $a=b$ and $g\in G_{U_a}$, we pick a
path $\gamma_0$ in $Y$ with $\gamma_0(0)=\hat{o}$, $\gamma_0(1)\in
\widehat{U_a}$, and pick a path $\gamma$ in $\widehat{U_a}$ connecting
$\gamma_0(1)$ to $g\cdot\gamma_0(1)$. Note that if we set $x=
\pi(\gamma_0(1))$, $\underline{x}=(x,U_a,\gamma_0(1))$, then $(\gamma,g)$
defines a loop $v$ in $(X,\underline{x})$ by Lemma 2.2.2. We set
$u=(\kappa\circ\gamma_0)\# v\#\nu(\kappa\circ\gamma_0)$, which is a
loop in $\Omega(X,\underline{o})$. The lifting of $u$ by $\kappa$ is
$\gamma_0\#\gamma\#\nu(g\circ\gamma_0)$, whose terminal point is
$g\cdot\hat{o}$. Hence the image of $[u]$ under the homomorphism
$\pi_1(X,\underline{o})\rightarrow G$ in $(2.2.10)$ is $g$.
On the other hand, the push-forward $\sigma\circ u$ is lifted by
$\kappa^\prime$ to the path $\ell(\sigma\circ u)=(f\circ\gamma_0)\#
\gamma^\prime\#\nu(Ad((\delta_{a,i_0}^\prime)^{-1})(\lambda_a(g))
\circ (f\circ\gamma_0))$ for some path $\gamma^\prime$ connecting
$f(\gamma_0(1))$ to $Ad((\delta_{a,i_0}^\prime)^{-1})(\lambda_a(g))
\cdot f(\gamma_0(1))$. The terminal point of $\ell(\sigma\circ u)$ is
$Ad((\delta_{a,i_0}^\prime)^{-1})(\lambda_a(g))\cdot\hat{o}^\prime$,
which implies that the image of $\Phi_\ast([u])$ under
$\pi_1(X^\prime,\underline{o^\prime})\rightarrow G^\prime$ in $(2.2.10)$
is $Ad((\delta_{a,i_0}^\prime)^{-1})(\lambda_a(g))$. Hence $\lambda(g)=
Ad((\delta_{a,i_0}^\prime)^{-1})(\lambda_a(g))$, from which $(2.4.17)$
for case (b) can be easily deduced.

\hfill $\Box$

\sectioni{An analog of CW-complex theory}

\subsection{Construction of mapping cylinders}

In this subsection we extend the mapping cylinder construction in the
ordinary homotopy theory (cf. e.g. \cite{Sw}) to the orbispace category.
Recall that for any map $f:Y\rightarrow X$ between topological spaces,
the mapping cylinder of $f$, denoted by $M_f$, is the topological space
obtained by identifying $(y,1)$ with $f(y)$ in the disjoint union
$Y\times I\bigsqcup X$ for all $y\in Y$. We set $[y,t]$ for the image
of $(y,t)\in Y\times I\bigsqcup X$ in $M_f$. There are embeddings
$i:Y\rightarrow M_f$, $j:X\rightarrow M_f$ given by $i(y)=[y,0]$ and
$j(x)=x$, realizing $Y$, $X$ as subspaces of $M_f$. The space $X$ is
a strong deformation retract of $M_f$, with the canonical retraction
$r:M_f\rightarrow X$ where $r([y,t])=f(y)$ and $r(x)=x$. The homotopy
$H:M_f\times [0,1]\rightarrow M_f$ between $j\circ r$ and $Id_{M_f}$
is given by $H([y,t],s)=[y,1-(1-t)s]$ and $H(x,s)=x$. Finally, given
any homotopy $F:Y\times [0,1]\rightarrow X$ between $f_1,f_2:Y
\rightarrow X$, there are canonical maps $\phi_F:M_{f_1}\rightarrow
M_{f_2}$, $\psi_F:M_{f_2}\rightarrow M_{f_1}$, such that $\psi_F\circ
\phi_F$ and $\phi_F\circ\psi_F$ are canonically homotopic to $Id_{M_{f_1}}$
and $Id_{M_{f_2}}$ respectively, relative to $X$ and $Y$. For instance,
we may define $\phi_F,\psi_F$ by
$$
\begin{array}{l}
\phi_F([y,t])=\left\{\begin{array}{ll}
[y,t] & 0\leq t\leq \frac{1}{2}\\
F(y,2-2t) & \frac{1}{2}\leq t\leq 1
\end{array} \right. \hspace{2mm}
\psi_F([y,t])=\left\{\begin{array}{ll}
[y,t] & 0\leq t\leq \frac{1}{2}\\
F(y,2t-1) & \frac{1}{2}\leq t\leq 1
\end{array} \right. \\
\phi_F(x)=x \hspace{5.5cm} \psi_F(x)=x
\end{array} \leqno (3.1.1)
$$
The canonical homotopy between $\psi_F\circ\phi_F$ and $Id_{M_{f_1}}$
may be taken to be $H:M_{f_1}\times [0,1]\rightarrow M_{f_1}$ where
$$
\begin{array}{l}
H([y,t],s)=\left\{\begin{array}{ll}
[y,(4-3s)t], & 0\leq t\leq \frac{1}{4-3s}\\
F(y,(4-3s)t-1), & \frac{1}{4-3s}\leq t\leq \frac{2-s}{4-3s}\\
F(y,\frac{1}{2}(4-3s)(1-t)), & \frac{2-s}{4-3s}\leq t\leq 1
\end{array} \right. \\
H(x,s)=x.
\end{array} \leqno (3.1.2)
$$
In the same vein, we obtain the canonical homotopy between
$\phi_F\circ\psi_F$ and $Id_{M_{f_2}}$.

We derive two technical lemmas first.

\begin{lem}
Let $\Phi:X\rightarrow X^\prime$ be a map, where $X^\prime=Y^\prime
/G^\prime$ is a global quotient, and $X$ is connected, locally
path-connected and semi-locally 1-connected. Then there exists a
connected, locally path-connected space $Y$ with a discrete group
action of $G$, such that $X=Y/G$. Moreover, $\Phi:X\rightarrow X^\prime$
is represented by a pair $(f,\lambda):(Y,G)\rightarrow
(Y^\prime,G^\prime)$ where $f$ is $\lambda$-equivariant. Such a $(Y,G)$
is uniquely determined up to an isomorphism, and $(f,\lambda)$ is unique up
to conjugation by an element $g\in G^\prime$, i.e., $(f,\lambda)\mapsto
(g\circ f,Ad(g)\circ\lambda)$.
\end{lem}

\pf
We take a based version of $\Phi$ with respect to some base-point
structures $\underline{o}$, $\underline{o^\prime}$ of $X$, $X^\prime$
respectively. Since $X$ is connected, locally path-connected
and semi-locally 1-connected, there is a covering space $\Pi:
(Y,\underline{q})\rightarrow (X,\underline{o})$ such that
$\Pi_\ast(\pi_1(Y,\underline{q}))=(\Phi_\ast)^{-1}
(\pi_1(Y^\prime,\hat{o}^\prime))$ by Proposition 2.4.7. Moreover, by
Lemma 2.4.6, there is a map $\Psi:(Y,\underline{q})\rightarrow
(Y^\prime,\hat{o}^\prime)$ satisfying $\Phi\circ\Pi=\Pi^\prime\circ
\Psi$, where $\Pi^\prime:(Y^\prime,\hat{o}^\prime)\rightarrow (X^\prime,
\underline{o^\prime})$ is the canonical covering map associated
to the global quotient $X^\prime=Y^\prime/G^\prime$.

Now recall that $\Phi$ satisfies the assumption in {\bf Convention}
in Introduction, so does the map $\Psi:(Y,\underline{q})\rightarrow
(Y^\prime,\hat{o}^\prime)$. This implies that $Y$ is actually a
topological space. We rename $\Psi$ by $f$, and set
$G=\pi_1(X,\underline{o})/\Pi_\ast(\pi_1(Y,\underline{q}))$.
Then $G$ acts on $Y$ through deck transformations, and
there is an injective homomorphism $\lambda:G\rightarrow G^\prime$
such that $f:Y\rightarrow Y^\prime$ is $\lambda$-equivariant.

We next verify that $X$ is isomorphic to the global quotient $Y/G$
and $\Pi:Y\rightarrow X$ is isomorphic to the canonical covering map.
To this end, we pick a representative of $\Pi:(Y,\underline{q})
\rightarrow (X,\underline{o})$, which can be described in the
following form:
$$
(\{V_{\alpha,i}|i\in I_\alpha\},\{U_\alpha\},\{\pi_{\alpha,i}\},
\{\rho_{(\beta,j)(\alpha,i)}^s\}), \leqno (3.1.3)
$$
where (1) $\{U_\alpha\}$ is a cover of $X$ by local charts, (2)
for each index $\alpha$, $\{V_{\alpha,i}|i\in I_\alpha\}$ is the set
of connected components of $\pi^{-1}(U_\alpha)\subset Y$, (3) each
$\pi_{\alpha,i}:V_{\alpha,i}\rightarrow\widehat{U_\alpha}$ is a
homeomorphism, (4) for any $V_{\alpha,i},V_{\beta,j}$ with
$V_{\alpha,i}\cap V_{\beta,j}\neq\emptyset$, there is a set
$\{\rho_{(\beta,j)(\alpha,i)}^s\in T(U_\alpha,U_\beta)\}$
labeled by the set $\{V_{(\alpha,i)(\beta,j)}^s\}$ of connected
components of $V_{\alpha,i}\cap V_{\beta,j}$, such that the following
equations are satisfied,
$$
\pi_{\beta,j}=\rho_{(\beta,j)(\alpha,i)}^s\circ\pi_{\alpha,i}
\mbox{ on } V_{(\alpha,i)(\beta,j)}^s, \mbox{ and }
\rho_{(\gamma,k)(\beta,j)}^t\circ\rho_{(\beta,j)(\alpha,i)}^s
(\pi_{\alpha,i}({\bf a}))=\rho_{(\gamma,k)(\alpha,i)}^r,
\leqno (3.1.4)
$$
where ${\bf a}$ is any connected component of $V_{(\alpha,i)(\beta,j)}^s
\cap V_{(\beta,j)(\gamma,k)}^t$ that is contained in
$V_{(\alpha,i)(\gamma,k)}^r$.

We claim: (1) Set $G_{\alpha,i}=\{g\in G\mid g\cdot V_{\alpha,i}=
V_{\alpha,i}\}$, then there is an isomorphism $\lambda_{\alpha,i}:
G_{\alpha,i}\rightarrow G_{U_\alpha}$ such that $\pi_{\alpha,i}$ is
$\lambda_{\alpha,i}$-equivariant. (2) For any $V_{\alpha,i},
V_{\alpha,j}$, there is a $g_{ji}^\alpha\in G$ satisfying
$\pi_{\alpha,j}\circ g_{ji}^\alpha=\pi_{\alpha,i}$ and
$g_{ki}^\alpha=g^\alpha_{kj}g_{ji}^\alpha$. To sketch a proof,
we pick, for each $V_{\alpha,i}$, a path $\gamma_i^\alpha$ in $Y$
satisfying $\gamma_i^\alpha(0)=\hat{q}$, $\gamma_i^\alpha(1)\in
V_{\alpha,i}$ and $\pi_{\alpha,i}(\gamma_i^\alpha(1))=
\pi_{\alpha,j}(\gamma_j^\alpha(1))$. We introduce base-point structures
$\underline{x_\alpha}=(x_\alpha,U_\alpha,\hat{x}_\alpha)$ where $x_\alpha
=\pi(\gamma_i^\alpha(1))$, $\hat{x}_\alpha=\pi_{\alpha,i}
(\gamma_i^\alpha(1))$, and denote by $\bar{\gamma_i}^\alpha$ the push-down
path $\Pi\circ\gamma_i^\alpha\in P(X,\underline{o},\underline{x_\alpha})$.
Then the isomorphism $\lambda_{\alpha,i}:G_{\alpha,i}\rightarrow
G_{U_\alpha}$ is defined as follows. For any $g\in G_{U_\alpha}$,
we pick a path $\gamma$ in $\widehat{U_\alpha}$ connecting $\hat{x}_\alpha$
to $g\cdot \hat{x}_\alpha$, and denote by ${\gamma}_g$ the loop in
$(X,\underline{x_\alpha})$ defined by $(\gamma,g)$. The inverse
of $\lambda_{\alpha,i}$ is defined by sending $g$ to the image of
$[\bar{\gamma_i}^\alpha\# {\gamma}_g\#\nu(\bar{\gamma_i}^\alpha)]
\in\pi_1(X,\underline{o})$ under the homomorphism $\pi_1(X,\underline{o})
\rightarrow G=\pi_1(X,\underline{o})/\Pi_\ast(\pi_1(Y,\underline{q}))$.
As for the set of elements $\{g_{ji}^\alpha\}$, we define $g_{ji}^\alpha$
to be the image of $[\bar{\gamma_j}^\alpha\#\nu(\bar{\gamma_i}^\alpha)]
\in\pi_1(X,\underline{o})$ under the homomorphism $\pi_1(X,\underline{o})
\rightarrow G$. With these definitions, the verification of the claim
is straightforward, hence we leave it to the reader.

For each $\alpha$, we fix an index $i_0\in I_\alpha$, and set $W_\alpha
=V_{\alpha,i_0}/G_{\alpha,i_0}$. Then an isomorphism between $X$ and
$Y/G$ is defined by $\sigma=(\{\pi_{\alpha,i_0}^{-1}\},
\{\eta_{\beta\alpha}\}):\Gamma\{U_\alpha\}\rightarrow\Gamma\{W_\alpha\}$,
where $\eta_{\beta\alpha}$ is the mapping uniquely determined by
the conditions $\eta_{\alpha\alpha}=\lambda_{\alpha,i_0}^{-1}$
and $\eta_{\beta\alpha}(\rho_{(\beta,j)(\alpha,i)}^s)=g_{j_0j}^\beta
g_{ii_0}^\alpha$. Under this isomorphism, the covering map $\Pi:
Y\rightarrow X$ is isomorphic to the canonical one associated to
the global quotient $Y/G$.

Finally, by Lemma 2.4.1, $\Phi:X\rightarrow X^\prime$ is represented
by $(f,\lambda):(Y,G)\rightarrow (Y^\prime,G^\prime)$. By Lemma 3.1.2
in \cite{C1}, $(f,\lambda)$ is uniquely determined up to conjugation.
To see that $(Y,G)$ is uniquely determined up to an isomorphism, we
simply observe that the injectivity of $\lambda:G\rightarrow G^\prime$
implies $\pi_1(Y,\hat{o})=(\Phi_\ast)^{-1}(\pi_1(Y^\prime,\hat{o}^\prime))$
for appropriate base-point structures $\underline{o}$, $\underline{o^\prime}$
of $X$, $X^\prime$. In other words, $\pi_1(Y,\hat{o})$ is uniquely determined,
hence so is the covering space $Y$ up to an isomorphism by Lemma 2.4.6.

\hfill $\Box$

Let $Y$ be a locally path-connected, semi-locally 1-connected
orbispace. For any map $\Phi:Y\rightarrow X$, we consider the set
$\V(\Phi)$ of connected open subsets $V$ of $Y$ such that $V\subset
\phi^{-1}(U)$, where $\phi:Y\rightarrow X$ is the induced map of $\Phi$
between underlying spaces, and $U\in\U$ is a local chart on
$X$. We apply the preceding lemma to $\Phi|_V:V\rightarrow U$,
so that for each $V\in\V(\Phi)$, there is a $(\widetilde{V},G^V)$
such that the subspace $V$ is isomorphic to the global quotient
$\widetilde{V}/G^V$, under which $\Phi|_V$ is represented by a
pair $(f_V,\lambda_V):(\widetilde{V},G^V)\rightarrow (\widehat{U},G_U)$
where $f_V$ is $\lambda_V$-equivariant. Set $\pi^V:\widetilde{V}
\rightarrow\widetilde{V}/G^V=V$. We remark that
\begin{itemize}
\item [{(1)}] $(\widetilde{V},G^V,\pi^V)$ depends only on $V$. This
is because, as we have seen in the proof of Lemma 3.1.1, the isomorphism
class of $(\widetilde{V},G^V)$ is determined by $(\Phi|_V)_\ast^{-1}
(\pi_1(\widehat{U}))$, and on the other hand, for any $U_1\subset U_2$,
$(\Phi|_V)_\ast^{-1}(\pi_1(\widehat{U_1}))=(\Phi|_V)_\ast^{-1}
(\pi_1(\widehat{U_2}))$. The last identity follows from the fact
that the inclusion $U_1\subset U_2$ induces an injective homomorphism
from $G_{U_1}$ to $G_{U_2}$ so that by $(2.2.10)$, the inverse image of
$\pi_1(\widehat{U_2})$ is $\pi_1(\widehat{U_1})$ under $\pi_1(U_1)
\rightarrow \pi_1(U_2)$.
\item [{(2)}] If $V_1,V_2\in\V(\Phi)$ such that $V_1\subset V_2$,
then $(\widetilde{V_1},G^{V_1},\pi^{V_1})$ is isomorphic to an induced
one from $(\widetilde{V_2},G^{V_2},\pi^{V_2})$. This is because
if we denote by $I:V_1\rightarrow V_2$ the inclusion as a subspace,
then $\Phi|_{V_1}=\Phi|_{V_2}\circ I$, and hence $I_\ast^{-1}\circ
(\Phi|_{V_2})_\ast^{-1}(\pi_1(\widehat{U}))=(\Phi|_{V_1})_\ast^{-1}
(\pi_1(\widehat{U}))$ where $V_1\subset V_2\subset \phi^{-1}(U)$,
so that $I_\ast^{-1}(\pi_1(\widetilde{V_2}))=
\pi_1(\widetilde{V_1})$ under $I_\ast:\pi_1(V_1)\rightarrow\pi_1(V_2)$.
By Lemma 2.4.1, $I$ is represented by a pair $(f,\lambda):
(\widetilde{V_1},G^{V_1})\rightarrow (\widetilde{V_2},G^{V_2})$
where $\lambda$ is injective and $f$ is a $\lambda$-equivariant open
embedding.
\item [{(3)}] When $V$ is a local chart on
$Y$ such that $\Phi|_V$ is represented by a pair $(f,\rho):
(\widehat{V},G_V)\rightarrow (\widehat{U},G_U)$, then $(\widehat{V},
G_V,\pi_V)$ is isomorphic to $(\widetilde{V},G^V,\pi^V)$. For
instance, suppose $\Phi$ is represented by $(\{f_\alpha\},
\{\rho_{\beta\alpha}\}):\Gamma\{V_\alpha\}\rightarrow
\Gamma\{U_{\alpha^\prime}\}$, then if a local chart $V\subset V_\alpha$
for some index $\alpha$, then $V\in\V(\Phi)$ and $(\widehat{V},
G_V,\pi_V)$ is isomorphic to $(\widetilde{V},G^V,\pi^V)$.
\end{itemize}

\begin{lem}
There is a set $\T(\Phi)=\{T(V_1,V_2)\mid V_1,V_2\in\V(\Phi),
\mbox{ s.t. } V_1\cap V_2\neq\emptyset\}$, which together with
the atlas of local charts $\{(\widetilde{V},G^V,\pi^V)\mid
V\in\V(\Phi)\}$ defines an orbispace structure on $Y$ that is
equivalent to the original one. Moreover, given any cover
$\{V_\alpha\}\subset\V(\Phi)$ of $Y$, where $V_\alpha\subset
\phi^{-1}(U_\alpha)$ for some local chart $U_\alpha$ on $X$,
there is a homomorphism $(\{f_\alpha\},\{\lambda_{\beta\alpha}\}):
\Gamma\{V_\alpha\}\rightarrow\Gamma\{U_\alpha\}$ whose equivalence
class is the given map $\Phi$, in which $(f_\alpha,\lambda_\alpha)
=(f_{V_\alpha},\lambda_{V_\alpha}):(\widetilde{V_\alpha},
G^{V_\alpha})\rightarrow (\widehat{U_\alpha},G_{U_\alpha})$.
\end{lem}

\pf
For any $V\in\V(\Phi)$, the subspace structure of $V\subset Y$ is
given by an atlas of local charts $\N(V)=\{(\widehat{W},G_W,\pi_W)\}$
and a set $\T(V)=\{T(W_1,W_2)\}$, where $W\subset V$ is a local chart
on $Y$. Since the subspace $V$ is isomorphic to $\widetilde{V}/G^V$,
we may further require that for each $W\in\N(V)$, $W\in\V(\Phi)$
and $(\widetilde{W},G^W,\pi^W)=(\widehat{W},G_W,\pi_W)$. Moreover,
a set $T(W,V)$ is defined, such that each $\xi\in T(W,V)$ is associated
with a $(\phi_\xi,\lambda_\xi):(\widehat{W},G_W)\rightarrow
(\widetilde{V},G^V)$ where $\lambda_\xi$ is injective and $\phi_\xi$
is a $\lambda_\xi$-equivariant open embedding, and for any $g\in G^V$,
$g\circ\xi\in T(W,V)$ and $(\phi_{g\circ\xi},\lambda_{g\circ\xi})
=(\phi_g,Ad(g))\circ (\phi_\xi,\lambda_\xi)$.

Let $V_1,V_2\in\V(\Phi)$ such that $V_1\cap V_2\neq\emptyset$.
The set $T(V_1,V_2)$ is defined as follows. First of all,
let $V$ be any connected component of $V_1\cap V_2$, and
$W\in\N(V_1)\cap\N(V_2)$, $W\subset V$.
To each $(\xi_1,\xi_2)\in T(W,V_1)\times T(W,V_2)$, we assign a pair
$\chi^W_V(\xi_1,\xi_2)=(f,\lambda)$, where $f:\widetilde{Z_1}
\rightarrow\widetilde{Z_2}$ is a homeomorphism from a connected
component of $(\pi^{V_1})^{-1}(V)\subset\widetilde{V_1}$ to a connected
component of $(\pi^{V_2})^{-1}(V)\subset\widetilde{V_2}$ satisfying
$\pi^{V_1}=\pi^{V_2}\circ f$, and $\lambda:G_{\widetilde{Z_1}}\rightarrow
G_{\widetilde{Z_2}}$ is an isomorphism from the subgroup of $G^{V_1}$
fixing $\widetilde{Z_1}$ to the subgroup of $G^{V_2}$ fixing
$\widetilde{Z_2}$, such that $f$ is $\lambda$-equivariant. To this
end, for $i=1,2$, we let $\widetilde{Z_i}$ be the connected component
of $(\pi^{V_i})^{-1}(V)\subset\widetilde{V_i}$ that contains
$\mbox{Range }(\phi_{\xi_i})$. Suppose $V\subset\phi^{-1}(U)$ for
some local chart $U$ on $X$. We pick base-point structures $\underline{q}
=(q,W,\hat{q})$ and $\underline{p}=(p,U,\hat{p})$ such that
$\Phi:(Y,\underline{q})\rightarrow (X,\underline{p})$. Now for $i=1,2$
we set $\hat{q_i}=\phi_{\xi_i}(\hat{q})$, and denote by $\Pi_i:
(\widetilde{Z_i},\hat{q_i})\rightarrow (V,\underline{q})$ the canonical
covering map associated to the isomorphism $V\cong \widetilde{Z_i}/
G_{\widetilde{Z_i}}$ which is defined by $\phi_{\xi_i}^{-1}$ at the
base-point structures. Then we have $(\Pi_1)_\ast(\pi_1(\widetilde{Z_1},
\hat{q_1}))=(\Phi|_V)_\ast^{-1}(\pi_1(\widehat{U},\hat{p}))=
(\Pi_2)_\ast(\pi_1(\widetilde{Z_2},\hat{q_2}))$ in $\pi_1(V,\underline{q})$.
By Lemma 2.4.6, there is a unique $f:(\widetilde{Z_1},\hat{q_1})\rightarrow
(\widetilde{Z_2},\hat{q_2})$ satisfying $\Pi_1=\Pi_2\circ f$. Furthermore,
by Proposition 2.4.8, there is a natural isomorphism $\lambda:
G_{\widetilde{Z_1}}\rightarrow G_{\widetilde{Z_2}}$ such that $f$ is
$\lambda$-equivariant. We define $\chi^W_V(\xi_1,\xi_2)=(f,\lambda)$.
It is easy to see that (1) $\chi^W_V(\xi_1,\xi_2)$ is independent of
the choices on $U,\underline{q}$ and $\underline{p}$, (2)
$\chi^W_V(\xi_1,\xi_2)=(\phi_{\xi_2},\lambda_{\xi_2})\circ
(\phi_{\xi_1}^{-1},\lambda_{\xi_1}^{-1})$ when restricted to
the domain of the latter, and (3) $\chi^W_V(\xi_1,\xi_2)=
\chi^W_V(g\circ\xi_1,\lambda(g)\circ\xi_2)$, $\forall g\in
G_{\widetilde{Z_1}}$.

Secondly, in $\bigsqcup_{W\in\N(V_1)\cap\N(V_2),W\subset V}T(W,V_1)\times
T(W,V_2)$ we introduce an equivalence relation $\sim$ generated as follows:
(a) for any $(\xi_1,\xi_2),(\eta_1,\eta_2)\in T(W,V_1)\times T(W,V_2)$,
$(\xi_1,\xi_2)\sim (\eta_1,\eta_2)$ if $\eta_1=g\circ\xi_1,
\eta_2=\lambda(g)\circ\xi_2$ for some $g\in\mbox{Domain }(\lambda)$,
where $\lambda$ is given in $\chi_V^W(\xi_1,\xi_2)=(f,\lambda)$, and
(b) for any
$(\xi_1,\xi_2)\in T(W,V_1)\times T(W,V_2)$, $(\eta_1,\eta_2)\in
T(W^\prime,V_1)\times T(W^\prime,V_2)$ where $W^\prime\subset W$,
$(\xi_1,\xi_2)\sim (\eta_1,\eta_2)$ if $\eta_1=\xi_1\circ\epsilon,
\eta_2=\xi_2\circ\epsilon$ for some $\epsilon\in T(W^\prime,W)$.
We define $T_V(V_1,V_2)=\bigsqcup_{W\in\N(V_1)\cap\N(V_2),W\subset V}
T(W,V_1)\times T(W,V_2)/\sim$, and for each $\xi=[(\xi_1,\xi_2)]$
where $(\xi_1,\xi_2)\in T(W,V_1)\times T(W,V_2)$, define $(\phi_\xi,
\lambda_\xi)=\chi_V^W(\xi_1,\xi_2)$.

Finally, we let $T(V_1,V_2)$ be the disjoint union of $T_V(V_1,V_2)$
for all components $V$ of $V_1\cap V_2$. Note that with this definition,
$T(V,V)$, $\forall V\in\V(\Phi)$, is naturally identified with $G^V$.

The composition is defined as follows. Let $\xi\in
T(V_1,V_2)$, $\eta\in T(V_2,V_3)$, and $x\in\phi_\xi^{-1}
(\mbox{Domain }(\phi_\eta))$. We pick a $W\in\N(V_1)\cap\N(V_2)\cap
\N(V_3)$ such that $\pi^{V_1}(x)\in W$. Then we can write
$\xi=[(\xi_1,\xi_2)]$ for some $(\xi_1,\xi_2)\in T(W,V_1)\times
T(W,V_2)$, $\eta=[(\eta_2,\eta_3)]$ for some $(\eta_2,\eta_3)\in
T(W,V_2)\times T(W,V_3)$ such that $\mbox{Range }(\phi_{\xi_2})=
\mbox{Range }(\phi_{\eta_2})$. The last condition allows us to
modify $(\xi_1,\xi_2)$ without changing its class $[(\xi_1,\xi_2)]$
but to further satisfy $\xi_2=\eta_2$. With these arranged, we define
$\eta\circ\xi(x)=[(\xi_1,\eta_3)]$. As for the inverse,
we simply define $\xi^{-1}=[(\xi_2,\xi_1)]$ if $\xi=[(\xi_1,\xi_2)]$.

The verification that $\V(\Phi)=\{(\widetilde{V},G^V,\pi^V)\}$,
$\T(\Phi)=\{T(V_1,V_2)\}$ satisfy the conditions in Proposition
2.1.1 of \cite{C1} is straightforward, which is left to the reader.
Thus they define an orbispace structure on $Y$ by Proposition 2.1.1
of \cite{C1}, which is clearly `equivalent' to the original
one on $Y$ in the sense of Remark 2.1.2 (5) in \cite{C1}.

Finally, given any cover $\{V_\alpha\}\subset\V(\Phi)$ of $Y$,
where $V_\alpha\subset\phi^{-1}(U_\alpha)$ for some local chart
$U_\alpha$ on $X$, the existence of a homomorphism $(\{f_\alpha\},
\{\lambda_{\beta\alpha}\}):\Gamma\{V_\alpha\}\rightarrow
\Gamma\{U_\alpha\}$, whose equivalence class is the given map $\Phi$
and in which $(f_\alpha,\lambda_\alpha)=(f_{V_\alpha},\lambda_{V_\alpha}):
(\widetilde{V_\alpha},G^{V_\alpha})\rightarrow (\widehat{U_\alpha},
G_{U_\alpha})$, is the content of Lemma 3.1.3 in \cite{C1}.

\hfill $\Box$

Let $\Gamma^0$ be an orbispace structure on $X$, which is given by
the data $\{(\widehat{U}_i^0,G_{U_i}^0,\pi_{U_i}^0)\}$ and $\T^0=
\{T^0(U_i,U_j)\}$ as described in Proposition 2.1.1 of \cite{C1}.
Suppose for each $U_i$, there is a triple $(\widehat{U_i},G_{U_i},
\pi_{U_i})$, where $\widehat{U_i}$ is not connected in general,
such that $\pi_{U_i}:\widehat{U_i}\rightarrow U_i$ induces
a homeomorphism $\widehat{U_i}/G_{U_i}\cong U_i$, and $(\widehat{U}_i^0,
G_{U_i}^0,\pi_{U_i}^0)$ is obtained from $(\widehat{U_i},G_{U_i},
\pi_{U_i})$ by restricting to a connected component $\widehat{U}_i^0$
of $\widehat{U_i}$. Then there is canonically an orbispace structure
$\Gamma$ on $X$, with $\Gamma^0\subset\Gamma$ being an equivalence,
and $\{(\widehat{U_i},G_{U_i},\pi_{U_i})\}$ being the atlas of local charts.
As for the set $\T=\{T(U_i,U_j)\}$ of $\Gamma$, each $T(U_i,U_j)$
is the orbit space of $G_{U_i}\times T^0(U_i,U_j)\times G_{U_j}$ modulo
the action $h\cdot (g_i,\xi_0,g_j)=(hg_i,\xi_0,\lambda_{\xi_0}(h)g_j),
\;\forall h\in\mbox{Domain }(\phi_{\xi_0})$. (It is instructive to think
$(g_i,\xi_0,g_j)$ as $g_j^{-1}\circ\xi_0\circ g_i$.)

Moreover, suppose $(\{f_\alpha^0\},\{\rho^0_{\beta\alpha}\}):
\Gamma^0\{U_\alpha\}\rightarrow\Gamma\{U_{\alpha^\prime}^\prime\}$ is
a homomorphism of groupoids, where each $\rho_\alpha^0:G_{U_\alpha}^0
\rightarrow G_{U_\alpha^\prime}$ is assumed to be injective. Then by
replacing $\Gamma^0$ with $\Gamma$ on $X$, we may canonically replace
$(\{f_\alpha^0\},\{\rho^0_{\beta\alpha}\})$ with a homomorphism
$(\{f_\alpha\},\{\rho_{\beta\alpha}\}):\Gamma\{U_\alpha\}\rightarrow
\Gamma\{U_{\alpha^\prime}^\prime\}$ such that (1) each $\rho_\alpha:
G_{U_\alpha}\rightarrow G_{U_\alpha^\prime}$ is isomorphic, and (2)
$f_\alpha^0,\rho^0_{\beta\alpha}$ are obtained by restricting $f_\alpha,
\rho_{\beta\alpha}$ to $\widehat{U}_\alpha^0,T^0(U_\alpha,U_\beta)$
respectively. It is done as follows. We define $\widehat{U_\alpha}=
(G_{U_\alpha^\prime}\times\widehat{U}_\alpha^0)/G_{U_\alpha}^0$ where
the action is given by $g^0\cdot (g^\prime,x)=(g^\prime\rho_\alpha^0
(g^0)^{-1},g^0\cdot x)$, we define $G_{U_\alpha}=G_{U_\alpha^\prime}$
with a left action on $\widehat{U_\alpha}$ induced by $g\cdot (g^\prime,x)
=(gg^\prime,x)$, and we define $\pi_{U_\alpha}:\widehat{U_\alpha}
\rightarrow U_\alpha$ by $[(g^\prime,x)]\mapsto\pi_{U_\alpha}^0(x)$,
which induces $\widehat{U_\alpha}/G_{U_\alpha}\cong U_\alpha$.
Note that $(\widehat{U}_\alpha^0,G_{U_\alpha}^0,\pi_{U_\alpha}^0)$
is obtained from $(\widehat{U_\alpha},G_{U_\alpha},\pi_{U_\alpha})$
under the mapping $(x,g^0)\mapsto ([(1,x)],\rho_\alpha^0(g^0))$.
The pair of maps $(f_\alpha,\rho_\alpha):(\widehat{U_\alpha},G_{U_\alpha})
\rightarrow (\widehat{U_\alpha^\prime},G_{U_\alpha^\prime})$ is defined by
$f_\alpha([(g^\prime,x)])=g^\prime\cdot f_\alpha^0(x)$ and $\rho_\alpha=Id$.
Its restriction to $(\widehat{U}_\alpha^0,G_{U_\alpha}^0)$ is clearly
$(f_\alpha^0,\rho_\alpha^0)$. The mappings $\rho_{\beta\alpha}:
T(U_\alpha,U_\beta)\rightarrow T(U_\alpha^\prime,U_\beta^\prime)$ are
defined by $\rho_{\beta\alpha}(g_\beta^{-1}\circ\xi_0\circ g_\alpha)
=\rho_\beta(g_\beta)^{-1}\circ\rho_{\beta\alpha}^0(\xi_0)\circ
\rho_{\alpha}(g_\alpha)=g_\beta^{-1}\circ\rho_{\beta\alpha}^0(\xi_0)
\circ g_\alpha$, which clearly satisfy $\rho_{\beta\alpha}|_{T^0(U_\alpha,
U_\beta)}=\rho_{\beta\alpha}^0$. One can check directly that
$(\{f_\alpha\},\{\rho_{\beta\alpha}\})$ is indeed a homomorphism. We
leave the details to the reader.

\vspace{2.5mm}

\noindent{\bf Proof of Proposition 1.5}

\vspace{1.5mm}

By Lemma 3.1.2, we may assume that the original orbispace structure on
$Y$ is given by $(\V(\Phi),\T(\Phi))$ for simplicity without loss of
generality.

Let $\{V_\alpha\}$ be the subset of $\V(\Phi)$ which consists of connected
components of $\phi^{-1}(U)$ for all local chart $U$ on $X$. For each
index $\alpha$, we pick a $U_\alpha$ such that $V_\alpha$ is a component
of $\phi^{-1}(U_\alpha)$. By Lemma 3.1.2, there is a homomorphism, denoted
by $(\{f_\alpha^0\},\{\rho^0_{\beta\alpha}\}):\Gamma^0\{V_\alpha\}
\rightarrow\Gamma\{U_\alpha\}$, whose equivalence class is the given
map $\Phi$. We further apply the trick described above to replace
$(\{f_\alpha^0\},\{\rho^0_{\beta\alpha}\})$ by a homomorphism
$(\{f_\alpha\},\{\rho_{\beta\alpha}\}):\Gamma\{V_\alpha\}
\rightarrow\Gamma\{U_\alpha\}$ such that each $\rho_\alpha$ is an
isomorphism and $(\{f_\alpha^0\},\{\rho^0_{\beta\alpha}\})$ is the
restriction of $(\{f_\alpha\},\{\rho_{\beta\alpha}\})$ to the sub-groupoid
$\Gamma^0\{V_\alpha\}$ of $\Gamma\{V_\alpha\}$.

Note that for each index $\alpha$ there may be more than one $U$ such that
$V_\alpha$ is a connected component of $\phi^{-1}(U)$. We shall modify
$(\{f_\alpha\},\{\rho_{\beta\alpha}\})$ by adding all such $U$ to
$\{U_\alpha\}$ and by allowing the $V_\alpha$'s in $\{V_\alpha\}$ to
repeat. The resulting homomorphism may be expressed as
$$
(\{V_{a,i}|i\in I_a\},\{U_a\},\{f_{a,i}\},\{\rho_{(b,j)(a,i)}\})
\leqno (3.1.5)
$$
where $\{U_a\}$ is the set of local charts on $X$ such that
$\phi^{-1}(U_a)\neq\emptyset$, and $\{V_{a,i}|i\in I_a\}$ is the set of
connected components of $\phi^{-1}(U_a)$. Note that $V_{a,i},V_{b,j}$ may
be identical even if $a\neq b$. Finally, we wish to emphasize that each
$\rho_{(a,i)}:G_{V_{a,i}}\rightarrow G_{U_a}$ is an isomorphism.

With these preparations, we shall construct a canonical orbispace
structure on the mapping cylinder $M_\phi$ as follows.

First of all, we specify the atlas of local charts on $M_\phi$. Observe
that the atlas of local charts on $X$ is the disjoint union of $\{U_a\}$
with a set $\{U_i\}$ where each $\phi^{-1}(U_i)=\emptyset$. If we set
$\phi_a=\phi|_{\phi^{-1}(U_a)}$, then it is readily seen that $\W=
\{M_{\phi_a}\}\cup\{U_i\}$ is a cover of the mapping cylinder $M_\phi$.
To each $W\in\W$, we assign a $(\widehat{W},G_W,\pi_W)$ such that
$\pi_W:\widehat{W}\rightarrow M_{\phi}$ induces a homeomorphism
$\widehat{W}/G_W\cong W$ as follows. If $W=U_i$ for some index $i$,
we simply put $(\widehat{W},G_W,\pi_W)=(\widehat{U_i},G_{U_i},\pi_{U_i})$.
If $W=M_{\phi_a}$ for some index $a$, we define $\widehat{W}=M_{f_a}$,
where $f_a=\sqcup_{i\in I_a} f_{a,i}:\bigsqcup_{i\in I_a}\widehat{V_{a,i}}
\rightarrow\widehat{U_a}$ (cf. $(3.1.5)$), and define $G_W=G_{U_a}$
with the action on $\widehat{W}$ given by the extension of the one on
$\widehat{U_a}$ by $g\cdot [y,t]=[\rho_{(a,i)}^{-1}(g)\cdot y,t]$,
$\forall g\in G_{U_a}, [y,t]\in M_{f_{a,i}}\subset\widehat{W}$, and
define $\pi_W$ by $[y,t]\mapsto [\pi_{V_{a,i}}(y),t], \forall [y,t]\in
M_{f_{a,i}}$ and $x\mapsto\pi_{U_a}(x), \forall x\in\widehat{U_a}$.

Secondly, we define the set $\T=\{T(W_1,W_2)\mid W_1,W_2\in\W,
\mbox{ s.t. } W_1\cap W_2\neq\emptyset\}$. We first look at the most
complicated case where $W_1=M_{\phi_a}$, $W_2=M_{\phi_b}$ for some
indexes $a,b$. In this case it is important to observe that the set
of connected components of $U_a\cap U_b$ is a subset $\{U_c\mid c\in
I_{a,b}\}\subset \{U_a\}$, and correspondingly the set of connected
components of $W_1\cap W_2$ is $\{W_c=M_{\phi_c}\mid c\in I_{a,b}\}$.
As a set, we shall define $T_{W_c}(W_1,W_2)=T_{U_c}(U_a,U_b)$, and
define $T(W_1,W_2)=\bigsqcup_{c\in I_{a,b}}T_{W_c}(W_1,W_2)=
\bigsqcup_{c\in I_{a,b}}T_{U_c}(U_a,U_b)=T(U_a,U_b)$. However, each
$\xi\in T_{W_c}(W_1,W_2)$ is assigned with a pair $(\bar{\phi}_\xi,
\bar{\lambda}_\xi)$ as follows. First, we regard $\xi\in T_{U_c}
(U_a,U_b)$ and write $\xi=\xi_2\circ\xi_1^{-1}$ for some $\xi_1\in
T(U_c,U_a)$, $\xi_2\in T(U_c,U_b)$. Second, set $\{\eta_\alpha\mid
\alpha\in I_{\xi_1}\}=\{\eta\in T(V_{c,s},V_{a,i})\mid s\in I_c,
i\in I_a, \mbox{ s.t. }\rho_{(a,i)(c,s)}(\eta)=\xi_1\}$. It is a routine
exercise to check that (1) both $\eta_\alpha\mapsto\mbox{Domain }
(\phi_{\eta_\alpha})$, $\eta_\alpha\mapsto\mbox{Range }(\phi_{\eta_\alpha})$
are bijections, and (2) $\bigsqcup_{\alpha\in I_{\xi_1}}\mbox{Domain }
(\phi_{\eta_\alpha})=\bigsqcup_{s\in I_c}\widehat{V_{c,s}}$ and
$\bigsqcup_{\alpha\in I_{\xi_1}}\mbox{Range }(\phi_{\eta_\alpha})
=f_a^{-1}(\mbox{Domain }(\phi_{\xi_1}))\subset\bigsqcup_{i\in I_a}
\widehat{V_{a,i}}$. It then follows that $\sqcup_{\alpha\in I_{\xi_1}}
\phi_{\eta_\alpha}:\bigsqcup_{s\in I_c}\widehat{V_{c,s}}\rightarrow f_a^{-1}
(\mbox{Domain }(\phi_{\xi_1}))$ is a homeomorphism which satisfies $f_a\circ
(\sqcup_{\alpha\in I_\xi}\phi_{\eta_\alpha})=\phi_{\xi_1}\circ f_c$.
We define $\bar{\phi}_{\xi_1}:M_{f_c}\rightarrow M_{f_a}$ to be the
corresponding open embedding between the mapping cylinders. Clearly
the range of $\bar{\phi}_{\xi_1}$ is a connected component of
$\pi_{W_1}^{-1}(W_c)$, and $\bar{\phi}_{\xi_1}$ is equivariant with
respect to $\lambda_{\xi_1}:G_{U_c}=G_{W_c}\rightarrow G_{U_a}=G_{W_1}$.
Similarly, one has an open embedding $\bar{\phi}_{\xi_2}:M_{f_c}\rightarrow
M_{f_b}$ which is $\lambda_{\xi_2}$-equivariant. Finally, we define
$(\bar{\phi}_\xi,\bar{\lambda}_\xi)=(\bar{\phi}_{\xi_2},\lambda_{\xi_2})
\circ (\bar{\phi}_{\xi_1}^{-1},\lambda_{\xi_1}^{-1})$, which is independent
of the choices on $\xi_1,\xi_2$. The definition for the remaining cases is
obvious, so we leave the details to the reader.

Finally, the composition and inverse for the elements in $\T$ are to
be inherited directly from those in the orbispace structure on $X$
under the natural identification described above. By Proposition 2.1.1
in \cite{C1}, $(\W,\T)$ defines an orbispace structure on $M_\phi$.

It remains to define the maps of orbispaces $i:Y\rightarrow M_\phi$,
$j:X\rightarrow M_\phi$ and $r:M_\phi\rightarrow X$ with the claimed
properties. First, $i:Y\rightarrow M_\phi$ is defined by the homomorphism
$$
(\{V_{a,i}|i\in I_a\},\{M_{\phi_a}\},\{i_{a,i}\},\{\rho_{(b,j)(a,i)}\})
\leqno (3.1.6)
$$
where $i_{a,i}:\widehat{V_{a,i}}\rightarrow M_{f_a}$ is $y\mapsto [y,0]$.
It realizes $Y$ as a subspace of $M_\phi$ essentially because by
Lemma 3.1.2, the orbispace structure given by $\V(\Phi),\T(\Phi)$
is equivalent to the original one on $Y$ in the sense of
Remark 2.1.2 (5) in \cite{C1}. Second, $j:X\rightarrow M_\phi$ is
defined by $(\{j_\alpha\},\{\delta_{\beta\alpha}\}):\Gamma\{U_\alpha\}
\rightarrow\{W_\alpha\}$, which is clearly a subspace, where (a)
$W_\alpha=U_i$ if $U_\alpha=U_i$, $W_\alpha=M_{\phi_a}$ if $U_\alpha=U_a$,
(b) $j_\alpha=Id$ if $U_\alpha=U_i$ and $j_\alpha=j_a:\widehat{U_a}
\rightarrow M_{f_a}$ if $U_\alpha=U_a$, and (c) each $\delta_{\beta\alpha}$
is the identity map under the natural identification.

Last, the retraction $r:M_\phi\rightarrow X$ is defined by
$(\{r_\alpha\},\{\delta_{\beta\alpha}\}):\Gamma\{W_\alpha\}
\rightarrow\Gamma\{U_\alpha\}$, where (a) $U_\alpha=U_i$ if $W_\alpha=U_i$,
$U_\alpha=U_a$ if $W_\alpha=M_{\phi_a}$, (b) $r_\alpha=Id$ if $W_\alpha=U_i$
and $r_\alpha=r_a:M_{f_a}\rightarrow\widehat{U_a}$ is the usual retraction
if $W_\alpha=M_{\phi_a}$, and (c) each $\delta_{\beta\alpha}$ is the
identity map. Clearly $r\circ j=Id_X$. On the other hand, there is a
canonical homotopy between $j\circ r$ and $Id_{M_\phi}$, defined by
$(\{H_\alpha\},\{\delta_{\beta\alpha}\}):\Gamma\{W_\alpha\times [0,1]\}
\rightarrow\Gamma\{U_\alpha\}$, where $H_\alpha=Id$ if $W_\alpha=U_i$,
and $H_\alpha$ is the usual homotopy between $j_a\circ r_a$ and
$Id_{M_{f_a}}$ if $W_\alpha=M_{\phi_a}$. Hence $j:X\rightarrow M_\phi$
is a strong deformation retract by $r:M_\phi\rightarrow X$. Finally,
we note that $\Phi=r\circ i$.

\hfill $\Box$

\subsection{Orbispaces via attaching cells of isotropy}

In this subsection we apply the mapping cylinder construction in
the preceding subsection to a special case where $Y=S^{k-1}(G)$,
$k\geq 1$, and make the meaning of `attaching a k-cell of isotropy
type $G$ to an orbispace' mathematically precise. To this end, we
have to impose further conditions on both the orbispace $X$ and the
attaching map $\Phi:S^{k-1}(G)\rightarrow X$, in order to deal with
two additional issues that are involved in the process.

The first one is how to construct an orbispace by collapsing a
subspace to a point. Suppose $A\subset X$ is a closed, connected
subspace of an orbispace $X$. We denote by $X/A$ the topological space
obtained by collapsing $A$ in the underlying space of $X$ to a point and
by $\ast\in X/A$ the image of $A$ under the canonical projection
$X\rightarrow X/A$.

\begin{lem}
Suppose $A$ is further contained in a local chart on $X$. Then
there is a canonical orbispace structure on $X/A$, and a map of
orbispaces $\pi:X\rightarrow X/A$ covering the canonical projection
$X\rightarrow X/A$ between the underlying spaces, such that the
restriction of $\pi$ to the open subspace $X\setminus A$ is an
isomorphism of orbispaces onto $(X/A)\setminus\{\ast\}$.
\end{lem}

\pf
Let $\U_1$ be the set of local charts $U$ on $X$ such that $U\cap
A=\emptyset$, and $\U_2$ be the set of local charts $U$ on $X$ such
that $A\subset U$. By the assumption, $\U_2\neq\emptyset$.

We shall define a canonical orbispace structure on $X/A$, where the atlas
of local charts $\V$ consists of connected open subsets $V$ such
that either $V\in\U_1$ as a subset of $X$, or $V=U/A$ for some
$U\in\U_2$. In the former case, we let $(\widehat{V},G_V,\pi_V)$ be
the one in the orbispace structure on $X$, while in the latter
case, we let $\widehat{V}$ be the space obtained by collapsing each
connected component of $\pi_U^{-1}(A)$ in $\widehat{U}$ to a point,
let $G_V=G_U$ with the induced action on $\widehat{V}$, and let
$\pi_V:\widehat{V}\rightarrow V$ be the map induced by $\pi_U$. The
set $\T=\{T(V_1,V_2)\mid V_1,V_2\in\V, \mbox{ s.t. } V_1\cap V_2
\neq\emptyset\}$ is identical to the one in the orbispace
structure on $X$ along with the assignment $\xi\mapsto\phi_\xi$ and
the composition and inverse, except for the case when $\xi\in
T_V(V_1,V_2)=T_U(U_1,U_2)$, where $V_i=U_i/A$, $i=1,2$, and
$V=U/A$, we instead assign $\xi$ with $\bar{\phi}_\xi$, whose
domain is the space obtained by collapsing each connected component
of $\pi_{U_1}^{-1}(A)$ in $\mbox{Domain }(\phi_\xi)$ to a point,
whose range is the space obtained by collapsing each connected
component of $\pi_{U_2}^{-1}(A)$ in $\mbox{Range }(\phi_\xi)$ to a
point, and $\bar{\phi}_\xi$ is the map induced by $\phi_\xi$. It is clear
that $\V,\T$ define an orbispace structure on $X/A$.

The map $\pi:X\rightarrow X/A$ is defined by $(\{\pi_\alpha\},
\{\delta_{\beta\alpha}\}):\Gamma\{U_\alpha\}\rightarrow
\Gamma\{V_\alpha\}$, where (a) $\{U_\alpha\}=\U_1\cup\U_2$, which is a
cover of $X$ since $A$ is closed, and $\{V_\alpha\}=\V$, (b)
$\pi_\alpha=Id$ if $U_\alpha\in\U_1$, and $\pi_\alpha$ is the
canonical projection if $U_\alpha\in\U_2$, and (c) each
$\delta_{\beta\alpha}$ is the identity map, which is clearly an
isomorphism onto $(X/A)\setminus\{\ast\}$ when restricted to the
open subspace $X\setminus A$.

\hfill $\Box$

In order to address the second preparatory issue, we note that
there is a natural map of orbispaces $k:Y\times I\rightarrow
M_\phi$, which is defined by the homomorphism
$$
(\{V_{a,i}\times I|i\in I_a\},\{M_{\phi_a}\},\{k_{a,i}\},
\{\rho_{(b,j)(a,i)}\}), \leqno (3.2.1)
$$
where $k_{a,i}:\widehat{V_{a,i}}\times I\rightarrow M_{f_a}$ is
$(y,t)\mapsto [y,t]$, cf. $(3.1.5), (3.1.6)$.

\begin{lem}
Let $\Phi:Y\rightarrow X$ be any map of orbispaces whose mapping
cylinder $M_\phi$ is defined. Suppose both $Y,X$ are compact and
Hausdorff, and for any $V\in\V(\Phi)$, $\widetilde{V}$ is locally
compact, Hausdorff, and the map $\pi^V:\widetilde{V}\rightarrow V$
is proper. Then for any maps of orbispaces $\Psi:X\rightarrow
X^\prime$, $\Upsilon:Y\times I\rightarrow X^\prime$ such that
$\Psi\circ\Phi=\Upsilon|_{Y\times\{1\}}$, there is a map $\Xi:
M_\phi\rightarrow X^\prime$ which satisfies $\Xi\circ j=\Psi$ and
$\Xi\circ k=\Upsilon$.
\end{lem}

\pf
We represent $\Psi:X\rightarrow X^\prime$ by a homomorphism
$\tau=(\{\psi_\alpha\},\{\lambda_{\beta\alpha}\}):\Gamma\{U_\alpha\}
\rightarrow\Gamma\{U_{\alpha^\prime}^\prime\}$, where $\{U_\alpha\}$
is a finite cover of $X$ such that each $U_\alpha$ is admissible as
defined in \S 3.2 of \cite{C1}, and each $\psi_\alpha$ can be extended
over the closure of $\widehat{U_\alpha}$. This is possible because
$X$ is compact and Hausdorff. For each $U_\alpha$, let
$\{V_{\alpha,i}\mid i\in I_\alpha\}$ be the set of connected components
of $\phi^{-1}(U_\alpha)$. We assign $U_\alpha$ to each $V_{\alpha,i}$,
then by Lemma 3.1.2, there is a homomorphism $\sigma=(\{f_{\alpha,i}\},
\{\rho_{(\beta,j)(\alpha,i)}\}):\Gamma\{V_{\alpha,i}\}\rightarrow
\Gamma\{U_\alpha\}$, whose equivalence class is the map $\Phi:
Y\rightarrow X$. The composition $\epsilon=(\{\varphi_{\alpha,i}\},
\{\delta_{(\beta,j)(\alpha,i)}\})$, where $\varphi_{\alpha,i}=\psi_\alpha
\circ f_{\alpha,i}$ and $\delta_{(\beta,j)(\alpha,i)}=\lambda_{\beta\alpha}
\circ\rho_{(\beta,j)(\alpha,i)}$, is a representative of $\Psi\circ\Phi$,
which is admissible as defined in \S 3.2 of \cite{C1} by the assumptions
made on $Y$, $\V(\Phi)$ and $\tau$.

On the other hand, we represent $\Upsilon:Y\times I\rightarrow X^\prime$
by a homomorphism $\kappa=(\{h_a\},\{\eta_{ba}\}):\Gamma\{V_a\times I_a\}
\rightarrow\Gamma\{U_{a^\prime}^\prime\}$, where $a\in\Lambda$, which
may be made admissible as defined in \S 3.2 of \cite{C1} by the
assumptions we have on $Y$ and $\V(\Phi)$. We may require $\#\Lambda<\infty$
since $Y$ is compact. Set $\Lambda_0=\{a\in\Lambda|1\in
I_a\}$, and we assume without loss of generality that $\overline{I_a}
\subset [0,1)$, $\forall a\in\Lambda\setminus\Lambda_0$. Note that
$\kappa(1)=(\{h_a(\cdot,1)\},\{\eta_{ba}\}):\Gamma\{V_a\}\rightarrow
\Gamma\{U_{a^\prime}^\prime\}$, where $a\in\Lambda_0$, also represents
$\Upsilon|_{Y\times \{1\}}=\Psi\circ\Phi$. By passing to an induced
homomorphism of $\kappa$, we may assume that $\kappa(1)$ is induced by
$\epsilon$ via some $\bar{\gamma}=(\theta,\{\xi_a\},\{\xi_a^\prime\})$.
By Lemma 3.2.2 in \cite{C1}, there is a local homeomorphism 
$\phi_{\bar{\gamma}}$
from an open neighborhood of $\kappa(1)$ onto an open neighborhood 
of $\epsilon$,
sending $\kappa(1)$ to $\epsilon$. Since $\Lambda$ is finite and 
$\overline{I_a}
\subset [0,1)$, $\forall a\in\Lambda\setminus\Lambda_0$, there is a
$t_0\in [0,1)$, $t_0\in I\setminus\overline{I_a},\forall a\in\Lambda\setminus
\Lambda_0$, such that for any $t\in (t_0,1]$, $\kappa(t)=(\{h_a(\cdot,t)\},
\{\eta_{ba}\}):\Gamma\{V_a\}\rightarrow\Gamma\{U_{a^\prime}^\prime\}$,
where $a\in\Lambda_0$, lies in the domain of $\phi_{\bar{\gamma}}$.
We write $\phi_{\bar{\gamma}}(\kappa(t))=\epsilon(t)=
(\{\varphi_{\alpha,i}^{(t)}\},\{\delta_{(\beta,j)(\alpha,i)}\})$,
$t\in (t_0,1]$. Note that $\varphi_{\alpha,i}^{(1)}=\psi_\alpha
\circ f_{\alpha,i}$.

Introduce the following notation: Let $f$ be a continuous map.
For any $t_0\in [0,1)$, we denote by $M_f(t_0)$ the open subset
$M_f\setminus \{[y,t]\mid t\in [0,t_0]\}$ in the mapping cylinder
$M_f$ of $f$. Then there is a homomorphism $\varsigma=(\{\chi_\alpha\},
\{\lambda_{\beta\alpha}\}):\Gamma\{W_\alpha\}\rightarrow\Gamma\{
U_{\alpha^\prime}^\prime\}$, where either $W_\alpha=U_\alpha$ and
$\chi_\alpha=\psi_\alpha$, or $W_\alpha=M_{\phi_\alpha}(t_0)$ and
$\chi_\alpha=\psi_\alpha$ on $j_\alpha(\widehat{U_\alpha})\subset
M_{f_\alpha}(t_0)=\widehat{W_\alpha}$, $\chi_\alpha([y,t])=
\varphi_{\alpha,i}^{(t)}(y)$, $\forall (y,t)\in\widehat{V_{\alpha,i}}
\times (t_0,1], i\in I_\alpha$. On the other hand, consider the
restriction of $\kappa$ to $Y\times (I\setminus \{1\})$, which is
still denoted by $\kappa=(\{h_a\},\{\eta_{ba}\}):\Gamma\{V_a\times I_a\}
\rightarrow\Gamma\{U_{a^\prime}^\prime\}$ for the sake of simplicity.
We can join $\varsigma$ and $\kappa$ via $\bar{\gamma}=(\theta,
\{\xi_a\},\{\xi_a^\prime\})$ to construct a homomorphism $\varsigma
\cup_{\bar{\gamma}}\kappa$ as follows. We add a set of mappings
$\lambda_{\alpha a}:T(V_a\times I_a,W_\alpha)\rightarrow T(U_a^\prime,
U_\alpha^\prime)$, where $(V_a\times I_a)\cap W_\alpha\neq\emptyset$,
to $\varsigma,\kappa$. Note that any $\xi\in T(V_a\times I_a,W_\alpha)$
may be regarded as an element of $T(V_a,V_{\alpha,i})$ for some
$i\in I_\alpha$. In this case we define (cf. $(3.2.8)$ in \S 3.2 of
\cite{C1})
$$
\lambda_{\alpha a}(\xi)=\delta_{(\alpha,i)\theta(a)}(\xi\circ\xi_a^{-1})
\circ\xi_a^\prime(x), \;\forall x\in h_a(\mbox{Domain }(\phi_\xi)).
\leqno (3.2.2)
$$
As seen in the proof of Lemma 3.2.2 in \cite{C1}, $\varsigma
\cup_{\bar{\gamma}}\kappa$ is indeed a homomorphism, whose equivalence
class is a map from the orbispace $M_\phi$ to $X^\prime$. We define
$\Xi=[\varsigma\cup_{\bar{\gamma}}\kappa]$, which clearly satisfies
$\Xi\circ j=\Psi$ and $\Xi\circ k=\Upsilon$.

\hfill $\Box$

We remark that the map $\Xi$ in the preceding lemma may not be
uniquely determined, but the ambiguity is caused only by the
different choices of $\gamma=[\bar{\gamma}]\in
\Gamma_{\epsilon\kappa(1)}$ (cf. Lemma 3.1.1 in \cite{C1}). If
we work with the based version, then $\#\Gamma_{\epsilon{\kappa(1)}}=1$
when $Y$ is connected, and $\Xi$ is uniquely determined in this case.

Now we are ready to describe the precise meaning of attaching a
k-cell $D^k(G)$ of isotropy type $G$ to an orbispace $X$ via a map
$\Phi:S^{k-1}(G)\rightarrow X$. First of all, we list the
additional conditions we need to impose on the orbispace $X$ and
the attaching map $\Phi:S^{k-1}(G)\rightarrow X$, where $k\geq 1$:

\begin{itemize}
\item $X$ is compact and Hausdorff.
\item The orbispace structure on $S^{k-1}$ given by $\V(\Phi),\T(\Phi)$
is contained in the standard one on $S^{k-1}(G)$.
\end{itemize}

Note that the last condition means that for each $V\in\V(\Phi)$,
$\widetilde{V}=V$, $G^V=G$ which acts on $\widetilde{V}$ trivially,
and $\pi^V:\widetilde{V}\rightarrow V$ is the identity map.
Moreover, each $T(V_1,V_2)\in\T(\Phi)$ is a disjoint union of a
number of copies of $G$ which is naturally labeled by the set of
connected components of $V_1\cap V_2$, and for each element $\xi\in
T(V_1,V_2)$, $\phi_\xi$ is the identity map, and
$\lambda_\xi=Ad(\xi):G\rightarrow G$.

Recall that the mapping cone of a continuous map $\phi:Y\rightarrow
X$, denoted by $C_\phi$, is the space $M_\phi/i(Y)$ obtained by
collapsing the subspace $i:Y\rightarrow M_\phi$ to a point. When
$Y=S^{k-1}$, $C_\phi$ is the result of ``attaching a k-cell $D^k$ to
$X$ via the map $\phi:S^{k-1}=\partial{D^k}\rightarrow X$''. The
interior of the k-cell in $C_\phi$ is the image of the open subset
$D_0=\{[y,t]\mid y\in S^{k-1},t\in [0,1)\}\subset M_\phi$ under the
canonical projection $M_\phi\rightarrow C_\phi=M_\phi/i(S^{k-1})$.

\begin{prop}
With the preceding understood, we assert that
\begin{itemize}
\item [{(1)}] There is a canonical orbispace structure on $C_\phi$
such that {\em (a)} there is a canonical map of orbispaces $\pi:
M_\phi\rightarrow C_\phi$ which is an open embedding when
restricted to $M_\phi\setminus i(S^{k-1}(G))$, and {\em (b)} the
interior of the attached k-cell in $C_\phi$ is an open k-cell of
isotropy type $G$ as a subspace of the orbispace $C_\phi$.
\item [{(2)}] The homotopy type of the orbispace $C_\phi$ depends
only on the homotopy class of the attaching map $\Phi$.
\item [{(3)}] The canonical embedding of the orbispace $X$ into
$C_\phi$ induced by $j:X\rightarrow M_\phi$ is a cofibration, i.e.,
it has the homotopy extension property with respect to all
orbispaces.
\item [{(4)}] A map $\Psi:X\rightarrow X^\prime$ can be extended
over the attached k-cell of isotropy type $G$ to $C_\phi$ iff
$\Psi\circ\Phi:S^{k-1}(G)\rightarrow X^\prime$ is null-homotopic,
and such an extension of $\Psi$ is given by a null-homotopy of
$\Psi\circ\Phi$. Moreover, suppose $F$ is a homotopy between
$\Psi_1$ and $\Psi_2$, which are extended to $C_\phi$ by the
null-homotopies $h_1,h_2$ of $\Psi_1\circ\Phi$, $\Psi_2\circ\Phi$
respectively. Then $F$ can be extended to $C_\phi\times [0,1]$ to a
homotopy between the corresponding extensions of $\Psi_1$ and
$\Psi_2$ iff the null-homotopies $h_1,h_2$ are homotopic via a
homotopy whose restriction to $S^{k-1}(G)\times [0,1]$ is $F\circ
(\Phi\times Id)$.
\end{itemize}
\end{prop}

\pf
(1) The canonical orbispace structure on $C_\phi$ will be the one
obtained from the canonical orbispace structure on the mapping cylinder
$M_\phi$ by collapsing the subspace $i(S^{k-1}(G))$ to a point. To this
end, we need to verify the hypothesis in Lemma 3.2.1. The case when
$k=1$ requires a separate, but similar argument because $i(S^{k-1}(G))$
is not connected. We shall only consider the cases when $k\geq 2$, the
remaining case is left to the reader for simplicity.

In order to apply Lemma 3.2.1, it suffices to show that we may add
the local chart $(\widehat{D_0},G_{D_0},\pi_{D_0})=(D_0,G,\pi_0)$,
where $D_0=\{[y,t]\mid y\in S^{k-1},t\in [0,1)\}\subset M_\phi$,
$G$ acts on $D_0$ trivially, and $\pi_0:D_0\rightarrow D_0$ is the
identity map, to the canonical orbispace structure on $M_\phi$ constructed
in Proposition 1.5, so that with respect to the new orbispace structure
on $M_\phi$, which contains the original one hence equivalent, the
subspace $i(S^{k-1}(G))$ is contained in a local chart, i.e.,
$(\widehat{D_0},G_{D_0},\pi_{D_0})=(D_0,G,\pi_0)$. The condition
that the orbispace structure on $S^{k-1}$ given by $\V(\Phi),\T(\Phi)$
is contained in the standard one on $S^{k-1}(G)$ guarantees this.
We shall continue to use the notations in the proof of Proposition 1.5
concerning the canonical orbispace structure on $M_\phi$.

More concretely, we need to define a set $\{T(D_0,W)\mid W\in\W,
\mbox{ s.t. } D_0\cap W\neq\emptyset\}$, and add it to the canonical
orbispace structure on $M_\phi$ along with $(\widehat{D_0},G_{D_0},
\pi_{D_0})$. Note that if $D_0\cap W\neq\emptyset$, then $W$ must be
$W_a=M_{\phi_a}$ for some index $a$, where $M_{\phi_a}$ is the mapping
cylinder of $\phi_a=\phi|_{\phi^{-1}(U_a)}:\bigsqcup_{i\in I_a}V_{a,i}
\rightarrow U_a$, and the set of connected components of $D_0\cap W_a$
is $\{W_{a,i}=V_{a,i}\times [0,1)|i\in I_a\}$.  We define $T(D_0,W_a)
=\bigsqcup_{i\in I_a}T_{W_{a,i}}(D_0,W_a)$ where $T_{W_{a,i}}(D_0,W_a)
=G_{U_a}$, and assign each $\xi\in T_{W_{a,i}}(D_0,W_a)$ with a pair
$(\phi_\xi,\lambda_\xi)$ as follows. Note that the assumption that the
orbispace structure on $S^{k-1}$ given by $\V(\Phi),\T(\Phi)$ is contained
in the standard one on $S^{k-1}(G)$ implies that the inverse image of $W_{a,i}$
in $\widehat{W_a}$ is $(G_{U_a}/\rho_{(a,i)}(G))\times W_{a,i}$, where
$G_{U_a}/\rho_{(a,i)}(G)$ is the set of right cosets. The map $\phi_\xi$
is the homeomorphism sending $x\in W_{a,i}$ to $([\xi],x)\in
(G_{U_a}/\rho_{(a,i)}(G))\times W_{a,i}$, and $\lambda_\xi=Ad(\xi)\circ
\rho_{(a,i)}:G\rightarrow Ad(\xi)(\rho_{(a,i)}(G))$. It is easily seen
that this will give rise to an orbispace structure on $M_\phi$, containing
the original one. Now the subspace $i(S^{k-1}(G))$ is closed, connected,
and is contained in the local chart $(\widehat{D_0},G_{D_0},\pi_{D_0})$.
Hence by Lemma 3.2.1, we can collapse $i(S^{k-1}(G))$ to a point,
and a canonical orbispace structure is resulted on $C_\phi=M_\phi/
i(S^{k-1})$, with a canonical map of orbispaces $\pi:M_\phi\rightarrow
C_{\phi}$ whose restriction to $M_{\phi}\setminus i(S^{k-1})$ is an
isomorphism of orbispaces onto $C_\phi\setminus \{\ast\}$. The attached
k-cell in $C_\phi$ is the image of $\overline{D_0}=\{[y,t]\mid y\in S^{k-1},
t\in I\}\subset M_\phi$ under the projection $M_\phi\rightarrow C_\phi=
M_\phi/i(S^{k-1})$. Hence its interior is $D=D_0/i(S^{k-1})$. Clearly,
as a subspace of the orbispace $C_\phi$, it is $D(G)$, the open k-cell
of isotropy type $G$, cf. Lemma 3.2.1.

\vspace{1.5mm}

(2) This is obtained by applying Lemma 3.2.2 along with $(3.1.1)$,
$(3.1.2)$. Note that the hypothesis in Lemma 3.2.2 is met by the
assumptions made in the present proposition.

\vspace{1.5mm}

(3) By definition the embedding of orbispace $i:X\rightarrow C_\phi$ is
called a cofibration, where $i$ is $j:X\rightarrow M_\phi$ followed by
$\pi:M_\phi\rightarrow C_\phi$, if the following is true: given any
orbispace $X^\prime$ and any map $\Psi:C_\phi\rightarrow X^\prime$, if
$H:X\times [0,1]\rightarrow X^\prime$ is a homotopy from its
restriction $\Psi|_{X}$ to another map from $X$ to $X^\prime$, then
there is a homotopy $F:C_\phi\times [0,1]\rightarrow X^\prime$ extending
$H$, such that $F|_{C_\phi\times\{0\}}=\Psi$. Consequently,
$\Psi:C_\phi\rightarrow X^\prime$ is homotopic through $F$ to another map
from $C_\phi$ to $X^\prime$, i.e., $F|_{C_\phi\times\{1\}}$, extending
that for the restriction $\Psi|_{X}:X\rightarrow X^\prime$.

We apply Lemma 3.2.2 to the mapping cylinder of $\Phi\times Id:S^{k-1}(G)
\times [0,1]\rightarrow X\times [0,1]$, with the map $H:X\times [0,1]
\rightarrow
X^\prime$ and a map $\Upsilon:S^{k-1}(G)\times [0,1]\times I\rightarrow
X^\prime$ constructed as follows. Let $R:[0,1]\times I\rightarrow
[0,1]\times\{1\}\cup \{0\}\times I$ be a retraction which sends $[0,1]
\times\{0\}$ to $\{0\}\times \{0\}$. We define
$$
\Upsilon=(H\circ (\Phi\times Id)\cup \Psi\circ\pi)\circ (Id\times R),
$$
which satisfies $\Upsilon|_{S^{k-1}(G)\times [0,1]\times\{1\}}
=H\circ (\Phi\times Id)$ and $\Upsilon|_{S^{k-1}(G)\times\{0\}\times I}
=\Psi\circ\pi|_{S^{k-1}(G)\times I}$, where $\pi:M_\phi\rightarrow C_\phi$,
and $\Upsilon|_{S^{k-1}(G)\times [0,1]\times
\{0\}}=\Psi\circ\pi|_{i(S^{k-1}(G))}$. Let $\Xi:M_{\phi\times Id}
\rightarrow X^\prime$ be the resulting map. Then the last property
of $\Upsilon$ implies that $\Xi$ factors through $\pi\times Id:
M_{\phi\times Id}=M_\phi\times [0,1]\rightarrow C_\phi\times [0,1]$, which
results in the homotopy $F:C_\phi\times [0,1]\rightarrow X^\prime$
with the claimed properties.

\vspace{1.5mm}

(4) Straightforward application of Lemma 3.2.2.

\hfill $\Box$

We will say that the orbispace $C_\phi$ is obtained by attaching a
k-cell of isotropy type $G$ to $X$ via the map $\Phi$.

Now we consider a subcategory ${\cal{C}}$ of the category of
orbispaces introduced in Part I of this series \cite{C1}, which is
the union of ${\cal{C}}^n$, $n\geq 0$. Each orbispace $X\in {\cal{C}}^n$
admits a canonical filtration of subspaces
$$
X_0\subset X_1\subset\cdots\subset X_n=X, \leqno (3.2.3)
$$
where $X_0$ is a set of finitely many points in $X$, and for each $k=1,2,
\cdots,n$, $X_k$ is obtained by attaching to $X_{k-1}$ finitely many
k-cells of various isotropy type. We shall call each subspace $X_k$
in $(3.2.3)$, $0\leq k\leq n$, the k-skeleton of $X$.

Note that in general it is not clear that for any $X_1,X_2\in
{\cal{C}}$, the product $X_1\times X_2$ is still an object in ${\cal{C}}$.
However, let us be content for now with the observation that for
any $X\in {\cal{C}}^n$, the product $Z=X\times I\in {\cal{C}}^{n+1}$.
In fact, let
$$
X_0\subset X_1\subset\cdots\subset X_n=X
$$
be the filtration of skeletons for $X$. Then correspondingly
$$
Z_0\subset Z_1\subset\cdots\subset Z_{n+1}=Z,
$$
where $Z_0=X_0\times\{0\}\cup X_0\times\{1\}$, $Z_k=X_k\times\{0\}
\cup X_k\times\{1\}\cup X_{k-1}\times I$ for $1\leq k\leq n$, and
$Z_{n+1}=X_n\times I$, is the canonical filtration of skeletons
for $Z=X\times I$.

\vspace{2.5mm}

The remaining of this subsection is occupied by the study of some
fundamental homotopy properties of the objects in ${\cal{C}}$. We
begin by observing

\begin{lem}
For $k\geq 1$, an element $\Phi\in [(D^k(G),S^{k-1}(G),\ast);(X,A,
\underline{o})]_\rho$ defines the trivial element in $\pi_k^{(G,\rho)}
(X,A,\underline{o})=[[(D^k(G),S^{k-1}(G),\ast);(X,A,\underline{o})]]_\rho$
iff there is a homotopy $H$ between $\Phi$ and a $\Phi^\prime$ as
elements of $[(D^k(G),\ast);(X,\underline{o})]_\rho$, such that
$H|_{S^{k-1}(G)\times\{t\}}=\Phi|_{S^{k-1}(G)}$, $\forall t\in [0,1]$,
and $\Phi^\prime\in [(D^k(G),\ast);(A,\underline{o}|_A)]_\rho$.
\end{lem}

\pf
Suppose there is such a homotopy $H$. Then $H|_{S^{k-1}(G)\times\{t\}}
=\Phi|_{S^{k-1}(G)}$ and $\Phi\in [(D^k(G),S^{k-1}(G),\ast);(X,A,
\underline{o})]_\rho$ imply that $H$ is also a homotopy between
$\Phi$ and $\Phi^\prime$ as elements of $[(D^k(G),S^{k-1}(G),\ast);(X,A,
\underline{o})]_\rho$. On the other hand, $\Phi^\prime\in [(D^k(G),\ast);
(A,\underline{o}|_A)]_\rho$ implies that $\Phi^\prime$ defines the
trivial element in $\pi_k^{(G,\rho)}(X,A,\underline{o})$, hence so
does $\Phi$.

Conversely, suppose $\Phi\in [(D^k(G),S^{k-1}(G),\ast);(X,A,
\underline{o})]_\rho$ defines the trivial element in $\pi_k^{(G,\rho)}
(X,A,\underline{o})$. Then there exists an $F\in [(CD^k(G),CS^{k-1}(G),
\ast);(X,A,\underline{o})]_\rho$ such that $F|_{D^k(G)}=\Phi$.
Denote by $C$ the subset $\{(x,t)\mid ||x||\leq 1-t\}$ of $D^k\times
[0,1]$. Then the map $(x,t)\mapsto ((1-t)x,t)$ from $D^k\times
[0,1]$ onto $C$ factors through $CD^k$, which defines an isomorphism
$\Psi:CD^k(G)\rightarrow C(G)$. On the other hand, there is a map
$\Upsilon:D^k(G)\times [0,1]\rightarrow C(G)$ defined by $(x,t)\mapsto
(x,t), \forall (x,t)\in C$, $(x,t)\mapsto (x,1-||x||), \forall (x,t)\in
(D^k\times [0,1])\setminus C$, and $Id:G\rightarrow G$. We simply let
$H=F\circ\Psi^{-1}\circ\Upsilon$, which clearly provides the desired
homotopy between $\Phi$ and $\Phi^\prime=H|_{D^k(G)\times\{1\}}$,
i.e, $\Phi^\prime\in [(D^k(G),\ast);(A,\underline{o}|_A)]_\rho$, and $H$
satisfies $H|_{S^{k-1}(G)\times\{t\}}=\Phi|_{S^{k-1}(G)}$, $\forall t\in
[0,1]$.

\hfill $\Box$

Recall that $\pi_0^G(X)=[[B_G;X]]$ where $B_G$ denotes the 0-cell
of isotropy type $G$. We shall say that a pair $(Y,B)$, where $B$ is a
subspace of $Y$ via $i:B\rightarrow Y$, is 0-connected, if $i_\ast:
\pi_0^G(B)\rightarrow\pi_0^G(Y)$ is a bijection for all $G$. For $n\geq
1$, we say that $(Y,B)$ is n-connected if it is 0-connected and
for all possible $\underline{o}$ and $(G,\rho)$, $\pi_k^{(G,\rho)}
(Y,B,\underline{o})$ is trivial for $1\leq k\leq n$.

\begin{lem}
Suppose $(Y,B)$ is a n-connected pair and $X\in {\cal{C}}^m$
where $m\leq n$. Then any map $\Phi:X\rightarrow Y$ is homotopic to
a map $\Phi^\prime:X\rightarrow B$. Moreover, if a subspace $A\subset X$
is a union of cells and $\Phi|_A\in [A;B]$, we may even arrange to
have $\Phi|_A=\Phi^\prime|_A$.
\end{lem}

\pf
Let $X_0\subset X_1\subset\cdots\subset X_m=X$ be the canonical
filtration of $X$. Given any map $\Phi:X\rightarrow Y$, the
restriction to $X_0$, $\Phi|_{X_0}$, is a map from a finite
disjoint union of 0-cells of various isotropy type into $Y$.
Since $i_\ast:\pi_0^G(B)\rightarrow\pi_0^G(Y)$ is a bijection for
all $G$, $\Phi|_{X_0}$ is homotopic to a map $\Psi_0:X_0\rightarrow B$
as maps into $Y$. By Proposition 3.2.3 (3), $X_0\subset X$ is a
cofibration, hence this homotopy can be extended to $X$, so that $\Phi$
is homotopic to a map $\Phi_0^\prime:X\rightarrow Y$ such that
$(\Phi^\prime_0)|_{X_0}:X_0\rightarrow B$. Now consider the
restriction of $\Phi_0^\prime$ to any of the 1-cells of various isotropy
type in $X_1$ that are attached to $X_0$. There are $\underline{o}$,
$(G,\rho)$ such that it defines an element in $\pi_1^{(G,\rho)}(Y,B,
\underline{o})$, which is trivial because $(Y,B)$ is 1-connected.
By Lemma 3.2.4 and Proposition 3.2.3 (4), the restriction of $\Phi_0^\prime$
to $X_1$ is homotopic to a map $\Psi_1:X_1\rightarrow B$. Again because
$X_1\subset X$ is a cofibration, there is a homotopy between $\Phi_0^\prime$
and a $\Phi_1^\prime:X\rightarrow Y$ such that the restriction of
$\Phi_1^\prime$ to $X_1$ maps into $B$. The lemma follows by repeating
this process under the condition $m\leq n$. It is clear that if
the restriction of $\Phi$ to a cell in $X$ is a map into $B$, one
may keep it unchanged in the above process.

\hfill $\Box$

\begin{prop}
Let $X\in {\cal{C}}$. For any $n\geq 0$, the pair $(X,X_n)$ is n-connected.
\end{prop}

\pf
It suffices to show that for all $G$, $i_\ast:\pi_0^G(X_0)\rightarrow
\pi_0^G(X)$ is surjective, and for all possible $\underline{o},(G,\rho)$,
$\pi_k^{(G,\rho)}(X,X_n,\underline{o})$ is trivial, $1\leq k\leq n$.

The key to the proof is the following fact by Proposition 3.2.3
(1): Denote by $\ast$ the cone point in $C_\phi$, i.e., the image
of $i(S^{k-1})$ under $\pi:M_\phi\rightarrow C_\phi$. Since when
restricted to $M_\phi\setminus i(S^{k-1})$, $\pi$ is an isomorphism
onto $C_\phi\setminus\{\ast\}$, there is a strong deformation
retraction which shrinks the attached k-cell of isotropy type $G$
to its boundary after a point in the interior is removed.

Given any map $\Phi:B_G\rightarrow X$, if its image lies in the
interior of an attached k-cell in $X_k$, $k\geq 1$, then for
dimensional reason it is in the complement of another interior point.
By the said strong deformation retraction, $\Phi$ is homotopic to another
map $\Psi:B_G\rightarrow X$ whose image lies in $X_{k-1}$. By
induction, $\Phi$ is homotopic to a map in $X_0$, hence $i_\ast:\pi_0^G(X_0)
\rightarrow\pi_0^G(X)$ is surjective.

Similarly, let $\Phi\in [(D^k(G),S^{k-1}(G),\ast);(X,X_n,\underline{o})]_\rho$,
$1\leq k\leq n$, whose image meets the interior of an attached m-cell
$e^m=D^m(H)$, $m\geq n+1$. In this case, we need to employ Theorem 1.4 (2)
in \cite{C1} to approximate $\Phi$ by a smooth map first. More concretely,
let $D_1$, $D_2$ be the closed ball of radius $\frac{1}{3},\frac{2}{3}$
centered at the cone point, and let $W_1=\phi^{-1}(D_1),W_2=\phi^{-1}(D_2)$
where $\phi$ is the induced map of $\Phi$ between underlying
spaces. By Theorem 1.4 (2) in \cite{C1}, there is a smooth map
$\Upsilon:W_2\rightarrow e^m$ such that $||\Upsilon-\Phi|_{W_2}||_{C^0}
<\frac{1}{100}$. Moreover, by Theorem 1.4 (1), the difference $\Upsilon
-\Phi|_{W_2}$ can be regarded as a $C^0$ section of a $C^0$
orbifold vector bundle over $W_2$, where $\Phi|_{W_2}$ is
identified with the zero section. Let $\beta$ be a smooth function on $W_2$
compactly supported in the interior $W_2$, such that $|\beta|\leq 1$
and $\beta=1$ on $W_1$. Then the map $\Phi^\prime:D^k(G)\rightarrow
X$, where $\Phi^\prime=\Phi+\beta(\Upsilon-\Phi|_{W_2})$, is
homotopic to $\Phi$ relative to $S^{k-1}(G)$, and the image of $\Phi^\prime$
will miss a point $z$ in the open ball of radius $\frac{1}{4}$
centered at the cone point. The reason for the latter: on the interior
$W_1$, $\Phi^\prime$ is smooth so that it will miss $z$ by transversality
(represent $\Phi^\prime|_{W_1}$ by a $(\{f_i\},\{\rho_{ji}\})$ and then
apply transversality argument to each $f_i$); on the complement,
$||\Phi^\prime-\Phi||_{C^0}<\frac{1}{100}$ so that it will be in the
complement of the ball of radius $\frac{1}{4}$ centered at the
cone point. By the strong deformation retraction mentioned
earlier, $\Phi^\prime$ is homotopic relative to $S^{k-1}(G)$ to a
map whose image will miss the entire interior of the m-cell $e^m$. By
induction, $\Phi$ is homotopic relative to $S^{k-1}(G)$ to a map
into $X_n$. By Lemma 3.2.4, $\pi_k^{(G,\rho)}(X,X_n,\underline{o})$
is trivial for $1\leq k\leq n$.

\hfill $\Box$

\begin{coro}
Let $X,X^\prime\in {\cal{C}}$. Then any map $\Phi:X\rightarrow X^\prime$
is homotopic to a map $\Psi:X\rightarrow X^\prime$, which is
`cellular' in the sense that $\Psi|_{X_n}\in [X_n;X_n^\prime]$ for
any $n\geq 0$. Moreover, if a subspace $A\subset X$ is a union of
cells and $\Phi|_{A\cap X_n}\in [A\cap X_n;X_n^\prime]$ for any $n\geq
0$, one may even arrange to have $\Phi|_A=\Psi|_A$.
\end{coro}

\noindent{\bf Proof of Theorem 1.6}

\vspace{1.5mm}

(1) Let $i:X\rightarrow M_\phi$, $j:X^\prime\rightarrow M_\phi$ be
the canonical embeddings into the mapping cylinder which realize $X,
X^\prime$ as a subspace, and let $r:M_\phi\rightarrow X^\prime$ be
the canonical strong deformation retraction, which satisfies $r\circ
i=\Phi$. By the exact sequence $(2.3.2)$ and the assumption that $\Phi_\ast$
is a weak homotopy equivalence, the pair $(M_\phi,X)$ is n-connected
for any $n\geq 0$.

Now let $Y\in {\cal{C}}$ be any element. We first show that $\Phi_\ast:
[[Y;X]]\rightarrow [[Y;X^\prime]]$ is surjective. In other words,
for any map $\Psi^\prime:Y\rightarrow X^\prime$, we will find a
map $\Psi:Y\rightarrow X$ such that $\Phi\circ\Psi\cong \Psi^\prime$.
This is done by applying Lemma 3.2.5 to the map $j\circ\Psi^\prime:
Y\rightarrow M_\phi$, which gives a map $\Psi:Y\rightarrow X$
satisfying $i\circ\Psi\cong j\circ\Psi^\prime$. Composing both
sides with $r:M_\phi\rightarrow X^\prime$, we obtain $\Phi\circ\Psi
=r\circ i\circ\Psi\cong r\circ j\circ\Psi^\prime=\Psi^\prime$.

As for the injectivity of $\Phi_\ast$, we apply the above argument
with $Y$ replaced by $Z=Y\times I$, which is also in ${\cal{C}}$.
More precisely, suppose $\Psi_1,\Psi_2:Y\rightarrow X$ are any two
maps such that $\Phi\circ\Psi_1\cong\Phi\circ\Psi_2$ via a homotopy
$H:Y\times I\rightarrow X^\prime$. Since $j\circ r\cong Id_{M_\phi}$,
we have $j\circ\Phi\circ\Psi_l=j\circ r\circ i\circ\Psi_l\cong i\circ\Psi_l$
for $l=1,2$. Combining with $j\circ H$, we obtain a homotopy $F^\prime:
Y\times I\rightarrow M_\phi$ between $i\circ\Psi_1$ and $i\circ\Psi_2$. On
the other hand, suppose $Y\in {\cal{C}}^n$, then $Z=Y\times
I\in {\cal{C}}^{n+1}$, with the canonical filtration of skeletons
$Z_0\subset Z_1\subset\cdots\subset Z_{n+1}=Z$, where
$Z_0=Y_0\times\{0\}\cup Y_0\times\{1\}$,
$Z_k=Y_k\times\{0\}\cup Y_k\times\{1\}\cup Y_{k-1}\times I$ for
$1\leq k\leq n$, and $Z_{n+1}=Y_n\times I$. Clearly, the subspace
$Y\times\{0\}\cup Y\times\{1\}\subset Z$ is a union of cells. Furthermore,
$F^\prime|_{Y\times\{0\}}=i\circ\Psi_1
\in [Y;X]$ and $F^\prime|_{Y\times\{1\}}=i\circ\Psi_2\in [Y;X]$. Hence
by Lemma 3.2.5, $F^\prime$ is homotopic to an $F:Y\times I\rightarrow X$
satisfying $F^\prime|_{Y\times\{0\}\cup Y\times\{1\}}=F|_{Y\times\{0\}
\cup Y\times\{1\}}$. It is easily seen that $\Psi_1\cong\Psi_2$ via $F$.
Thus $\Phi_\ast:[[Y;X]]\rightarrow [[Y;X^\prime]]$ is injective,
and hence a bijection.

\vspace{1.5mm}

(2) Suppose $\Phi:X\rightarrow X^\prime$ is a homotopy equivalence.
We need to show that (1) for all $G$, $\Phi_\ast:\pi_0^G(X)\rightarrow
\pi_0^G(X^\prime)$ is a bijection, which is trivial because $\pi_0^G(X)
=[[B_G;X]]$, and (2) for all possible data $\underline{o},
\underline{o^\prime},(G,\rho)$ and
$(G,\rho^\prime)$, $\Phi_\ast:\pi_k^{(G,\rho)}(X,\underline{o})\rightarrow
\pi_k^{(G,\rho^\prime)}(X^\prime,\underline{o^\prime})$
is isomorphic for all $k\geq 0$. Here possible $\underline{o},
\underline{o^\prime}$ are meant to be those with respect to which $\Phi$
has a based version
$\Phi\in [(X,\underline{o});(X^\prime,\underline{o^\prime})]_\eta$
for some injective homomorphism $\eta:G_{\hat{o}}\rightarrow
G_{\hat{o}^\prime}$. For instance, this is the case when $\Phi$ can
be represented by a homomorphism which also defines an element in
$[(X,\underline{o});(X^\prime,\underline{o^\prime})]_\eta$.
Note that not all base-point structures are being considered here. But in
light of Proposition 1.3 (3), no generality is lost. Now let $\Phi_0,\Phi_1:
X\rightarrow X^\prime$ be any maps which have based versions $\Phi_0\in
[(X,\underline{o});(X^\prime,\underline{o^\prime}_0)]_{\eta_0}$ and
$\Phi_1\in [(X,\underline{o});(X^\prime,\underline{o^\prime}_1)]_{\eta_1}$
for some $\underline{o},\underline{o^\prime}_0$ and $\underline{o^\prime}_1$.
Furthermore, $\Phi_0,\Phi_1$ are homotopic through a homotopy $F$. Then
there is a guided path $u\in [(I(G_{\hat{o}}),0,1);(X^\prime,
\underline{o^\prime}_0,\underline{o^\prime}_1)]_{(\eta_0,\eta_1)}$, which
is defined by the restriction of $F$ to $\{o\}\times I$. We claim
that $(\Phi_0)_\ast=u_\ast\circ (\Phi_1)_\ast:\pi_k^{(G,\rho)}(X,\underline{o})
\rightarrow\pi_k^{(G,\eta_0\circ\rho)}(X^\prime,\underline{o^\prime}_0)$,
where $u_\ast$ is the isomorphism associated to the guided path
$u$, cf. Proposition 1.3 (1). It follows easily from the claim
that homotopy equivalence implies weak homotopy equivalence. To
prove the claim, we observe that $(u_t)_\ast\circ (\Phi_t)_\ast:
\pi_k^{(G,\rho)}(X,\underline{o})\rightarrow\pi_k^{(G,\eta_0\circ\rho)}
(X^\prime,\underline{o^\prime}_0)$ is locally constant in $t$,
where $\Phi_t=F|_{X\times\{t\}}\in [(X,\underline{o});(X^\prime,
\underline{o^\prime}_t)]_{\eta_t}$ and $u_t\in [(I(G_{\hat{o}}),0,1);
(X^\prime,\underline{o^\prime}_0,\underline{o^\prime}_t)]_{(\eta_0,\eta_t)}$,
which is defined by the restriction of $F$ to $\{o\}\times [0,t]$.
Clearly, $(u_0)_\ast\circ (\Phi_0)_\ast=(\Phi_0)_\ast$ and
$(u_1)_\ast\circ (\Phi_1)_\ast=u_\ast\circ (\Phi_1)_\ast$.
Hence the claim.

Conversely, suppose $\Phi:X\rightarrow X^\prime$ is a weak
homotopy equivalence. We need to find a homotopy inverse
$\Psi:X^\prime\rightarrow X$ of $\Phi$.
First of all, since any orbispace in ${\cal{C}}$ is locally
path-connected and semi-locally 1-connected, the mapping cylinder $M_\phi$
of $\Phi:X\rightarrow X^\prime$ is defined. Hence by (1) above, $\Phi_\ast:
[[X^\prime;X]]\rightarrow [[X^\prime;X^\prime]]$ is a bijection.
In particular, there is a map $\Psi:X^\prime\rightarrow X$ such that
$\Phi\circ\Psi\cong Id_{X^\prime}\in [X^\prime;X^\prime]$. On the
other hand, note that $\Phi\circ\Psi\cong Id_{X^\prime}$ implies
$\Psi_\ast=\Phi_\ast^{-1}$ so that $\Psi$ is also a weak homotopy
equivalence. Thus there is a map $\Upsilon:X\rightarrow X^\prime$ such
that $\Psi\circ\Upsilon\cong Id_X$. But
$\Upsilon\cong \Phi\circ\Psi\circ\Upsilon\cong \Phi$. Hence $\Psi$ is a
homotopy inverse of $\Phi$, and $\Phi:X\rightarrow X^\prime$ is a
homotopy equivalence.

\hfill $\Box$

\subsection{CW-complex of groups and its geometric realization}

This final subsection is concerned with a subcategory of
${\cal{C}}$, where the objects have much more pleasant geometrical
properties. Most importantly, this subcategory is large enough that
all compact smooth orbifolds are contained in it. In a certain sense,
this subcategory of ${\cal{C}}$ consists of those orbispaces which are 
the geometric analogs of finite CW-complexes. We denote this subcategory 
of ${\cal{C}}$ by ${\cal{G}}$.

We give a description of ${\cal{G}}$ first. For the purpose here, we
need to impose an additional condition on the general CW-complexes, for
instance, as defined in \cite{Sw}, to which we refer the reader for the
basic definitions and properties of CW-complexes. The condition is: for
each attaching map, if its image meets the interior of a cell,
it contains the whole cell. Note that simplicial complexes
satisfy this condition. Immediate consequences of this assumption
include that every face of a cell is an immediate face, and that
each cell is a disjoint union of the interiors of finitely many
cells.

An open CW-complex is an open subset of a CW-complex which is a
disjoint union of the interiors of a subset of cells. Note that
under the additional assumption here, the closure of an open
CW-complex is also the smallest sub-complex containing the open
CW-complex. We will assume that every open CW-complex is
associated with a CW-complex which contains it as the interior,
and any map between open CW-complexes is the restriction of a map
between the closures. An open CW-complex is called an open
sub-complex of a CW-complex $K$ if the associated closure is the
closure in $K$.

For any cell $\sigma$, the star of $\sigma$, denoted by $\mbox{St}(\sigma)$,
is the smallest open sub-complex that contains the interior of $\sigma$. It is
also the disjoint union of the interiors of those cells which have $\sigma$
as a face. In the case of simplicial complex, this definition coincides with
the usual one. We remark that the star of a cell is connected, and
the underlying space of a CW-complex is locally connected.

Group actions on a CW-complex are required to satisfy the
following conditions: (1) the image of a cell is a cell, (2) if an
interior point of a cell is fixed, the whole cell must be fixed.
Let $G$ be a discrete group acting on a CW-complex $K$. Then $K/G$
is naturally a CW-complex, and the orbit map $\pi:K\rightarrow K/G$
is cellular.

\begin{defi}
The subcategory ${\cal{G}}$ consists of orbispaces $X$ where
\begin{itemize}
\item [{(1)}] $X$ is the underlying space of a finite CW-complex $K$,
\item [{(2)}] for each cell $\sigma\in K$, $\mbox{St}(\sigma)$ is a local
chart on $X$, such that in $(\widehat{\mbox{St}(\sigma)},G_{{St}(\sigma)},
\pi_{{St}(\sigma)})$, $\widehat{\mbox{St}(\sigma)}$ is an open CW-complex,
$\hat{\sigma}=\pi_{{St}(\sigma)}^{-1}(\sigma)$ is a cell in the associated
closure of $\widehat{\mbox{St}(\sigma)}$, and $\mbox{St}(\hat{\sigma})=
\widehat{\mbox{St}(\sigma)}$.
\end{itemize}
\end{defi}

For the sake of simplicity, we denote $G_{{St}(\sigma)}$ by $G_\sigma$
and $\pi_{{St}(\sigma)}$ by $\pi_\sigma$. We call $G_\sigma$ the
isotropy group of the cell $\sigma$. Note that if $X, X^\prime\in {\cal{G}}$,
the product $X\times X^\prime$ is also in ${\cal{G}}$.

The objects of ${\cal{G}}$ are closely related to the notion `CW-complex
of groups' described in Introduction. More precisely, to each $X\in
{\cal{G}}$, which is the underlying space of a finite CW-complex
$K$, one can associate an equivalence class of CW-complexes of
groups on $K$ as follows.

\begin{itemize}
\item [{(1)}] Each cell $\sigma\in K$ is associated with its
isotropy group $G_\sigma$.
\item [{(2)}] For each arrow $a$, note that $\mbox{St}(i(a))\subset
\mbox{St}(t(a))$ because $t(a)$ is a face of $i(a)$. Assign the arrow
$a$ with $\psi_a:G_{i(a)}\rightarrow G_{t(a)}$, where $\psi_a$ is
given by $\lambda_{\xi_a}$ for a fixed choice of $\xi_a\in
T(\mbox{St}(i(a)),\mbox{St}(t(a)))$.
\item [{(3)}] To each pair of composable arrows $a,b$, an element
$g_{a,b}\in G_{t(a)}$ is assigned, which is the unique element
satisfying $\xi_a\circ\xi_b=g_{a,b}\circ\xi_{ab}$. The equation
$$
Ad(g_{a,b})\circ\psi_{ab}=\psi_a\circ\psi_b \leqno(3.3.1)
$$
follows immediately from the definition. The cocycle condition for
a triple of composable arrows $a,b,c$
$$
\psi_a(g_{b,c})g_{a,bc}=g_{a,b}g_{ab,c} \leqno (3.3.2)
$$
is a consequence of the associativity of composition in the groupoid.
\end{itemize}

One can easily check that a different choice of $\{\xi_a\}$
will result in an equivalent CW-complex of groups on $K$. The
CW-complex of groups thus obtained is called associated to the
orbispace $X\in {\cal{G}}$, and the orbispace $X\in {\cal{G}}$ is
called the geometric realization of the associated CW-complex of groups.

We remark that the relationship between CW-complex of groups and
the geometric realization is analogous to that between complex of
groups and the associated orbihedron in Haefliger \cite{Ha3}. In
fact, when the underlying CW-complex is simplicial, a CW-complex
of groups is simply a complex of groups, and the corresponding
geometric realization is isomorphic as an orbispace to the associated
orbihedron. However, there is a minor difference between these two
concepts which lies in the fact that in the definition of an orbihedron
in Haefliger \cite{Ha3}, the star of a cell is not chosen to be
the one in the simplicial complex, but rather in the barycentric
subdivision of it. For details, see \cite{Ha3}.

The central result of this subsection is ${\cal{G}}\subset {\cal{C}}$.

Let $X\in {\cal{G}}$, with the finite CW-complex structure $K$,
and $\Phi:S^{k-1}(G)\rightarrow X$ be any map, such that the
induced map $\phi:S^{k-1}\rightarrow X$ is an attaching map so
that the mapping cone $C_\phi$ supports a canonical finite CW-complex
$L$ obtained by attaching a k-cell to $K$ via $\phi$. Recall
that there is a canonical orbispace structure on $S^{k-1}$ defined by
$\V(\Phi),\T(\Phi)$ as constructed in Lemma 3.1.2, which is equivalent to
the standard one on $S^{k-1}(G)$. In order to attach a k-cell of isotropy type
$G$ to $X$ via $\Phi$, we assume further that the orbispace structure
$(\V(\Phi),\T(\Phi))$ is contained in the standard one on $S^{k-1}(G)$,
cf. Proposition 3.2.3. In the present case, in order to ensure that the
orbispace $C_\phi$ constructed in Proposition 3.2.3 is an object of the
subcategory ${\cal{G}}$, we need to further impose two additional conditions:

\begin{itemize}
\item Let $(\{V_{\sigma,i}|i\in I_\sigma\},\{\mbox{St}(\sigma)\},
\{f_{\sigma,i}\},\{\rho_{(\tau,j)(\sigma,i)}\})$ be the canonical
representative of $\Phi$ constructed in Lemma 3.1.2, where $\sigma\in K$,
and $\phi^{-1}(\mbox{St}(\sigma))=\bigsqcup_{i\in I_\sigma} V_{\sigma,i}$.
(cf. $(3.1.5)$ also.) We further assume that each $\rho_{(\tau,j)(\sigma,i)}:
G\rightarrow T(\mbox{St}(\sigma),\mbox{St}(\tau))$ is independent of the
indexes $i\in I_\sigma, j\in I_\tau$.
\item For any cell $\sigma\subset \phi(S^{k-1})$, we require that the map
$\pi_\sigma:\widehat{\mbox{St}(\sigma)}\rightarrow \mbox{St}(\sigma)$,
which is defined over the associated closures of the open CW-complexes by
our earlier assumption, is one to one when restricted to the closure
of $\bigcup_{i\in I_\sigma} f_{\sigma,i}(\widehat{V_{\sigma,i}})$.
\end{itemize}

\begin{lem}
With the preceding understood, the orbispace $C_\phi$ constructed in
Proposition 3.2.3 is canonically an object of the subcategory ${\cal{G}}$.
\end{lem}

\pf
First of all, we recall the following notation: the open subset
$M_f\setminus \{[y,t]|t\in [0,t_0]\}$ of a mapping cylinder $M_f$
is denoted by $M_f(t_0)$.

Now observe that in Proposition 3.2.3 (1), a local chart on $C_\phi$
is either a local chart $W$ on the open subspace $C_\phi\setminus\{\ast\}$
where $\{\ast\}$ is the cone point, or the interior $D$ of the attached k-cell
of isotropy type $G$. If $W\cap D\neq\emptyset$, $W$ must be 
$M_{\phi_\sigma}(0)$
for some cell $\sigma\in K$, where $M_{\phi_\sigma}$ is the mapping cylinder
of $\phi_\sigma=\phi|_{\sqcup_{i\in I_\sigma}V_{\sigma,i}}$. Moreover,
$\widehat{W}=M_{f_\sigma}(0)$ where $f_\sigma$ is a map into
$\widehat{\mbox{St}(\sigma)}$ defined as follows. $\mbox{Domain
}(f_\sigma)=\bigsqcup_{i\in I_\sigma} (G_\sigma\times\widehat{V_{\sigma,i}})
/\rho_{(\sigma,i)}(G)$, where the action is given
by $\rho_{(\sigma,i)}(g)\cdot (g^\prime,x)=(g^\prime\rho_{(\sigma,i)}(g)^{-1},
g\cdot x),\; \forall g\in G$, and $f_\sigma$ is defined by $f_\sigma([(g^\prime,x)])
=g^\prime\cdot f_{\sigma,i}(x),\; \forall g^\prime\in G_\sigma, x\in
\widehat{V_{\sigma,i}}, i\in I_\sigma$. With this understood, the action of
$G_W=G_\sigma$ on $\widehat{W}$ is given by $h\cdot [(g^\prime,x)]=[(hg^\prime,x)]$.
Now by the first imposed condition, $\rho_{(\sigma,i)}=\rho_\sigma:G\rightarrow
G_\sigma$ is independent of the index $i\in I_\sigma$. This allows us to define
$D_\sigma=M_{\phi_\sigma}(0)\cup D$, $\widehat{D_\sigma}=M_{f_\sigma}(0)\cup
(G_\sigma\times D)/\rho_\sigma(G)$, and $G_{D_\sigma}=G_\sigma$ with
the natural action by multiplication from the left. Again the first imposed
condition that each $\rho_{(\tau,j)(\sigma,i)}:G\rightarrow T(\mbox{St}(\sigma),
\mbox{St}(\tau))$ is independent of the indexes $i\in I_\sigma, j\in I_\tau$
allows us to define naturally a set $\{T(D_\sigma,D_\tau)\}$, such that a set
of new local charts $\{(\widehat{D_\sigma},G_{D_\sigma},\pi_{D_\sigma})\}$ may
be added consistently to the canonical orbispace structure on $C_\phi$.

Next we prove that with the modified orbispace structure which is equivalent
to the original one, $C_\phi$ belongs to ${\cal{G}}$. Recall that the finite
CW-complex structure $L$ on $C_\phi$ is the one obtained by attaching the k-cell
$\overline{D}$ to $K$ via $\phi$. Thus the following is true for $L$: (1) For
any $\sigma\in L$ such that $\sigma\neq \overline{D}$ and $\sigma$ is not
a face of $\overline{D}$, the star of $\sigma$ in $L$ is $\mbox{St}(\sigma)$,
the star of $\sigma$ in $K$. (2) If $\sigma=\overline{D}$, then the star of
$\sigma$ in $L$ is the interior of $\sigma$. (3) If 
$\sigma\subset\phi(S^{k-1})$,
then the star of $\sigma$ in $L$ is $\mbox{St}(\sigma)\sqcup D=D_\sigma$.
It is clear that in order to show that $C_\phi$ belongs to ${\cal{G}}$,
it suffices to verify that $\widehat{D_\sigma}$, which is
$\widehat{\mbox{St}(\sigma)}\sqcup (G_\sigma/\rho_\sigma(G))\times D$, has a
natural open CW-complex structure, and that it is the star of $\hat{\sigma}=
\pi_\sigma^{-1}(\sigma)$ in the associated closure. To this end, we observe 
that
the closure of $D_\sigma$ in $L$ is $\overline{\mbox{St}(\sigma)}\sqcup
(\overline{D}\setminus \overline{\mbox{St}(\sigma)})$. The second imposed
condition then implies that
$$
(\mbox{closure of  }\widehat{{\mbox{St}(\sigma)}})
\sqcup (G_\sigma/\rho_\sigma(G))\times
(\overline{D}\setminus\overline{\mbox{St}(\sigma)}) \leqno (3.3.3)
$$
is naturally a CW-complex, which contains $\widehat{D_\sigma}$ as the
interior, and admits a natural action of $G_\sigma$ extending that of
$G_\sigma$ on the interior, such that the orbit space is the closure of
$D_\sigma$ in $L$. Hence $\widehat{D_\sigma}$ has a natural open CW-complex
structure. Finally, it is easily seen that $\widehat{D_\sigma}$ is the
star of $\hat{\sigma}=\pi_\sigma^{-1}(\sigma)$ in the associated closure
$(3.3.3)$.

\hfill $\Box$

Now we are ready for a proof of Proposition 1.7.

\begin{prop}
\begin{itemize}
\item [{(1)}] Let $(K,G_\sigma,\psi_a,g_{a,b})$ be a CW-complex of groups,
and $X$ be the underlying space of $K$. Then there is a canonical
orbispace structure on $X$ such that the orbispace $X$ belongs to ${\cal{G}}$,
and the associated CW-complex of groups is $(K,G_\sigma,\psi_a,g_{a,b})$.
Moreover, the orbispace $X$ also belongs to ${\cal{C}}$ with the canonical
filtration of skeletons
$$
X_0\subset X_1\subset\cdots\subset X_n=X,
$$
where $X_k$, $0\leq k\leq n$, is the underlying space of the k-skeleton of
$K$, such that the orbispace structure on $X_k$ is the one
canonically determined by the restriction of $(K,G_\sigma,\psi_a,g_{a,b})$
to the k-skeleton of $K$.
\item [{(2)}] The orbispace structure of an object in ${\cal{G}}$
is uniquely determined by the equivalence class of the associated
CW-complexes of groups.
\end{itemize}
\end{prop}

\pf
(1) The statement is true when the dimension of $K$ is zero. We shall
prove that if it is true when the dimension is $n$, it is also
true when the dimension is $n+1$.

First of all, by the induction assumption, there is a canonical orbispace
structure on $X_n$ such that the orbispace $X_n\in {\cal{G}}$ and
the associated CW-complex of groups is $(K_n,G_\sigma,\psi_a,g_{a,b})$.
In particular, the latter means that there are $\xi_a\in T(\mbox{St}(i(a)),
\mbox{St}(t(a)))$ satisfying $\xi_a\circ\xi_b=g_{a,b}\circ\xi_{ab}$
and $\lambda_{\xi_a}=\psi_a$. We also observe, from the proof of
Lemma 3.3.2, that $\pi_\sigma|_{\hat{\sigma}}$
is a homeomorphism onto $\sigma$ for any cell $\sigma$.

In order to construct an orbispace structure on $X_{n+1}$ such
that $X_{n+1}\in {\cal{G}}$, we shall define, for each (n+1)-cell
$e$, a map $\Phi_e:S^n(G_e)\rightarrow X_n$ whose induced map $\phi_e:
S^n\rightarrow X_n$ is the attaching map for the (n+1)-cell $e$,
such that the hypothesis in Lemma 3.3.2 holds for $\Phi_e$. To
this end, we set $\eta_a=g_{a,d}^{-1}\circ\xi_a$ for any arrow $a$
such that both $i(a), t(a)$ are contained in the image of the
attaching map $\phi_e$ of $e$, where $d$ is the arrow $(e,i(a))$.
It is easy to check that $\eta_a\circ\eta_b=\eta_{ab}$ holds for
any composable arrows $a,b$. Now for any cell $\sigma\subset
\mbox{Im }\phi_e$, we define a map $u_\sigma:\mbox{Im }\phi_e\cap
\mbox{St}(\sigma)\rightarrow\widehat{\mbox{St}(\sigma)}$ by
setting $u_\sigma|_{\sigma^0}=\pi_\sigma^{-1}|_{\sigma^0}$ and
$u_\sigma|_{\tau^0}=\phi_{\eta_a}\circ\pi_\tau^{-1}|_{\tau^0}$ for any
cell $\tau$ such that $\tau^0\subset \mbox{Im }\phi_e\cap\mbox{St}(\sigma)$
and $\tau\neq\sigma$, where $\sigma^0$, $\tau^0$ denote the
interior of the corresponding cell, and the arrow $a=(\tau,\sigma)$.
We note that (1) $\pi_\sigma\circ u_\sigma=Id$, (2) $\phi_{\eta_a}\circ
u_{i(a)}=u_{t(a)}|_{Dom(u_{i(a)})}$ for any arrow $a$ such that
both $i(a),t(a)$ are contained in $\mbox{Im }\phi_e$, and (3) the
image of $u_\sigma$ lies in the fixed-point set of $\psi_d(G_e)\subset
G_\sigma$ where $d$ is the arrow $(e,\sigma)$. Furthermore, it
follows from (2) above that each $u_\sigma$ is continuous.

For each cell $\sigma\subset\mbox{Im }\phi_e$, let $\{V_{\sigma,i}|
i\in I_\sigma\}$ be the set of connected components of $\phi_e^{-1}
(\mbox{St}(\sigma))$. We define $f_{\sigma,i}=u_\sigma\circ \phi_e|
_{V_{\sigma,i}}$ (with $\widehat{V_{\sigma,i}}=V_{\sigma,i}$ understood),
$\rho_{(\sigma,i)}=\psi_d$ where $d=(e,\sigma)$,
and for any $\sigma,\tau$ such that $\tau$ is a face of $\sigma$,
define $\rho_{(\tau,j)(\sigma,i)}$ by $g\mapsto\eta_a\circ\psi_d(g),\;
\forall g\in G_e$, where $a=(\sigma,\tau)$ and $d=(e,\sigma)$. One
can easily check, using the properties of $\{u_\sigma\}$ and $\{\eta_a\}$
established in the preceding paragraph, that $(\{f_{\sigma,i}\},
\{\rho_{(\tau,j)(\sigma,i)}\})$ is a homomorphism of groupoids.
The equivalence class of $(\{f_{\sigma,i}\},\{\rho_{(\tau,j)(\sigma,i)}\})$
is defined to be $\Phi_e$, which satisfies the hypothesis in Lemma
3.3.2 by the nature of construction. By Lemma 3.3.2, there is a
canonical orbispace structure on $X_{n+1}$ such that $X_{n+1}\in
{\cal{G}}$.

To see that the associated CW-complex of groups is $(K_{n+1},G_\sigma,
\psi_a,g_{a,b})$, we need to find a $\xi_d\in T(\mbox{St}(e),
\mbox{St}(\sigma))$ for each arrow $d=(e,\sigma)$, where $e$ is a
(n+1)-cell, such that $\lambda_{\xi_d}=\psi_d$, and for any arrow
$a$ satisfying $i(a)=t(d)$, $\xi_a\circ\xi_d=g_{a,d}\circ\xi_{ad}$.
Observe that, according to Lemma 3.3.2, the inverse image of $e^0$
(the interior of $e$) in $\widehat{\mbox{St}(\sigma)}$ is
$(G_\sigma/\psi_d(G_e))\times e^0$, and $T(\mbox{St}(e),\mbox{St}(\sigma))
=G_\sigma$ with $\lambda_\xi=Ad(\xi)\circ\psi_d:G_e\rightarrow
G_\sigma$. With this understood, $\xi_d=1\in G_\sigma$ will work
for the purpose here. This completes the induction step, and hence
(1) of the proposition.

\vspace{1.5mm}

(2) The proof goes by induction on the dimension of $K$. Assume the
statement is true when the dimension is $n$. Then it is easy to see
that the induction step boils down to the verification that for each
(n+1)-cell $e$ and any arrow $d=(e,\sigma)$, the attaching maps of
$\pi_\sigma^{-1}(e)$ in $\widehat{\mbox{St}(\sigma)}$ are canonically
determined by the CW-complex of groups. For this we simply observe: When
$d=(e,\sigma)$ is primitive, i.e., there is no $\tau$ such that $d$ is the
composition of $(e,\tau)$ with $(\tau,\sigma)$, the attaching maps are unique.
When $d=ad^\prime$ where $d^\prime$ is primitive, the attaching maps are
determined by $g_{a,d^\prime}^{-1}\circ\xi_a$ from those of
$\pi_{t(d^\prime)}^{-1}(e)$ in $\widehat{\mbox{St}(t(d^\prime))}$.

\hfill $\Box$

The preceding proposition clearly established the one to one
correspondence between equivalence classes of CW-complexes
of groups and isomorphism classes of orbispaces in ${\cal{G}}$,
and in particular, the inclusion ${\cal{G}}\subset {\cal{C}}$.

\vspace{2mm}

{\Small Current Address: Department of Mathematics and Statistics, 
University of Massachusetts at Amherst, Amherst, MA 01003,
{\it e-mail:} wchen@@math.umass.edu}


\begin{thebibliography}{}

\bibitem {Bre} G.E. Bredon, {\em Equivariant cohomology theories},
                            Lecture Notes in Math. {\bf 34} (1967),
                            Springer-Verlag.

\bibitem {BH} M.R. Bridson and A. Haefliger,
                       {\em Metric spaces of non-positive curvature},
                       Grundlehren der mathematischen Wissenschaften
                       {\bf 319}, Springer-Verlag, 1999.


\bibitem {C1} W. Chen,
      {\em On a notion of maps between orbifolds, I. function spaces},
             Communications in Contemporary
             Mathematics {\bf 8} no.5 (2006), 569-620.



\bibitem {C2} W. Chen, {\em Orbifold adjunction formula and symplectic
                            cobordisms between lens spaces}, Geometry and
                        Topology {\bf 8} (2004), 701-734.
                            

\bibitem {C3} W. Chen, {\em Smooth $s$-cobordisms of elliptic 
                            $3$-manifolds}, J. Diff. Geom. {\bf 73} no. 3
                            (2006), 413-490.

\bibitem {C4} W. Chen, {\em Pseudoholomorphic curves in four-orbifolds and some
                   applications}, in Geometry and Topology of
                   Manifolds, Boden, H.U. et al ed., Fields Institute 
                   Communications {\bf 47}, pp. 11-37. Amer. Math. Soc.,
                                  Providence, RI, 2005


\bibitem {CR1} W. Chen and Y. Ruan, {\em Orbifold Gromov-Witten theory},
                                  in ``Orbifolds in Mathematics
                                  and Physics'', 25-85, Edited by Adem et
                                  al., Contemporary Mathematics
                                  {\bf 310}, Amer. Math. Soc.,
                                  Providence, RI, 2002.

\bibitem {CR2} W. Chen and Y. Ruan, {\em A new cohomology theory of 
                               orbifolds}, Communications in Mathematical 
                                Physics {\bf 248} no.1 (2004), 1-31.


\bibitem {DHVW} L. Dixon, J. Harvey, C. Vafa and E. Witten, {\em
                                  Strings on orbifolds, I,II},
                                  Nucl. Phys. B261 (1985) 678-,
                                  B274 (1986) 285-.


\bibitem {Ha1} A. Haefliger, {\em Homotopy and integrability},
                    in ``Manifolds -- Amsterdam 1970'', Lecture
                    Notes in Math. {\bf 197}, Springer-Verlag.

\bibitem {Ha2} A. Haefliger, {\em Groupoides d'holonomie et
                                  classifiants}, Ast\'{e}risque
                                  {\bf 116} (1984), 70-97.


\bibitem {Ha3} A. Haefliger, {\em Complexes of groups and
                                  orbihedra}, in ``Group theory
                                  from a geometrical viewpoint, 26
                                  March --- 6 April 1990, ICTP,
                                  Trieste'', World Scientific (1991),
                                  504-540.

\bibitem {Pr} D. A. Pronk, {\em Groupoid representations for
                             sheaves on orbifolds}, Ph.D. Thesis,
                             Utrecht 1995.

\bibitem {Sw}  R. M. Switzer,
               {\em Algebraic topology -- homotopy and homology},
               Springer-Verlag, 1975.



\bibitem {Th} W. Thurston, {\em The geometry and topology of
                               three-manifolds}, Princeton lecture
                               notes.

\bibitem {Y} C. T. Yang, {\em The triangulability of the orbit
                          space of a differentiable transformation
                          group}, Bull. Amer. Math. Soc.
                          {\bf 69} (1963), 405-408.

\end{thebibliography}
\end{document}